\documentclass[a4paper,12pt]{article}
\usepackage[cp1251]{inputenc}
\usepackage{amssymb,amsmath,amsthm}
\usepackage[russian]{babel}
\usepackage{indentfirst}
\usepackage{amssymb,amsmath,amsthm}
\usepackage{hyperref,hypbmsec}
\usepackage[pdftex]{graphicx}

\numberwithin{equation}{section}

\newtheorem{Th}{\hskip\parindent Theorem}[section]
\newtheorem{Le}{\hskip\parindent Lemma}[section]
\newtheorem{Sl}{\hskip\parindent Corollary}[section]
\newtheorem{Zam}{\hskip\parindent Remark}[section]
\newtheorem{Hyp}{\hskip\parindent Conjecture}[section]

\newcommand{\A}{\mathcal{A}}
\newcommand{\E}{\mathfrak{C}}
\newcommand{\R}{\mathfrak{R}}
\newcommand{\s}{\mathbb{S}}
\newcommand{\D}{\mathfrak{D}}
\newcommand{\N}{\mathbb{N}}
\newcommand{\M}{\mathfrak{M}}
\newcommand{\NN}{\mathfrak{N}}
\newcommand{\1}{\mathbf{1}}
\newcommand{\jj}{\hat{j}}
\newcommand{\h}{\hat{h}}
\newcommand{\Om}{\widetilde{\Omega}}
\newcommand{\q}{\mathbf{q}}
\newcommand{\KK}{\overline{K}}
\newcommand{\K}{|K|}
\newcommand{\Z}{\mathbb{Z}}
\newcommand{\rr}{\mathbb{R}}
\newcounter{propet}

\renewcommand{\le}{\leqslant}\renewcommand{\ge}{\geqslant}
\renewcommand{\partname}{Chapter}

\begin{document}
\author{ D.\,A.\,Frolenkov\footnote{Research is supported by RFFI (grant № 11-01-00759-а)} \quad
I.\,D.\,Kan,\footnote{Research is supported by RFFI  (grant № 12-01-00681-а)} }
\title{
\begin{flushright}
\small{Dedicated to the memory of\\professor N.M. Korobov.}
\end{flushright}
A reinforcement of the Bourgain-Kontorovich's theorem.
}
\date{}
\maketitle
\begin{abstract}
Zaremba's conjecture (1971) states that every positive integer number $d$ can be represented as a denominator (continuant) of a finite continued fraction  $\frac{b}{d}=[d_1,d_2,\ldots,d_{k}],$ with all partial quotients $d_1,d_2,\ldots,d_{k}$ being bounded by an absolute constant $A.$ Recently (in 2011) several new theorems concerning this conjecture were proved by Bourgain and Kontorovich. The easiest of them states that the set of numbers satisfying Zaremba's conjecture with $A=50$ has positive proportion in $\N.$ In this paper the same theorem is proved with $A=7.$

Bibliography: 15 titles.

 \par
\end{abstract}

\renewcommand{\proofname}{{\bf Proof}}
\renewcommand{\partname}{Chapter}
\setcounter{Zam}0
\tableofcontents
\part{\large{Introduction}}
\section{Historical background}
Let $\R_{A}$ be the set of rational numbers whose continued fraction expansion has all partial quotients being bounded by  $A:$
$$
\R_{A}=\left\{\frac{b}{d}=[d_1,d_2,\ldots,d_{k}]\Bigl| 1\le d_{j}\le A \,\mbox{для}\, j=1,\ldots,k\right\},
$$
where
\begin{equation}\label{cont.fraction}
[d_1,\ldots,d_k]=\cfrac{1}{d_1+{\atop\ddots\,\displaystyle{+\cfrac{1}{d_k}}}}
\end{equation}
is a finite continued fraction, $d_1,\ldots,d_{k}$ are partial quotients. Let $\D_{A}$ be the set of denominators of numbers in $\R_{A}$:
$$
\D_{A}=\left\{d\Bigl| \exists b:\, (b,d)=1, \frac{b}{d}\in\R_{A}\right\}.
$$
And set
$$
\D_{A}(N)=\left\{d\in\D_{A}\Bigl| d\le N \right\}.
$$
\begin{Hyp}(\textbf{Zaremba's conjecture} ~\cite[p. 76]{Zaremba},\,1971 ).
For sufficiently large $A$  one has
$$\D_{A}=\N.$$
\end{Hyp}
That is, every $d\ge1$ can be represented as a denominator of a finite continued fraction $\frac{b}{d}$ whose partial quotients are bounded by $A.$ In fact Zaremba conjectured that $A=5$ is already large enough. A bit earlier, in the 1950-th, studying problems concerning  numerical integration Bahvalov, Chensov and N.M. Korobov also made the same assumption. But they did not publish it anywhere. Korobov ~\cite{Korobov2} proved that for a prime $p$ there exists $a,$ such that the  greatest partial quotient of $\frac{a}{p}$ is smaller than $\log p.$ A detailed survey on results concerning Zaremba's conjecture can be found in ~\cite{BK},~\cite{NG}.\par
Bourgain and Kontorovich suggested that the problem should be generalized in the following way. Let $\A \in \N$ be any finite alphabet ($|\A|\ge2$) and let $\R_{\A}$ and $\E_{\A}$  be the set of finite and infinite continued fractions
whose partial quotients belong to $\A:$
$$\R_{\A}=\left\{[d_1,\ldots,d_k]: d_j\in\A,\,j=1,\ldots,k\right\},$$
$$\E_{\A}=\left\{[d_1,\ldots,d_j,\ldots]: d_j\in\A,\,j=1,\ldots\right\}.$$
And let
$$\D_{\A}(N)=\left\{d\Bigl| d\le N,\, \exists b: (b,d)=1,\, \frac{b}{d}\in\R_{\A}\right\}$$
be the set of denominators bounded by $N.$ Let $\delta_{\A}$ be the Hausdorff dimension of $\E_{\A}.$
Then the Bourgain-Kontorovich's theorem ~\cite[p. 13, Theorem 1.25]{BK} ia as follows
\begin{Th}\label{BKtheorem}
For any alphabet $\A$  with
\begin{equation}\label{BKcondition}
\delta_{\A}>1-\frac{5}{312}=0,983914\ldots,
\end{equation}
the following inequality (positive proportion)
\begin{equation}\label{BKresult}
\#\D_{\A}(N)\gg N.
\end{equation}
holds.
\end{Th}
For some alphabets the condition \eqref{BKcondition} can be verified by two means. For an alphabet $\A=\left\{1,2,\ldots,A-1,A\right\}$ Hensley ~\cite{Hen1} proved that
\begin{equation}\label{hen-hausdorff-dimension}
\delta_{\A}=\delta_{A}=1-\frac{6}{\pi^2}\frac{1}{A}-\frac{72}{\pi^4}\frac{\log A}{A^2}+O\left(\frac{1}{A^2}\right).
\end{equation}
Moreover Jenkinson ~\cite{Jenkinson} obtained approximate values for some $\delta_{\A}.$ In view of these results the alphabet
\begin{equation}\label{specalpha}
\left\{1,2,\ldots,A-1,A\right\}
\end{equation}
with $A=50$ is assumed to satisfy \eqref{BKcondition}. Several results improving \eqref{BKresult} were also proved in ~\cite{BK}.
However, we do not consider them in our work.

\section{Statement of the main result}
First of all we must state one lemma ~\cite[p. 46, Lemma 7.1.]{BK}. The proof of our main result is essentially based on this lemma.\par
To begin with we describe all necessary objects. Let $K, X, Y \ge1$ be real numbers and $q$ be a positive integer.
Moreover, let $\eta=(x,y)^t, \eta'=(u,v)^t \in \mathbb{Z}^2$ be vectors such as
$$|\eta|\asymp\frac{X}{Y}, |\eta'|\asymp X, (x,y)=1,(u,v)=1.$$
\begin{Le}\label{BK lemma 7.1}(~\cite[p. 46, Lemma 7.1.]{BK})
If the following inequality
\begin{equation}\label{BK lemma 7.1 condition}
(qK)^{\frac{13}{5}}<Y<X,
\end{equation}
holds, then
\begin{equation}\label{BK lemma 7.1 statement}
\#\left\{\gamma\in SL_2(\mathbb{Z})\Bigl| \|\gamma\|\asymp Y, |\gamma\eta-\eta'|<\frac{X}{K}, \gamma\eta\equiv\eta'\pmod{q}\right\}\ll\frac{Y^2}{(qK)^2}.
\end{equation}
\end{Le}
The main result of the paper is the following theorem.
\begin{Th}\label{uslov}
For any alphabet $\A$ with
\begin{equation}\label{KFcondition1}
\delta_{\A}>1-\frac{27-\sqrt{633}}{16}=0,8849\ldots,
\end{equation}
the following inequality (positive proportion) holds
\begin{equation}\label{KFresult}
\#\D_{\A}(N)\gg N.
\end{equation}
\end{Th}
\begin{Zam}
It is proved ~\cite{Jenkinson} that $\delta_{7}=0,8889\ldots.$ From this follows that the alphabet $\left\{1,2,\ldots,7\right\}$
satisfies the condition of Theorem \ref{uslov}. It is also proved in~\cite{Jenkinson} that for the alphabet $\A=\left\{1,2,\ldots,6,8\right\}$ one has $\delta_{\A}=0,8851\ldots.$ Consequently, the alphabet $\A=\left\{1,2,\ldots,6,8\right\}$ also satisfies the condition of Theorem \ref{uslov}.
\end{Zam}
\begin{Zam}\label{zam1}
It seems to follow from the proof of the Lemma \ref{BK lemma 7.1} that the condition \eqref{BK lemma 7.1 condition} can be replaced by a weaker one
\begin{equation}\label{BK lemma 7.1 condition-2}
(qK)^{\frac{64}{25}+\epsilon}<Y<X.
\end{equation}
Then the statement of the Theorem \ref{uslov} holds with
\begin{equation}\label{KFcondition1-2}
\delta_{\A}>1-\frac{25}{114+2\sqrt{2274}}=0,8805\ldots.
\end{equation}
\end{Zam}
\begin{Zam}
The proof of Lemma \ref{BK lemma 7.1} given in the article~\cite{BK} is based on the paper~\cite{BKS}. In this article  using the spectral theory of automorphic forms statements similar to Lemma \ref{BK lemma 7.1} are proved. In our next article we are going to prove a result similar to Lemma \ref{BK lemma 7.1} using only estimates for Kloosterman sums. Certainly we are able to prove only a weaker result, however, it will suffice to prove that for the alphabet $\left\{1,2,\ldots,12,13\right\}$ the inequality \eqref{KFresult} holds.
\end{Zam}
\begin{Zam}
The proof of the Theorem \ref{uslov} is  based on the method of Bourgain-Kontorovich~\cite{BK}. In this article a new technique  for estimating exponential sums taking over a set of continuants is devised. We improve on this method by refining the main set $\Omega_N.$ In ~\cite{BK} this set is named as ensemble.
\end{Zam}
\section{Acknowledgment}
We thank prof. N.G. Moshchevitin for numerous discussions of our results. It was he who took our notice of ~\cite{BK}. We are also grateful to prof. I.D. Shkredov and prof. I.S. Rezvyakova for their questions and comments during our talks.
\section{Notation}
Throughout $\epsilon_0=\epsilon_0(\A)\in(0,\frac{1}{2500}).$ For two functions $f(x), g(x)$ the Vinogradov notation $f(x)\ll g(x)$  means that there exists a constant $C,$ depending on $A,\epsilon_0,$ such that $|f(x)|\le Cg(x).$  The notation $f(x)=O(g(x))$ means the same. The notation $f(x)=O_1(g(x))$ means that $|f(x)|\le g(x).$ The notation $f(x)\asymp g(x)$ means that
$f(x)\ll g(x)\ll f(x).$ Also a traditional notation $e(x)=\exp(2\pi ix)$ is used. The cardinality of a finite set $S$ is denoted either $|S|$ or $\#S.$ $[\alpha]$ and $\|\alpha\|$ denote the integral part of $\alpha$ and the distance from $\alpha$ to the nearest integer:
$[\alpha]=\max\left\{z\in\mathbb{Z}|\,z\le \alpha \right\},$
$\|\alpha\|=\min\left\{|z-\alpha|\Bigl|\,z\in\mathbb{Z}\right\}. $
The following sum $\sum_{d|q}1$ is denoted as $\tau(q).$
\part{\large{Preparation for estimating exponential sums}}
\section{Continuants and matrices}
In this section we recall the simplest techniques concerning continuants. As a rule, all of them can be found in any research dealing with continued fractions. To begin with we define several operations on finite sequences. Let
\begin{equation}\label{def of seq}
D=\{d_1,d_2,\ldots,d_k\},\, B=\{b_1,b_2,\ldots,b_n\}
\end{equation}
then $D,B$ denotes the following sequence
$$D,B=\{d_1,d_2,\ldots,d_k,b_1,b_2,\ldots,b_n\}.$$
For every $D$ from \eqref{def of seq} let define $D_{-},D^{-},\overleftarrow{D}$as follows
\begin{equation*}
D_{-}=\{d_2,d_3,\ldots,d_k\},\,
D^{-}=\{d_1,d_2,\ldots,d_{k-1}\},\,
\overleftarrow{D}=\{d_k,d_{k-1},\ldots,d_2,d_1\}.
\end{equation*}
We denote by $[D]$ the continued fraction \eqref{cont.fraction}, that is $[D]=[d_1,\ldots,d_k].$ And by $\langle D\rangle$ we denote its denominator $\langle D\rangle=\langle d_1,\ldots,d_k\rangle.$ This denominator is called the continuant of the sequence $D.$ The continuant of the sequence can also be defined as follows
\begin{equation*}
\langle\,\rangle=1,\,\langle d_1\rangle=d_1,
\end{equation*}
\begin{equation*}
\langle d_1,\ldots,d_k\rangle=\langle d_1,\ldots,d_{k-1}\rangle d_k+\langle d_1,\ldots,d_{k-2}\rangle,\,\mbox{для}\, k\ge2.
\end{equation*}
It is well known ~\cite{knut} that
\begin{equation}\label{continuant properties}
\langle D\rangle=\langle\overleftarrow{D}\rangle,\,
[D]=\frac{\langle D_{-}\rangle}{\langle D\rangle},\,
[\overleftarrow{D}]=\frac{\langle D^{-}\rangle}{\langle D\rangle},
\end{equation}
\begin{equation}\label{continuant identity}
\langle D,B\rangle=\langle D\rangle\langle B\rangle(1+[\overleftarrow{D}][B]).
\end{equation}
It follows from this that
\begin{equation}\label{continuant inequality}
\langle D\rangle\langle B\rangle\le
\langle D,B\rangle\le
2\langle D\rangle\langle B\rangle,
\end{equation}
and that the elements of the matrix
\begin{equation}\label{matrix}
\gamma=
\begin{pmatrix}
a & b \\
c & d
\end{pmatrix}=
\begin{pmatrix}
0 & 1 \\
1 & d_1
\end{pmatrix}
\begin{pmatrix}
0 & 1 \\
1 & d_2
\end{pmatrix}\ldots
\begin{pmatrix}
0 & 1 \\
1 & d_k
\end{pmatrix}
\end{equation}
can be expressed by continuants
\begin{gather}\label{matrix=continuant}
a=\langle d_2,d_3,\ldots,d_{k-1}\rangle,\,
b=\langle d_2,d_3,\ldots,d_{k}\rangle,\,
c=\langle d_1,d_2,\ldots,d_{k-1}\rangle,\,
d=\langle d_1,d_2,\ldots,d_{k}\rangle.
\end{gather}
For the matrix $\gamma$ from \eqref{matrix} we use the following norm
\begin{equation*}
\|\gamma\|=\max\{|a|,|b|,|c|,|d|\},
\end{equation*}
It follows from \eqref{matrix=continuant} that
\begin{equation}\label{matrix norm}
\|\gamma\|=d=\langle d_1,d_2,\ldots,d_{k}\rangle.
\end{equation}
For $\gamma$ from \eqref{matrix} we have $\det \gamma=(-1)^k.$ So for an even $k$ one has $\det \gamma=1,$ that is $\gamma\in SL\left(2,\mathbb{Z}\right).$\par
Let $\Gamma_{\A}\subseteq SL\left(2,\mathbb{Z}\right)$ be a semigroup generated by
\begin{equation*}
\begin{pmatrix}
0 & 1 \\
1 & a_i
\end{pmatrix}
\begin{pmatrix}
0 & 1 \\
1 & a_j
\end{pmatrix}=
\begin{pmatrix}
1   & a_j \\
a_i & a_ia_j+1
\end{pmatrix},
\end{equation*}
where $a_i,a_j \in\A.$ It follows from \eqref{matrix norm} that to prove a positive density of continuants in $\N$
\begin{equation}\label{BKresult1,1}
\#\D_{\A}(N)\gg N,
\end{equation}
it is enough to obtain the same property of the set
\begin{equation}\label{set matrix norm}
\|\Gamma_{\A}\|=\left\{\|\gamma\|\Bigl| \gamma\in\Gamma_{\A}\right\}.
\end{equation}
In fact, for proving inequality  \eqref{BKresult1,1} only a part of the semigroup $\Gamma_{\A},$ a so called ensemble $\Omega_N,$ will be used. A preparation for constructing $\Omega_N$ will start in the next section.
\section{Hensley's method}\label{method of hensley}
Before estimating the amount of continuants not exceeding $N$ it might be well to assess the amount of continued fractions with denominator being bounded by $N.$ Though these problems are similar, in the second case every continuant should be counted in view of its multiplicity. Let $V_{\A}(k)$ be the set of words with the length $k$
$$
V_{\A}(k)=\left\{(d_1,d_2,\ldots,d_{k})\Bigl| 1\le d_{j}\le A ,\, j=1,\ldots,k\right\},
$$
and let $V_{\A}=\bigcup_{k\ge1}V_{\A}(k)$ be the set of all finite words. Let
\begin{gather}\label{maxA}
A=\max\A.
\end{gather}
We will consider on alphabet $\A$ only words having an even length. They can also be treated as words on the alphabet $(\A,\A),$ that is, consisting of  pairs $(a,b)$ where $a,b\in\A.$ Let us denote the alphabet $(\A,\A)$ by $\A^2.$ Let $V_{\A^2}$ be the set of words on $\A$ having an even length. Let $\R_{\A^2}$ be the set of finite continued fractions constructed from sequences in  $V_{\A^2}.$ We also write
$$\R_{\A^2}(N)=\left\{D\in V_{\A^2}\Bigl|\langle D\rangle\le N\right\},$$
$$
F_{\A}(x)=\#\left\{D\in V_{\A^2}\Bigl|\langle D\rangle\le x\right\}=\#\R_{\A^2}(x).
$$
And let $D_{\A^2}(N)$ be the set of denominators of fractions from $\R_{\A^2}(N).$ Generalizing Hensley's method~\cite{Hen2} one can prove that
$2\delta_{\A^2}=2\delta_{\A}$ and that the following theorem holds.
\begin{Th}\label{Kan-Hen}
Let $\delta_{\A}>\frac{1}{2},$ then for any $x\ge4A^2$ one has
\begin{gather}\label{F_A less}
\frac{1}{32A^4}x^{2\delta_{\A}}\le F_{\A}(x)-F_{\A}\left(\frac{x}{4A^2}\right)\le F_{\A}(x)\le8 x^{2\delta_{\A}}.
\end{gather}
\end{Th}
Hensley~\cite{Hen2} proved this theorem for the alphabet $\A$ of the form \eqref{specalpha}.
\section{The basic ideas of the Bourgain-Kontorovich's method}\label{section idea BK}
In this section a notion of constructions necessary for proving Theorem  \ref{uslov} will be given.\par
In view of exponential sums, for studying the density of the set \eqref{set matrix norm} it is natural to estimate the absolute value of the sum
\begin{gather}\label{def Sn 1}
S_N(\theta):=\sum_{\gamma\in\Omega_{N}}e(\theta\|\gamma\|)
\end{gather}
where $\Omega_{N}\subseteq\Gamma_{\A}\cap\{\|\gamma\|\le N\}$ is a proper set of matrices (ensemble), $\theta\in[0,1],$ and the norm $\|\gamma\|$ is defined in \eqref{matrix norm}. As usual, the Fourier coefficient of the function $S_N(\theta)$ is defined by
\begin{gather*}
\widehat{S}_N(n)=\int_0^1 S_N(\theta)e(-n\theta)d\theta=\sum_{\gamma\in \Omega_{N}}\1_{\{\|\gamma\|=n\}}.
\end{gather*}
Note that if ${\widehat{S}_N(n)>0}$ then $n\in\D_{\A^2}(N).$ Since
\begin{gather*}
S_N(0)=\sum_{\gamma\in\Omega_{N}}1=\sum_{n=1}^N \sum_{\gamma\in \Omega_{N}}\1_{\{\|\gamma\|=n\}}=
\sum_{n=1}^N \widehat{S}_N(n)=\sum_{n=1}^N \widehat{S}_N(n)\1_{\{\widehat{S}_N(n)\neq0\}},
\end{gather*}
then applying Cauchy-Schwarz inequality one has
\begin{gather}\label{8-4}
(S_N(0))^2\le \sum_{n=1}^N \1_{\{\widehat{S}_N(n)\neq0\}}\sum_{m=1}^N \left(\widehat{S}_N(m)\right)^2.
\end{gather}
The first factor of the right hand side of the inequality \eqref{8-4} is bounded from above by $\#\D_{\A^2}(N).$ Applying Parseval for the second factor one has
\begin{gather*}
\sum_{n=1}^N \left(\widehat{S}_N(n)\right)^2=\int_0^1\left|S_N(\theta)\right|^2d\theta.
\end{gather*}
Consequently
\begin{gather}\label{8-5}
(S_N(0))^2\le\#\D_{\A^2}(N)\int_0^1\left|S_N(\theta)\right|^2d\theta.
\end{gather}
Thus a lower bound on the magnitude of the set $\D_{\A^2}(N)$ follows from \eqref{8-5}
\begin{gather}\label{8-6}
\#\D_{\A^2}(N)\ge\frac{(S_N(0))^2}{\int\limits_0^1\left|S_N(\theta)\right|^2d\theta}.
\end{gather}
Thus, the estimate
\begin{gather*}
\#\D_{\A}(N)\gg N
\end{gather*}
will be proved, if we are able to assess exactly as possible the integral from \eqref{8-6}
\begin{gather}\label{8-7}
\int_0^1\left|S_N(\theta)\right|^2d\theta\ll\frac{(S_N(0))^2}{N}=\frac{|\Omega_{N}|^2}{N}.
\end{gather}
It follows from the Dirichlet's theorem that for any $\theta\in[0,1]$ there exist $a,q\in\N\cup\{0\}$ and $\beta\in\rr$ such that
\begin{gather*}
\theta=\frac{a}{q}+\beta,\;(a,q)=1,\; 0\le a\le q\le N^{1/2},\;\beta=\frac{K}{N},\; |K|\le\frac{N^{1/2}}{q},
\end{gather*}
with $a=0$ and $a=q$ being possible if only $q=1.$
Following the article ~\cite{BK}, to obtain the estimate \eqref{8-7} we represent the integral as the sum of integrals over different domains in variables $(q,K).$  Each of them will be estimated in a special way depending on the domain.\par
It remains to define ensemble $\Omega_{N}.$  To begin with we determine a concept "pre-ensemble". The subset $\Xi$ of matrices $\gamma\in\Gamma_{\A}$ is referred to as $N$~--\textbf{pre-ensemble}, if the following conditions hold
\begin{enumerate}
  \item for any matrix $\gamma\in\Xi$ its norm is of the order of $N:$
\begin{gather}\label{8-18a}
\|\gamma\|\asymp N;
\end{gather}
  \item for any $\epsilon>0$ the set $\Xi$ contains $\epsilon$~--full amount of elements, that is
\begin{gather}\label{8-18b}
\#\Xi\gg_{\epsilon} N^{2\delta_{\A}-\epsilon}.
\end{gather}
\end{enumerate}
By the \textbf{product} of two pre-ensembles $\Xi^{(1)}\Xi^{(2)}$ we mean the set of all possible  products of matrices $\gamma_1\gamma_2$ such that $\gamma_1\in\Xi^{(1)},\,\gamma_2\in\Xi^{(2)}.$ The product of pre-ensembles has an \textbf{unique} \textbf{expansion} if it follows from the relations
\begin{gather}\label{8-19}
\gamma_1\gamma_2=\gamma_1'\gamma_2', \quad \gamma_1,\gamma_1'\in\Xi^{(1)},\,\gamma_2,\gamma_2'\in\Xi^{(2)}
\end{gather}
that
\begin{gather}\label{8-20}
\gamma_1=\gamma_1',\quad \gamma_2=\gamma_2'.
\end{gather}
Let $\epsilon_0$ be a fixed number such that $0<\epsilon_0\le\frac{1}{2}.$ Then $N$~-- pre-ensemble $\Omega$ is called the right (left) $(\epsilon_0,N)$~-- \textbf{ensemble} if for any $M,$ such that
\begin{gather}\label{8-21}
1\ll M\le N^{\frac{1}{2}},
\end{gather}
there exist positive numbers $N_1$ and $N_2$ such that
\begin{gather}\label{8-22}
N_1N_2\asymp N,\quad N_2^{1-\epsilon_0}\ll M\ll N_2,
\end{gather}
and there exist $N_1$~--pre-ensemble $\Xi^{(1)}$ and $N_2$~--pre-ensemble $\Xi^{(2)}$ such that the pre-ensemble $\Omega$ is equal to the product $\Xi^{(1)}\Xi^{(2)}$ ($\Xi^{(2)}\Xi^{(1)}$ respectively) having an unique expansion. Such terminology allows us to say that in the article ~\cite{BK} the $(\frac{1}{2},N)$~--ensemble has been constructed while we will construct $(\epsilon_0,N)$~--ensemble, being simultaneously the right and the left (bilateral) ensemble, for $\epsilon_0\in\left(0,\frac{1}{2500}\right).$
\begin{Zam}\label{zam813}
Observe that there is no use to require an upper bound in \eqref{8-18b}. According to the Theorem \ref{Kan-Hen} it follows from \eqref{8-18a} that
\begin{gather*}
\#\Xi\ll N^{2\delta_{\A}}.
\end{gather*}
\end{Zam}
\section{Pre-ensemble $\Xi(M).$}\label{predansambl}
Let $\delta:=\delta_{\A}>\frac{1}{2},$ $\Gamma:=\Gamma_{\A}$ and as usual $A=\max\A.$
Let also $M$ be a fixed parameter satisfying the inequality
\begin{gather}\label{9-1}
M\ge2^9A^3\log^3M.
\end{gather}
In this section we construct a pre-ensemble $\Xi(M)\subset\Gamma$ It is the key element which will be used to construct the ensemble $\Omega_N.$ To generate $\Xi(M)$ we use an algorithm. The number $M$ is an input parameter. During the algorithm we generate the following numbers
\begin{gather*}
L=L(M)\asymp M,\,p=p(M)\asymp\log\log M,\, k=k(M)\asymp\log M,
\end{gather*}
being responsible for the properties of the elements of $\Xi(M).$ We now proceed to the description of the algorithm consisting of four steps.
\begin{description}
  \item[Step 1]
First consider the set $S_1\subset\Gamma$ of matrices $\gamma\in\Gamma,$ such that
\begin{gather}\label{9-2}
\frac{M}{64A^2}\le\|\gamma\|\le M.
\end{gather}
According to the Theorem \ref{Kan-Hen}\, $\#S_1\asymp M^{2\delta}.$
  \item[Step 2]
Let $F_0=0,\,F_1=1,\,F_{n+1}=F_{n}+F_{n-1}$ for $n\ge1$ be Fibonacci numbers. числа Фибоначи. Define an integer number $p=p(M)$ by the relation
\begin{gather}\label{9-4}
F_{p-1}\le\log^{\frac{1}{2}}M\le F_p.
\end{gather}
Note that then
\begin{gather}\label{9-5}
F_{p}\le2F_{p-1}\le2\log^{\frac{1}{2}}M.
\end{gather}
Let consider the set $S_2\subset S_1$ of matrices $\gamma\in S_1$ of the form \eqref{matrix} for witch
\begin{gather}\label{9-6}
d_1=d_2=\ldots=d_p=1,\quad
d_k=d_{k-1}=\ldots=d_{k-p+1}=1.
\end{gather}
To put it another way, the first $p$ and the last $p$ elements of the sequence $D=D(\gamma)=\{d_1,d_2,\ldots,d_k\}$ are equal to one. At the moment we have to interrupt for a while the description of the algorithm in order to prove the following lemma.
\begin{Le}\label{lemma9-1}
One has the estimate
\begin{gather}\label{9-8}
\#S_2\ge\frac{M^{2\delta}}{2^{13}A^4\log^2M}.
\end{gather}
\begin{proof}
The sequence $D$ can be represented in the form
Последовательность $D$  представим как
\begin{gather*}
D=\{\underbrace{1,1,\ldots,1}_p,b_1,b_2,\ldots,b_n,\underbrace{1,1,\ldots,1}_p\}
\end{gather*}
then the required inequality \eqref{9-2} can be represented in the form
\begin{gather}\label{9-9}
\frac{M}{64A^2}\le\langle\underbrace{1,1,\ldots,1}_p,b_1,b_2,\ldots,b_n,\underbrace{1,1,\ldots,1}_p\rangle\le M.
\end{gather}
Let prove that the inequality \eqref{9-9} follows from the inequality
\begin{gather}\label{9-10}
\frac{M}{64A^2\log M}\le\langle b_1,b_2,\ldots,b_n\rangle\le \frac{M}{16\log M}.
\end{gather}
Indeed, let the inequality \eqref{9-10} be true. Then on the one side it follows from inequalities \eqref{continuant inequality} and \eqref{9-5} that
\begin{gather*}
\langle\underbrace{1,1,\ldots,1}_p,b_1,b_2,\ldots,b_n,\underbrace{1,1,\ldots,1}_p\rangle\le4F_p^2
\langle b_1,b_2,\ldots,b_n\rangle\le16\langle b_1,b_2,\ldots,b_n\rangle\log M\le M,
\end{gather*}
and on the other side it follows in a similar way from the inequality \eqref{9-4} that
\begin{gather*}
\langle\underbrace{1,1,\ldots,1}_p,b_1,b_2,\ldots,b_n,\underbrace{1,1,\ldots,1}_p\rangle\ge F_p^2
\langle b_1,b_2,\ldots,b_n\rangle\ge\langle b_1,b_2,\ldots,b_n\rangle\log M\ge\frac{M}{64A^2}.
\end{gather*}
Thus the implication $\eqref{9-10}\Rightarrow\eqref{9-9}$ is proved. It remains to obtain a lower bound for the amount of sequences $B=\langle b_1,b_2,\ldots,b_n\rangle$ satisfying the inequality \eqref{9-10}. We set $x=\frac{M}{16\log M}$ and note that the condition $x\ge4A^2$ in Theorem \ref{Kan-Hen} follows from the inequality \eqref{9-1}. Consequently, considering \eqref{F_A less} one has
\begin{gather*}
\#S_2\ge F_{\A}(x)-F_{\A}\left(\frac{x}{4A^2}\right)\ge\frac{1}{32A^4}\left(\frac{M}{16\log M}\right)^{2\delta}\ge
\frac{M^{2\delta}}{2^{13}A^4\log^2M},
\end{gather*}
since $\delta<1.$ This completes the proof of the lemma.
\end{proof}
\end{Le}
Let return to the description of the algorithm.
  \item[Step 3]
For any $L$ in the interval
\begin{gather}\label{9-11}
\left[\frac{M}{64A^2}, M\right]
\end{gather}
consider the set $S_3(L)\subset S_2$ of matrices $\gamma\in S_2,$ for which the following inequality holds
\begin{gather}\label{9-12}
\max\left\{\frac{M}{64A^2}, L(1-\log^{-1}L)\right\}\le\|\gamma\|\le L.
\end{gather}
Here, we also have to interrupt the description of the algorithm in order to prove the following lemma.
\begin{Le}\label{lemma9-2}
There is a number $L$ in the interval \eqref{9-11} such that
\begin{gather}\label{9-13}
|S_3(L)|\ge\frac{L^{2\delta}}{2^{16}A^5\log^3L}.
\end{gather}
\begin{proof}
Let $t$ be the minimal positive integer number satisfying the inequality
\begin{gather}\label{9-14}
(1-\log^{-1}M)^t\le\frac{1}{64A^2}.
\end{gather}
Note that $t\le8A\log M.$ For $j=1,2,\ldots,t$ consider sets $s(j)$ each of them consists of matrices  $\gamma\in S_2,$ such that
\begin{gather*}
M(1-\log^{-1}M)^j\le\|\gamma\|\le M(1-\log^{-1}M)^{j-1}
\end{gather*}
Since $S_2\subset\bigcup_{1\le j\le t}s(j),$ by the pigeonhole principle there is a set among $s(j)$ containing at least $\frac{|S_2|}{t}$ matrices. Let $$L=M(1-\log^{-1}M)^{j_0-1},$$ then $L$ belongs to the segment \eqref{9-11} and $s(j_0)\subset S_3(L).$ Hence $|S_3(L)|\ge\frac{|S_2|}{t}.$ Using the bound \eqref{9-8} and the restriction on $t$ one has
\begin{gather}\label{9-17}
|S_3(L)|\ge\frac{M^{2\delta}}{(2^{13}A^4\log^2M)8A\log M}=\frac{M^{2\delta}}{2^{16}A^5\log^3M}.
\end{gather}
Because the function $f(x)=x^{2\delta}\log^{-3}x$ increases and since $M\ge L,$ then replacing in \eqref{9-17} the parameter $M$ by $L$ one has the inequality \eqref{9-13}. This completes the proof of the lemma.
\end{proof}
\end{Le}
Returning to the algorithm we choose in the interval \eqref{9-11} any $L$ (for example the maximal one) satisfying the inequality \eqref{9-13} and fix it. Now let $S_3:=S_3(L).$
  \item[Step 4]
For $\gamma\in S_3$ let $k(\gamma)$ be the length of the sequence $D(\gamma).$ Represent the set $S_3$ as the union of the sets $S_4(k),$ consisting of those matrices $\gamma\in S_3$ for which $k(\gamma)=k$ is fixed.
\begin{Le}\label{lemma9-3}
There exists $k$ for which
\begin{gather}\label{9-18}
|S_4(k)|\ge\frac{L^{2\delta}}{2^{18}A^5\log^4L}.
\end{gather}
\begin{proof}
Since for all $D\in V_{\A}(r)$ the inequality
$
\langle D\rangle\ge\langle \,\underbrace{1,1,\ldots,1}_{r}\,\rangle,
$
holds, then
\begin{equation*}
\langle D\rangle\ge\left(\frac{\sqrt{5}+1}{2}\right)^{r-1}
\end{equation*}
and consequently
\begin{gather}\label{9-19}
k\le\frac{\log\|\gamma\|}{\log\frac{\sqrt{5}+1}{2}}+1\le4\log\|\gamma\|\le4\log L.
\end{gather}
Hence, by the pigeonhole principle, there is a $k,$ for which
\begin{gather*}
|S_4(k)|\ge\frac{|S_3|}{4\log L}\ge\frac{L^{2\delta}}{(4\log L)2^{16}A^5\log^3L}=\frac{L^{2\delta}}{2^{18}A^5\log^4L}
\end{gather*}
by \eqref{9-13} and \eqref{9-19}. This completes the proof of the lemma.
\end{proof}
\end{Le}
Returning to the algorithm we fix any $k,$ satisfying the inequality \eqref{9-18} and write $S_4:=S_4(k).$
\textbf{Algorithm is completed.}
\end{description}
Now we write $\Xi(M):=S_4.$ Recall the properties of matrices $\gamma\in\Xi(M).$ For any $\gamma\in\Xi(M)$ we have from the construction:
\begin{description}
  \item[i] the first and the last $p$ elements of the sequence $D(\gamma)$ are equal to 1, where $p$ is defined by the inequality \eqref{9-4};
  \item[ii] $L(1-\log^{-1}L)\le\|\gamma\|\le L;$
  \item[iii] $k(\gamma)=const,$ that is, the length of $D(\gamma)$ is fixed fir all $\gamma\in\Xi(M).$
\end{description}
Besides, we have proved that
\begin{gather}\label{9-20}
\#\Xi(M)\ge\frac{L^{2\delta}}{2^{18}A^5\log^4L}.
\end{gather}
The first property allows us to prove an important lemma
\begin{Le}\label{pred9-1}
For every matrix $\gamma\in\Xi(M)$ of the form
$\gamma=
\begin{pmatrix}
a & b \\
c & d
\end{pmatrix}
$
the following inequalities hold
\begin{gather}\label{9-21}
\left|\frac{b}{d}-\frac{1}{\varphi}\right|\le\frac{2}{\log L},\quad
\left|\frac{c}{d}-\frac{1}{\varphi}\right|\le\frac{2}{\log L},
\end{gather}
where $\varphi$ is the golden ratio
\begin{gather}\label{9-22}
\varphi=1+[1,1,\ldots,1,\ldots]=\frac{\sqrt{5}+1}{2}.
\end{gather}
\begin{proof}
It follows from \eqref{continuant properties} and \eqref{matrix=continuant} that
$\frac{b}{d}=[D(\gamma)],\quad \frac{c}{d}=[\overleftarrow{D(\gamma)}].$ Hence the fraction
$\alpha=[\underbrace{1,1,\ldots,1}_p]$
is a convergent fraction to $\frac{b}{d}$ and to $\frac{c}{d}.$ The denominator of $\alpha$ is equal to $F_p$ and
it follows from the choice of parameters \eqref{9-4}, \eqref{9-11} that
$F_p\ge\log^{\frac{1}{2}}L.$
Hence,
\begin{gather}\label{9-24}
\left|\frac{b}{d}-\alpha\right|\le\frac{1}{\log L},\quad
\left|\frac{c}{d}-\alpha\right|\le\frac{1}{\log L}.
\end{gather}
But $\alpha$ is also a convergent fraction to $\frac{1}{\varphi},$ thus
$\left|\alpha-\frac{1}{\varphi}\right|\le\frac{1}{\log L}.$
Applying the triangle inequality we obtain the desired inequalities. This completes the proof of the lemma.
\end{proof}
\end{Le}
\section{Parameters and their properties}
Let $N\ge N_{min}=N_{min}(\epsilon_0,\A)$ and write
\begin{gather}\label{10-2}
J=J(N)=\left[\frac{\log\log N-4\log(10A)+2\log\epsilon_0}{-\log(1-\epsilon_0)}\right],
\end{gather}
where as usual $A\ge|\A|\ge2$ and require the following inequality $J(N_{min})\ge10$ to hold.
Using the definition \eqref{10-2}, one has
\begin{gather}\label{10-4}
\frac{10^4A^4}{\log N}\le\epsilon_0^2(1-\epsilon_0)^J\le\frac{10^5A^4}{\log N}.
\end{gather}
Now let define a finite sequence
\begin{gather}\label{10-6}
\left\{N_{-J-1},N_{-J},\ldots,N_{-1},N_0,N_1,\ldots,N_{J+1}\right\},
\end{gather}
having set $N_{J+1}=N$ and
\begin{gather}\label{10-8}
N_j=
\left\{
              \begin{array}{ll}
                N^{\frac{1}{2-\epsilon_0}(1-\epsilon_0)^{1-j}}, & \hbox{если $-1-J\le j\le1$;} \\
                N^{1-\frac{1}{2-\epsilon_0}(1-\epsilon_0)^{j}}, & \hbox{если $0\le j\le J$.}
              \end{array}
\right.
\end{gather}
It is obvious that the sequence is well-defined for $j=0$ and $j=1.$
\begin{Le}\label{lemma10-2}
\begin{enumerate}
  \item For $-J\le m\le J-1$ the following equation holds
  \begin{gather}\label{10-9}
N_{-m}N_{m+1}=N.
\end{gather}
  \item For $-J-1\le m\le J-1$ the following relations hold
\begin{gather}\label{10-10}
\frac{N_{m+1}}{N_m}=
\left\{
              \begin{array}{ll}
                N_{m+1}^{\epsilon_0}, & \hbox{if $m\le0$;} \\ \\
                \left(\frac{N}{N_m}\right)^{\epsilon_0}, & \hbox{if $m\ge0$,}
              \end{array}
\right.
\end{gather}
\begin{gather}\label{10-11}
\frac{N_{m+1}}{N_m}= N^{\frac{\epsilon_0}{2-\epsilon_0}(1-\epsilon_0)^{|m|}},
\end{gather}
\begin{gather}\label{10-12}
N_{m}\ge N_{m+1}^{1-\epsilon_0}.
\end{gather}
  \item Для $-1\le j<h\le J+1$ выполнено
\begin{gather}\label{10-13}
N_{h-J}^{(1-\epsilon_0)^{h-j}}=N_{j-J}.
\end{gather}
\end{enumerate}
\begin{proof}
All propositions follow directly from the definition  \eqref{10-8}. This completes the proof of the lemma.
\end{proof}
\end{Le}
\begin{Le}\label{lemma10-3}
For $-J\le m\le J-1$ the following estimate holds
\begin{gather}\label{10-16}
\frac{N_{m+1}}{N_{m}}\ge\exp\left(\frac{10^4A^4}{2\epsilon_0}\right);
\end{gather}
moreover
\begin{gather}\label{10-17}
\exp\left(\frac{10^4A^4}{2\epsilon_0^2}\right)\le\frac{N}{N_{J}}=N_{1-J}\le\exp\left(\frac{10^5A^4}{\epsilon_0^2}\right).
\end{gather}
\begin{proof}
The inequality \eqref{10-16} follows from \eqref{10-11} and the lower bound in \eqref{10-4}. Now let prove the inequality \eqref{10-17}. The equation $\frac{N}{N_{J}}=N_{1-J}$ follows from \eqref{10-9} with $m=-J.$  The estimate of $N_{1-J}$ follows from \eqref{10-8} and \eqref{10-4}. This completes the proof of the lemma.
\end{proof}
\end{Le}
\begin{Le}\label{lemma10-4}
For any $M,$ such that
\begin{gather}\label{10-19}
N_{1-J}\le M\le N_{J},
\end{gather}
there exist indexes  $j$ and $h,$ such that
\begin{gather}\label{10-20}
2\le j\le 2J,\quad 1\le h\le 2J-1,\quad h=2J-j+1,
\end{gather}
for which the following inequalities hold
\begin{gather}\label{10-21}
N_{j-J}^{1-\epsilon_0}\le M\le N_{j-J},
\end{gather}
\begin{gather}\label{10-22}
\left(\frac{N}{N_{h-J}}\right)^{1-\epsilon_0}\le M\le\frac{N}{N_{h-J}}.
\end{gather}
\begin{proof}
Since the sequence $\{N_j\}$ is increasing there exists the index $j$ in \eqref{10-20} such that
\begin{gather}\label{10-23}
N_{j-1-J}\le M\le N_{j-J},
\end{gather}
then \eqref{10-21} follows from \eqref{10-12}.
The inequality \eqref{10-22} can be obtained by substituting the equation \eqref{10-9} into \eqref{10-21}. This completes the proof of the lemma.
\end{proof}
\end{Le}
For a nonnegative integer number $n$ we write
\begin{gather}\label{10-26}
\widetilde{N}_{n-J}=
\left\{
              \begin{array}{ll}
                N_{n-J}, & \hbox{if $n\ge1$;} \\
                1, & \hbox{if $n=0.$}
              \end{array}
\right.
\end{gather}
Moreover for integers $j$ and $h$ such that
\begin{gather}\label{10-27}
0\le j<h\le 2J+1,
\end{gather}
we write
\begin{gather}\label{10-28}
j_0(j,h)=\min\left\{|n-J-1|\,\Bigl|\, j+1\le n\le h\right\}.
\end{gather}
Note that there are only three alternatives for the value of $j_0$
\begin{gather}\label{10-29}
j_0\in\left\{j-J,\,0,\,J+1-h\right\}.
\end{gather}
\begin{Le}\label{lemma10-5}
For integers $j$ and $h$ from \eqref{10-27} the following estimate holds
\begin{gather}\label{10-30}
\prod_{n=j+1}^{h}\left(2^9A^3\log\frac{N_{n-J}}{\widetilde{N}_{n-1-J}}\right)\le
\left((1-\epsilon_0)^{j_0}\log N\right)^{\frac{7}{4}(h-j)},
\end{gather}
where
$j_0=j_0(j,h).$
\begin{proof}
Consider two cases depending on the value of $j.$
\begin{description}
  \item[1)\,$j>0.$]
Using \eqref{10-16} if $n<2J+1$ and \eqref{10-17} if $n=2J+1,$ for any $n$ in the segment $j+1\le n\le h$ we obtain
\begin{gather}\label{10-31}
2^9A^3\le\frac{2^9}{10^3}(2\epsilon_0)^{3/4}\left(\log\frac{N_{n-J}}{N_{n-1-J}}\right)^{3/4}
\le\left(\log\frac{N_{n-J}}{N_{n-1-J}}\right)^{3/4}.
\end{gather}
Hence, since for $j>0$ one has $\widetilde{N}_{n-1-J}=N_{n-1-J},$ we obtain
\begin{gather}\label{10-32}
\prod_{n=j+1}^{h}\left(2^9A^3\log\frac{N_{n-J}}{\widetilde{N}_{n-1-J}}\right)\le
\prod_{n=j+1}^{h}\left(\log\frac{N_{n-J}}{N_{n-1-J}}\right)^{7/4}.
\end{gather}
Making allowance for $\epsilon_0\in\left(0,\frac{1}{2500}\right),$ it follows from \eqref{10-11} if  $n<2J+1$
and from \eqref{10-8} if $n=2J+1,$ that
\begin{gather}\label{10-33}
\log\frac{N_{n-J}}{N_{n-1-J}}\le
\frac{1}{2-\epsilon_0}(1-\epsilon_0)^{|n-J-1|}\log N\le(1-\epsilon_0)^{j_0}\log N
\end{gather}
Substituting the estimate \eqref{10-33} into \eqref{10-32} we obtain the statement of the theorem in case $j>0.$
  \item[2)\,$j=0.$]
Using the lower bound from \eqref{10-17} we obtain
\begin{gather}\label{10-34}
2^9A^3\le\frac{2^9}{10^3}(2\epsilon_0^2)^{3/4}\left(\log N_{1-J}\right)^{3/4}
\le\left(\log N_{1-J}\right)^{3/4}.
\end{gather}
It follows from the definition \eqref{10-26} and the result of the previous item that
\begin{gather}\notag
\prod_{n=j+1}^{h}\left(2^9A^3\log\frac{N_{n-J}}{\widetilde{N}_{n-1-J}}\right)=2^9A^3\log N_{1-J}
\prod_{n=2}^{h}\left(2^9A^3\log\frac{N_{n-J}}{\widetilde{N}_{n-1-J}}\right)\le\\\le
(\log N_{1-J})^{7/4}\left((1-\epsilon_0)^{j_0(1,h)}\log N\right)^{\frac{7}{4}(h-j-1)}.\label{10-35}
\end{gather}
We obtain by the definition \eqref{10-8} that
$$\log N_{1-J}=\frac{1}{2-\epsilon_0}(1-\epsilon_0)^{J}\log N\le(1-\epsilon_0)^{J}\log N.$$
Substituting this estimate into \eqref{10-35} we obtain
\begin{gather}\label{10-36}
\prod_{n=j+1}^{h}\left(2^9A^3\log\frac{N_{n-J}}{\widetilde{N}_{n-1-J}}\right)\le
\left((1-\epsilon_0)^{j_0(1,h)}\log N\right)^{\frac{7}{4}(h-j)}(1-\epsilon_0)^{\frac{7}{4}(J-j_0(1,h))}.
\end{gather}
Taking account of$j_0\le J$ and $1-\epsilon_0<1,$ we obtain
\begin{gather}\label{10-37}
\prod_{n=j+1}^{h}\left(2^9A^3\log\frac{N_{n-J}}{\widetilde{N}_{n-1-J}}\right)\le
\left((1-\epsilon_0)^{j_0(1,h)}\log N\right)^{\frac{7}{4}(h-j)}.
\end{gather}
The fact that $j_0(1,h)=j_0(0,h)$ completes the proof.
\end{description}
Lemma is proved.
\end{proof}
\end{Le}
\section{The ensemble: constructing the set $\Omega_N.$ }\label{postoenie}
In this section we construct a set $\Omega_N.$ It will be proved in §\ref{ansambl-property} that this set is $(\epsilon_0,N)$~--ensemble. We construct the set by the inductive algorithm with the steps numbered by indexes $1,2,\ldots,2J+1.$
\begin{enumerate}
  \item The first (the starting) step.\par
We set
\begin{gather}\label{11-1}
M=M_1=N_{1-J}.
\end{gather}
Because of the lower bound in \eqref{10-17} the condition \eqref{9-1} obviously holds. So we can run the algorithm of §\ref{predansambl} to generate the set $\Xi(M).$ During the algorithm we also obtain the numbers $L=L(M),p=p(M),k=k(M).$ By the construction the number $L$ belongs to the segment $\left[\frac{M}{64A^2}, M\right],$ so we can set
\begin{gather}\label{11-2}
L=\alpha_1M_1=\alpha_1 N_{1-J},
\end{gather}
where $\alpha_1$ is a number from
\begin{gather}\label{11-3}
\left[\frac{1}{64A^2}, 1\right].
\end{gather}
Let rename the returned pre-ensemble $\Xi(M)$ and numbers $L,p$ and $k$ to
\begin{gather*}
\Xi(M)=\Xi_1,\,L=L_1,\,p=p_1,\,k=k_1.
\end{gather*}
For the next step of the algorithm we define the number $M_2$
\begin{gather}\label{11-4}
M_2=\frac{N_{2-J}}{(1+\varphi^{-2})\alpha_1 N_{1-J}},
\end{gather}
where $\varphi$ id from \eqref{9-22}.
  \item The step with the number $j,$ where $2\le j\le2J+1$ (the inductive step).\par
Write $M=M_j.$ According to the inductive assumption the number $M_j$ has been defined on the previous step by the formula
\begin{gather}\label{11-5}
M_j=\frac{N_{j-J}}{(1+\varphi^{-2})\alpha_{j-1} N_{j-1-J}},
\end{gather}
where $\alpha_{j-1}$ is a number from \eqref{11-3}. To verify for such $M$ the condition \eqref{9-1} it is sufficient to apply bounds of Lemma \ref{lemma10-3} having put $m=j-J.$
Hence we can run the algorithm of §\ref{predansambl} to generate $\Xi(M).$ Besides there exists a number $\alpha_j$ from the interval \eqref{11-3} such that for the parameter $L$ the following equation holds
\begin{gather}\label{11-6}
L=\alpha_jM=\frac{\alpha_j N_{j-J}}{(1+\varphi^{-2})\alpha_{j-1} N_{j-1-J}}.
\end{gather}
We rename $\Xi(M)$ to $\Xi_j,$ the number $L$ to $L_j,$ the quantity $p$ to $p_j,$ the length $k$ to $k_j.$  If $j\le 2J,$ then the number $M_{j+1},$ which will be used in the next step, is defined by the equation
\begin{gather*}
M_{j+1}=\frac{N_{j+1-J}}{(1+\varphi^{-2})\alpha_{j} N_{j-J}}.
\end{gather*}
If $j=2J+1,$ then the notation $M_{j+1}$ is of no use, as the algorithm is completed.
\end{enumerate}
We now define the ensemble $\Omega_N$ writing all the sets generated in the algorithm for one another
\begin{gather*}
\Omega_N=\Xi_1\Xi_2\ldots\Xi_{2J}\Xi_{2J+1}.
\end{gather*}
It means that the set $\Omega_N$ consists of all possible products of the form

$\gamma_1\gamma_2\ldots\gamma_{2J}\gamma_{2J+1},$ with $\gamma_1\in\Xi_1,\,
\gamma_2\in\Xi_2,\ldots,
\gamma_{2J+1}\in\Xi_{2J+1}.$\par
To prove that $\Omega_N$ is really an ensemble we need two technical lemmas concerning quantities $L_j.$
We will use the following notation
\begin{gather}\label{12-4}
f(x)=O_1(g(x)),\,\mbox{if}\,|f(x)|\le  g(x).
\end{gather}

\begin{Le}\label{lemma12-2}
The following inequality holds
\begin{gather}\label{12-7}
\sum_{n=1}^{2J+1}\frac{1}{\log L_n}\le\frac{1}{16000}.
\end{gather}
\begin{proof}
Let prove that numbers  $L_j$ satisfy the following inequality
\begin{gather}\label{12-8}
L_j\ge\frac{1}{64A^2(1+\varphi^{-2})}\frac{N_{j-J}}{N_{j-J-1}}\ge\frac{N_{j-J}}{100A^2N_{j-J-1}}.
\end{gather}
Actually, for $j>1$ the inequality \eqref{12-8} follows from the definition \eqref{11-6}.
For $j=1$ to deduce the same inequality \eqref{12-8} from \eqref{11-2} it is sufficient to know that $N_{-J}\ge1.$ The last inequality follows from inequalities \eqref{10-12} and \eqref{10-17} with $m=-J:$
\begin{gather*}
N_{-J}\ge N_{1-J}^{1-\epsilon_0}\ge\exp\left(\frac{10^4A^4}{2\epsilon_0^2}(1-\epsilon_0)\right)\ge1.
\end{gather*}
Thus, the inequality \eqref{12-8} is proved. It follows from the bound \eqref{10-16} that
\begin{gather*}
\left(\frac{N_m}{N_{m-1}}\right)^{1/2}\ge\exp\left(\frac{10^4A^4}{4}\right)>100A^2,\quad -J\le m\le J-1.
\end{gather*}
In that case, if $j\le2J,$ then using \eqref{10-11} the estimate \eqref{12-8} can be resumed
\begin{gather}\label{12-9}
L_j\ge\left(\frac{N_{j-J}}{N_{j-J-1}}\right)^{\frac{1}{2}}\ge N^{\frac{1}{4}\epsilon_0(1-\epsilon_0)^{|j-J-1|}},
\end{gather}
Hence
\begin{gather}\label{12-10}
\log L_j\ge\frac{1}{4}\epsilon_0(1-\epsilon_0)^{|j-J-1|}\log N,\quad j\le2J.
\end{gather}
If $j=2J+1,$ then from the lower bound in \eqref{10-17} we obtain in a similar way
\begin{gather*}
L_{2J+1}\ge\left(\frac{N_{J+1}}{N_{J}}\right)^{\frac{1}{2}}\ge \exp\left(\frac{10^4A^4}{4\epsilon_0^2}\right),
\end{gather*}
whence it follows that
\begin{gather}\label{12-11}
\log L_j\ge\frac{10^4A^4}{4\epsilon_0^2}.
\end{gather}
Substituting the estimates \eqref{12-10} and \eqref{12-11} into the sum in \eqref{12-7} one has
\begin{gather}\label{12-12}
\sum_{n=1}^{2J+1}\frac{1}{\log L_n}\le\frac{4}{\epsilon_0}\sum_{n=0}^{2J}\frac{1}{(1-\epsilon_0)^{|n-J|}\log N}+\frac{4\epsilon_0^2}{10^4A^4}\le\frac{8}{\epsilon_0}\sum_{n=0}^{J}\frac{1}{(1-\epsilon_0)^{n}\log N}+\frac{4\epsilon_0^2}{10^4A^4}.
\end{gather}
Let estimate the geometric progression from  \eqref{12-12}:
\begin{gather*}
\sum_{n=0}^{J}\frac{1}{\epsilon_0(1-\epsilon_0)^{n}}\le\frac{1}{\epsilon_0(1-\epsilon_0)^{J}}\sum_{n=0}^{\infty}(1-\epsilon_0)^{n}\le
\frac{1}{\epsilon_0^2(1-\epsilon_0)^{J}}\le\frac{\log N}{10^4A^4},
\end{gather*}
since \eqref{10-4}. Substituting this bound into \eqref{12-12} we obtain
\begin{gather*}
\sum_{n=1}^{2J+1}\frac{1}{\log L_n}\le\frac{8}{\log N}\frac{\log N}{10^4A^4}+\frac{4\epsilon_0^2}{10^4A^4}=
\frac{4(2+\epsilon_0^2)}{10^4A^4}<\frac{1}{10^3A^4}\le\frac{1}{16000}.
\end{gather*}
This completes the proof of the lemma.
\end{proof}
\end{Le}
To state the following lemma we suppose the real numbers
\begin{gather*}
\Pi_1,\Pi_2,\ldots,\Pi_{2J+1}
\end{gather*}
to satisfy the relations
\begin{gather}\label{12-13}
\Pi_j=\left(1+2O_1(\log^{-1}L_j)\right)^2\prod_{n=1}^{j-1}\left(1+2O_1(\log^{-1}L_n)\right)^3,
\end{gather}
where the product over the empty set is regarded to be equal to one.
\begin{Le}\label{lemma12-3}
For any $j=1,2,\ldots,2J+1$ the following bound holds
\begin{gather}\label{12-14}
\exp(-10^{-3})\le\Pi_j\le\exp(10^{-3}).
\end{gather}
\begin{proof}
Taking the logarithm of the equation \eqref{12-13} and bounding from above the absolute value of the sum by the sum of absolute values we obtain
\begin{gather}\notag
|\log\Pi_j|\le2|\log\left(1+2O_1(\log^{-1}L_j)\right)|+3\sum_{n=1}^{j-1}|\log\left(1+2O_1(\log^{-1}L_n)\right)|\le\\\le
3\sum_{n=1}^{2J+1}|\log\left(1+2O_1(\log^{-1}L_n)\right)|.\label{12-15}
\end{gather}
In view of Lemma \ref{lemma12-2}, every number $\log^{-1}L_n$ for $n=1,2,\ldots,2J+1$ is less than $\frac{1}{16000};$ and in particular every number $2O_1(\log^{-1}L_n)$ belongs to the segment $\left[-\frac{1}{2},\frac{1}{2}\right].$ But for any $z$ in the segment $-\frac{1}{2}\le z\le\frac{1}{2}$ the inequality $|\log(1+z)|\le |z|\log4$ holds. Then by \eqref{12-15} we obtain
\begin{gather}\label{12-21}
|\log\Pi_j|\le3\sum_{n=1}^{2J+1}\left|2O_1(\log^{-1}L_n)\right|\log4<12\sum_{n=1}^{2J+1}\log^{-1}L_n.
\end{gather}
Using Lemma \ref{lemma12-2} we obtain
\begin{gather}\label{12-22}
|\log\Pi_j|\le\frac{12}{16000}<10^{-3}.
\end{gather}
The inequality \eqref{12-14} follows immediately from  \eqref{12-22}. This completes the proof of the lemma.
\end{proof}
\end{Le}
\section{Properties of $\Omega_N.$ It is really an ensemble!}\label{ansambl-property}
In this section we prove that the constructed set $\Omega_N$ is an ensemble, that is, it satisfies the definition of ensemble in §\ref{section idea BK}. Unique expansion is the easiest property to verify. Actually, if
$$\Omega^{(1)}=\Xi_1\Xi_2\ldots\Xi_{j},$$
$$\Omega^{(2)}=\Xi_{j+1}\Xi_{j+2}\ldots\Xi_{2J+1},$$
then firstly $\Omega_N=\Omega^{(1)}\Omega^{(2)}.$ Secondly, as the representation of a matrix in the form \eqref{matrix} is unique then the implication $\eqref{8-19}\Rightarrow\eqref{8-20}$ holds since the length $D(\gamma)=k_j$ is fixed for all $\gamma\in\Xi_{j}$ (the property (iii) in §\ref{section idea BK}), for each $j=1,2,\ldots,2J+1.$ \par
The next purpose is to prove that $\Omega_N$ is a pre-ensemble.
\begin{Le}\label{lemma13-1}
For any $j$ in the segment
\begin{gather}\label{13-1}
1\le j\le2J+1,
\end{gather}
for any  collection of matrices
\begin{gather}\label{13-2}
\xi_1\in\Xi_1,\,\xi_2\in\Xi_2,\ldots,\xi_j\in\Xi_j,
\end{gather}
one can find a number $\Pi_j,$ satisfying the equality \eqref{12-13}, such that
\begin{gather}\label{13-3}
\|\xi_1\xi_2\ldots\xi_j\|=\alpha_jN_{j-J}\Pi_j.
\end{gather}
\begin{proof}
Let first $j=1.$ Then, by the construction of the pre-ensemble $\Xi_1$ (§\ref{predansambl}) and by the equation \eqref{11-2}, the following equation holds
\begin{gather}\label{13-4}
\|\xi_1\|=\alpha_1N_{1-J}(1+O_1(\log^{-1}L_1)).
\end{gather}
Since
\begin{gather*}
1+O_1(\log^{-1}L_1)=(1+2O_1(\log^{-1}L_1))^2=\Pi_1,
\end{gather*}
then substituting the last equation into \eqref{13-4} one has
\begin{gather}\label{13-5}
\|\xi_1\|=\alpha_1N_{1-J}\Pi_1,
\end{gather}
and in the case $j=1$ lemma is proved.\par
We now assume that lemma is proved for some $j,$ such that $1\le j\le2J,$ and prove that it holds for $j+1.$ It follows from
\eqref{continuant identity} that
\begin{gather}\label{13-6}
\|\xi_1\xi_2\ldots\xi_j\xi_{j+1}\|=\|\xi_1\xi_2\ldots\xi_j\|\|\xi_{j+1}\|
(1+[\overleftarrow{D}(\xi_j),\overleftarrow{D}(\xi_{j-1}),\ldots,\overleftarrow{D}(\xi_1)][D(\xi_{j+1})]),
\end{gather}
where $D(\gamma),$ as usual, denotes the sequence $D(\gamma)=\{d_1,d_2,\ldots,d_k\},$ where
\begin{equation*}
\gamma=
\begin{pmatrix}
0 & 1 \\
1 & d_1
\end{pmatrix}
\begin{pmatrix}
0 & 1 \\
1 & d_2
\end{pmatrix}\ldots
\begin{pmatrix}
0 & 1 \\
1 & d_k
\end{pmatrix}.
\end{equation*}
It follows immediately from Lemma \ref{pred9-1} that
\begin{gather}\label{13-7}
[D(\xi_{j+1})]=\varphi^{-1}+2O_1(\log^{-1}L_{j+1}),\,
[\overleftarrow{D}(\xi_j)\ldots,\overleftarrow{D}(\xi_1)]=\varphi^{-1}+2O_1(\log^{-1}L_{j}).
\end{gather}
Substituting \eqref{13-7} into \eqref{13-6}, we obtain
\begin{gather}\label{13-9}
\|\xi_1\xi_2\ldots\xi_j\xi_{j+1}\|=\|\xi_1\xi_2\ldots\xi_j\|\|\xi_{j+1}\|
\left(1+\varphi^{-2}\right)\left(1+2O_1(\log^{-1}L_j)\right)\left(1+2O_1(\log^{-1}L_{j+1})\right).
\end{gather}
By the inductive hypothesis we have
\begin{gather}\label{13-10}
\|\xi_1\xi_2\ldots\xi_j\|=\alpha_jN_{j-J}\Pi_j.
\end{gather}
And by the construction of the ensemble $\Omega_N$ (§\ref{postoenie}, "The step with the number $j+1.$" ) the following equation holds
\begin{gather}\label{13-11}
\|\xi_{j+1}\|=L_{j+1}(1+O_1(\log^{-1}L_{j+1})=
\frac{\alpha_{j+1} N_{j+1-J}}{(1+\varphi^{-2})\alpha_{j} N_{j-J}}(1+O_1(\log^{-1}L_{j+1})).
\end{gather}
Substituting \eqref{13-10} and \eqref{13-11} into \eqref{13-9} and making cancelations, we obtain
\begin{gather}\label{13-12}
\|\xi_1\xi_2\ldots\xi_j\xi_{j+1}\|=
\alpha_{j+1} N_{j+1-J}\widetilde{\Pi}_{j+1},
\end{gather}
where
\begin{gather}\label{13-13}
\widetilde{\Pi}_{j+1}=(1+O_1(\log^{-1}L_{j+1}))\left(1+2O_1(\log^{-1}L_j)\right)\left(1+2O_1(\log^{-1}L_{j+1})\right)\Pi_j.
\end{gather}
Using the definition of $\Pi_j$ by the equation \eqref{12-13} we obtain
\begin{gather*}
\widetilde{\Pi}_{j+1}=(1+2O_1(\log^{-1}L_{j+1}))^2\left(1+2O_1(\log^{-1}L_j)\right)\Pi_j=\Pi_{j+1}
\end{gather*}
and hence
\begin{gather*}
\|\xi_1\xi_2\ldots\xi_j\xi_{j+1}\|=
\alpha_{j+1} N_{j+1-J}\Pi_{j+1}.
\end{gather*}
The lemma is proved.
\end{proof}
\end{Le}
\begin{Le}\label{lemma13-2}
For any  collection of matrices \eqref{13-2}, for any numbers $j,\,h$ in the interval
\begin{gather}\label{13-14}
1\le j\le2J+1,\quad j<h\le 2J+1
\end{gather}
the following inequalities hold
\begin{gather}\label{13-15}
\frac{1}{70A^2}N_{j-J}\le\|\xi_1\xi_2\ldots\xi_j\|\le1,01N_{j-J},
\end{gather}
\begin{gather}\label{13-16}
\frac{1}{70A^2}N\le\|\xi_1\xi_2\ldots\xi_{2J+1}\|\le1,01N,
\end{gather}
\begin{gather}\label{13-17}
\frac{1}{150A^2}\frac{N_{h-J}}{N_{j-J}}\le\|\xi_{j+1}\xi_{j+2}\ldots\xi_h\|\le73A^2\frac{N_{h-J}}{N_{j-J}};
\end{gather}
moreover, for $j\le2J$ one has
\begin{gather}\label{13-18}
\frac{1}{150A^2}\frac{N}{N_{j-J}}\le\|\xi_{j+1}\xi_{j+2}\ldots\xi_{2J+1}\|\le73A^2\frac{N}{N_{j-J}}.
\end{gather}
\begin{proof}
First we prove the inequality  \eqref{13-15}. Recall that by the construction of the set $\Omega_N$ the following inequality holds
\begin{gather}\label{13-19}
\frac{1}{64A^2}\le\alpha_{j}\le1,
\end{gather}
and by Lemma \ref{lemma12-3} we also have
\begin{gather}\label{13-20}
\exp(-10^{-3})\le\Pi_j\le\exp(10^{-3}).
\end{gather}
Substituting \eqref{13-19} and \eqref{13-20} into \eqref{13-3}, we obtain \eqref{13-15}. In particular, as $N_{J+1}=N,$ then by using \eqref{13-15} for $j=2J+1$ we obtain \eqref{13-16}.\par
Now we prove the inequality \eqref{13-17}. To do this we denote
\begin{gather*}
W(j,h)=\|\xi_{j}\xi_{j+1}\ldots\xi_h\|
\end{gather*}
and rewrite the inequality \eqref{continuant inequality} in the form
\begin{gather*}
W(1,j)W(j+1,h)\le W(1,h)\le2W(1,j)W(j+1,h).
\end{gather*}
Hence, applying the inequality \eqref{13-15} twice we obtain
\begin{gather*}
W(j+1,h)\ge\frac{W(1,h)}{2W(1,j)}\ge\frac{N_{h-J}/(70A^2)}{2,02N_{j-J}}\ge\frac{1}{150A^2}\frac{N_{h-J}}{N_{j-J}},
\end{gather*}
and in the same way
\begin{gather*}
W(j+1,h)\le\frac{W(1,h)}{W(1,j)}\le\frac{1,01N_{h-J}}{N_{j-J}/(70A^2)}\le73A^2\frac{N_{h-J}}{N_{j-J}}.
\end{gather*}
These prove the inequality \eqref{13-17}. Putting $h=2J+1$ in it we obtain \eqref{13-18}.
The lemma is proved.
\end{proof}
\end{Le}
For integers $j$ and $h,$ such that
\begin{gather}\label{13-21}
0\le j<h\le 2J+1,
\end{gather}
we put
\begin{gather}\label{13-22}
\Omega(j,h)=\Xi_{j+1}\Xi_{j+2}\ldots\Xi_{h}.
\end{gather}
\begin{Le}\label{lemma13-3}
The following estimate holds
\begin{gather}\label{13-23}
\left|\Omega(0,j)\right|\le9N_{j-J}^{2\delta}.
\end{gather}
\begin{proof}
By definition,  for $\gamma\in\Omega(0,j)$ one has $\gamma=\xi_{1}\xi_{2}\ldots\xi_{j}$ for a collection of matrices \eqref{13-2}. So, it follows from the inequality  \eqref{13-15} that
\begin{gather}\label{13-24}
\|\gamma\|\le1,01N_{j-J}.
\end{gather}
The number of matrices $\gamma,$ satisfying the inequality \eqref{13-24} can be bounded by Theorem \ref{Kan-Hen}. Estimating the result from above we obtain \eqref{13-23}. This completes the proof of the lemma.
\end{proof}
\end{Le}
Recall that parameters $\widetilde{N}_{n-J}$ and $j_0(j,h)$ were introduced by formulae \eqref{10-26} and \eqref{10-28}.  Note that the restrictions \eqref{13-21} on $j$ and $h$ coincide with the restrictions \eqref{10-27}.
\begin{Le}\label{lemma13-4}
For $j$ and $h$ in \eqref{13-21} the following bound holds
\begin{gather}\label{13-25}
\left|\Omega(j,h)\right|\ge\frac{1}{\left((1-\epsilon_0)^{j_0}\log N\right)^{7(h-j)}}
\left(\frac{N_{h-J}}{\widetilde{N}_{j-J}}\right)^{2\delta},
\end{gather}
where $j_0=j_0(j,h).$
\begin{proof}
Multiplying the lower bounds \eqref{9-20} we obtain
\begin{gather}\label{13-26}
\left|\Omega(j,h)\right|\ge\prod_{n=j+1}^{h}\left|\Xi_n\right|\ge\prod_{n=j+1}^{h}\frac{L_n^{2\delta}}{2^{18}A^5\log^4L_n}.
\end{gather}
It follows from formulae \eqref{11-2} and \eqref{11-6} that
\begin{gather}\label{13-27}
\frac{N_{n-J}}{c_1\widetilde{N}_{n-1-J}}\le L_n\le\frac{c_1N_{n-J}}{\widetilde{N}_{n-1-J}},
\end{gather}
where
\begin{gather}\label{13-28}
c_1=64A^2(1+\varphi^{-2})\le 2^7A^2.
\end{gather}
Let estimate the product of the numerators in \eqref{13-26}. Applying \eqref{13-27} and \eqref{13-28} we obtain
\begin{gather*}
\prod_{n=j+1}^{h}L_n^{2\delta}\ge\prod_{n=j+1}^{h}\left(\frac{N_{n-J}}{2^7A^2\widetilde{N}_{n-1-J}}\right)^{2\delta}
\end{gather*}
After the cancelations we obtain
\begin{gather}\label{13-29}
\prod_{n=j+1}^{h}L_n^{2\delta}\ge
\left(\frac{N_{h-J}}{\widetilde{N}_{j-J}}\right)^{2\delta}\prod_{n=j+1}^{h}(2^7A^2)^{-2\delta}.
\end{gather}
So the estimate \eqref{13-26} can be resumed in such a way
\begin{gather}\notag
\left|\Omega(j,h)\right|\ge
\left(\frac{N_{h-J}}{\widetilde{N}_{j-J}}\right)^{2\delta}
\prod_{n=j+1}^{h}\frac{(2^7A^2)^{-2\delta}}{2^{18}A^5\log^4L_n}\ge\\\ge
\left(\frac{N_{h-J}}{\widetilde{N}_{j-J}}\right)^{2\delta}
\left(\prod_{n=j+1}^{h}(2^8A^3\log L_n)\right)^{-4}.\label{13-30}
\end{gather}
The last product in  \eqref{13-30} will be estimated separately. In view of the upper bound in \eqref{13-27} we have
\begin{gather}\label{13-31}
\prod_{n=j+1}^{h}(2^8A^3\log L_n)\le\prod_{n=j+1}^{h}\left(2^8A^3\left(
\log(2^7A^2)+\log\frac{N_{n-J}}{\widetilde{N}_{n-1-J}}\right)\right).
\end{gather}
Applying Lemma \ref{lemma10-3}, we obtain
\begin{gather}\label{13-31}
\log\frac{N_{n-J}}{\widetilde{N}_{n-1-J}}\ge\frac{10^4A^4}{2\epsilon_0}\ge\log(2^7A^2),
\end{gather}
hence, applying Lemma \ref{lemma10-5}, we obtain
\begin{gather}\label{13-32}
\prod_{n=j+1}^{h}(2^8A^3\log L_n)\le\prod_{n=j+1}^{h}\left(2^9A^3\log\frac{N_{n-J}}{\widetilde{N}_{n-1-J}}\right)\le
\left((1-\epsilon_0)^{j_0}\log N\right)^{\frac{7}{4}(h-j)}.
\end{gather}
Substituting the estimate \eqref{13-32} into \eqref{13-30}, we obtain \eqref{13-25}. This completes the proof of the lemma.
\end{proof}
\end{Le}
\begin{Th}\label{theorem13-1}
For $j$ and $h$ in \eqref{13-21} the following estimate holds
\begin{gather}\label{13-34}
\left|\Omega(j,h)\right|\ge
\left(\frac{N_{h-J}}{\widetilde{N}_{j-J}}\right)^{2\delta}
\exp\left(-\left(\frac{\log\log N}{\log(1-\epsilon_0)}+j_0\right)^2\right),
\end{gather}
where $j_0=j_0(j,h).$
\begin{proof}
It follows from \eqref{13-25} that it is enough to prove the inequality
\begin{gather*}
\exp\left(\left(\frac{\log\log N}{\log(1-\epsilon_0)}+j_0\right)^2\right)\ge
\exp\left(7(h-j)\left(\log\log N+j_0\log(1-\epsilon_0)\right)\right).
\end{gather*}
Hence, it is sufficient to prove
\begin{gather*}
7(h-j)\left(\log\log N+j_0\log(1-\epsilon_0)\right)\le
\left(\frac{\log\log N+j_0\log(1-\epsilon_0)}{\log(1-\epsilon_0)}\right)^2.
\end{gather*}
One can readily obtain from \eqref{10-2} that
\begin{gather}\label{13-35}
J\le\frac{\log\log N}{-\log(1-\epsilon_0)}-1,
\end{gather}
and since $j_0\le J,$ so one has
\begin{gather*}
\log\log N+j_0\log(1-\epsilon_0)>0.
\end{gather*}
Thus, it is sufficient to prove that
\begin{gather}\label{13-36}
7(h-j)\le
\frac{\log\log N+j_0\log(1-\epsilon_0)}{\log^2(1-\epsilon_0)}=
\frac{1}{-\log(1-\epsilon_0)}\left(\frac{\log\log N}{-\log(1-\epsilon_0)}-j_0\right).
\end{gather}
We observe that as $\epsilon_0\in\left(0,\frac{1}{2500}\right),$ so $$-7\log(1-\epsilon_0)\le\frac{1}{7}.$$ It follows from this, \eqref{13-35} and \eqref{13-36} that it is sufficient to prove
\begin{gather}\label{13-37}
\frac{1}{7}(h-j)\le J+1-j_0.
\end{gather}
For $j_0=0$ the inequality \eqref{13-37} follows from the trivial bound $h-j\le2J+1.$ For $j_0\neq0$ $\left( \mbox{hence,} j>J\,\mbox{or}\,h\le J\right)$ it follows from \eqref{10-28} and \eqref{10-29} that
\begin{gather}\label{13-38}
h-j\le
\left\{
              \begin{array}{ll}
                2J+1-j=J+1-j_0, & \hbox{если $j>J$;} \\
                h=J+1-j_0, & \hbox{если $h\le J$.}
              \end{array}
\right.
\end{gather}
Hence the inequality \eqref{13-37} follows. This completes the proof of the theorem.
\end{proof}
\end{Th}
\begin{Sl}\label{conseq13-1}
For $N\ge\exp(\epsilon_0^{-5})$ the following estimate holds
\begin{gather}\label{13-39}
N^{2\delta-\epsilon_0}\le N^{2\delta}
\exp\left(-\left(\frac{\log\log N}{\epsilon_0}\right)^2\right)\le
\#\Omega_N\le9N^{2\delta}.
\end{gather}
\begin{proof}
To prove the upper bound we apply Lemma \ref{lemma13-3} with $j=0,\,h=2J+1.$ To prove the lower bound we use Theorem \ref{theorem13-1} putting $j=0,\,h=2J+1$ in it (hence,  $j_0=0.$)
As $\epsilon_0\in\left(0,\frac{1}{2500}\right),$ so we have $\log^2(1-\epsilon_0)>\epsilon_0^2$ and obtain the lower bound in \eqref{13-39}. The corollary is proved.
\end{proof}
\end{Sl}
\begin{Sl}\label{conseq13-2}
The set $\Omega_N$ is a $N$~--pre-ensemble.
\begin{proof}
It follows from the inequality \eqref{13-16} that the property $\|\gamma\|\asymp N$ holds for each matrix $\gamma\in\Omega_N.$ By the Corollary \ref{conseq13-1} we have $\#\Omega_N\ge_{\epsilon}N^{2\delta-\epsilon}.$  To put it another way, both items of the definition of pre-ensemble hold. The corollary is proved.
\end{proof}
\end{Sl}
We now verify the last property of ensemble related to the relations  \eqref{8-21} and \eqref{8-22}.
\begin{Le}\label{lemma13-5}
For any $M$ in the interval \eqref{10-19} there exist $j$ and $h$ in the intervals \eqref{10-20}, such that for any collection of matrices \eqref{13-2} the following inequalities hold
\begin{gather}\label{13-40}
0,99\|\xi_1\xi_2\ldots\xi_j\|^{1-\epsilon_0}\le M\le
70A^2\|\xi_1\xi_2\ldots\xi_j\|,
\end{gather}
\begin{gather}\label{13-41}
\frac{1}{73A^2}\|\xi_{h+1}\xi_{h+2}\ldots\xi_{2J+1}\|^{1-\epsilon_0}\le M\le
150A^2\|\xi_{h+1}\xi_{h+2}\ldots\xi_{2J+1}\|.
\end{gather}
\begin{proof}
Let $M$ be fixed in the interval \eqref{10-19}. Then using Lemma \ref{lemma10-4} we find $j$ and $h$ in \eqref{10-20}, such that the inequalities  \eqref{10-21} and \eqref{10-22} hold. We note that a bilateral bound on the number $N_{j-J}$ in terms of $\|\xi_1\xi_2\ldots\xi_j\|$ follows from \eqref{13-15}. To obtain the inequality \eqref{13-40} one should substitute this bilateral bound into \eqref{10-21}. To obtain the inequality \eqref{13-41} one should substitute a bilateral bound on $N/N_{h-J},$ following from the inequality \eqref{13-18} with $j=h,$ into \eqref{10-22}. This completes the proof of the lemma.
\end{proof}
\end{Le}
\begin{Sl}\label{conseq13-3}
For any $M,$ satisfying the inequality
\begin{gather}\label{13-42}
\exp\left(\frac{10^5A^4}{\epsilon_0^2}\right)\le M\le N\exp\left(-\frac{10^5A^4}{\epsilon_0^2}\right),
\end{gather}
there exist indexes $j$ and $h$ in the intervals \eqref{10-20}, such that for any collection of matrices \eqref{13-2}the inequalities \eqref{13-40} and \eqref{13-41} hold.
\begin{proof}
By applying Lemma \ref{lemma10-3}, it follows from the inequality \eqref{13-42} that the inequality \eqref{10-19} holds. It is sufficient now to apply Lemma \ref{lemma13-5}. The corollary is proved.
\end{proof}
\end{Sl}
\begin{Th}\label{theorem13-2}
The set $\Omega_N$ is a bilateral $(\epsilon_0,N)$~--ensemble.
\begin{proof}
The unique expansion of the products, which are equal to $\Omega_N,$ has been proved in the beginning of §\ref{ansambl-property}. The property of $\Omega_N$ to be a $N$~--pre-ensemble has been proved in the Corollary \ref{conseq13-2}. The right $(\epsilon_0,N)$~--ensemble property is proved by the inequality \eqref{13-41}, the left~--by \eqref{13-40}, since it follows from the Corollary \ref{conseq13-3} that these inequalities hold for any $M,$ satisfying the inequality \eqref{13-42}. This completes the proof of the theorem.
\end{proof}
\end{Th}
Thus the main purpose of the section is achieved. However, we formulate a few more properties of a bilateral ensemble. These properties will be of use while estimating exponential sums.\par
For $M,$ satisfying the inequality \eqref{13-42}, we denote by $\jj$ and $\h$ the numbers $j$ and $h$ from Lemma \ref{lemma10-4} corresponding to $M$. For brevity we will write that numbers $\jj,\,\h$ corresponds to $M.$ In the following theorem $\jj^{(1)}$ corresponds to $M^{(1)},$ and $\h^{(3)}$ corresponds to $M^{(3)}.$
\begin{Le}\label{lemma13-6}
Let the inequality
\begin{gather}\label{13-43}
M^{(1)}M^{(3)}< N^{1-\epsilon_0}.
\end{gather}
holds for $M^{(1)}$ and $M^{(3)}$ in the interval \eqref{13-42}. Then $\jj^{(1)}<\h^{(3)},$ and for $M=M^{(1)}$ the inequality \eqref{13-40} holds and for $M=M^{(3)}.$ the inequality \eqref{13-41} holds.
\begin{proof}
It is sufficient to verify the condition $\jj^{(1)}<\h^{(3)}.$ Then the statement of Lemma \ref{lemma13-6} will follow from Lemma \ref{lemma13-5} and Corollary \ref{conseq13-3}.\par
We recall that $\h^{(3)}=2J-\jj^{(3)}+1.$ Hence, it is sufficient to prove that $$\jj^{(1)}+\jj^{(3)}<2J+1.$$ Assume the contrary.
It follows from \eqref{10-23} that
\begin{gather*}
N_{\jj^{(1)}-1-J}N_{\jj^{(3)}-1-J}\le M^{(1)}M^{(3)}.
\end{gather*}
Since $\jj^{(1)}+\jj^{(3)}\ge2J+1,$ one can consider the inequality $\jj^{(1)}\ge J+1$ to hold.  It follows from the relation $\jj^{(3)}-1-J\ge-\jj^{(1)}+J$ and the increasing of the sequence $\left\{N_j\right\}_{j=-J-1}^{J+1},$ that
\begin{gather*}
N_{\jj^{(1)}-1-J}N_{-\jj^{(1)}+J}\le M^{(1)}M^{(3)},\quad \jj^{(1)}\ge J+1.
\end{gather*}
We put $m=\jj^{(1)}-J-1,$  then
\begin{gather*}
N_{m}N_{-m-1}\le M^{(1)}M^{(3)},\quad m\ge0.
\end{gather*}
Because of \eqref{10-8}, for $ m\ge0$ we obtain
\begin{gather*}
N_{m}N_{-m-1}=N^{1-\frac{1}{2-\epsilon_0}(1-\epsilon_0)^{m}+\frac{1}{2-\epsilon_0}(1-\epsilon_0)^{m+2}}=
N^{1-\epsilon_0(1-\epsilon_0)^{m}}\ge N^{1-\epsilon_0}.
\end{gather*}
Hence, the inequality $M^{(1)}M^{(3)}\ge N^{1-\epsilon_0}$ holds, contrary to \eqref{13-43}. This completes the proof of the lemma.
\end{proof}
\end{Le}
We write
\begin{gather}\label{13-44}
j_0(M)=j_0(0,\jj)=j_0(\h,2J+1),\quad\Omega_2(M)=\Omega(\h,2J+1),\quad\Omega_1(M)=\Omega(0,\jj),
\end{gather}
where $j_0(j,h)$ is defined in \eqref{10-28} and $\Omega(j,h)$ is defined in \eqref{13-22}. Let verify that $j_0(M)$ is well-defined. Actually, since $\h=2J-\jj+1,$ so we are to prove that
$$j_0(0,\jj)=j_0(2J-\jj+1,2J+1)\,\mbox{for }2\le\jj\le2J.$$
If $\jj\le J+1,$ then
$$j_0(0,\jj)=\left|\jj-J-1\right|=J+1-\jj=\left|2J+2-\jj-J-1\right|=j_0(2J-\jj+1,2J+1).$$
If $\jj>J+1,$ then
$$j_0(0,\jj)=j_0(2J-\jj+1,2J+1)=0.$$
Hence, $j_0(M)$ is well-defined.
\begin{Th}\label{theorem13-3}
For any $M,$ satisfying the inequality \eqref{13-42}, the ensemble $\Omega_N$ can be represented in two ways может быть представлен двояким образом в виде
\begin{gather}\label{13-46-1}
\Omega_N=\Omega^{(1)}\Omega^{(2)}\quad\mbox{or}\quad\Omega_N=\Omega^{(4)}\Omega^{(3)},
\end{gather}
where
\begin{gather*}
\Omega^{(1)}=\Omega_1(M)=\Xi_{1}\Xi_{2}\ldots\Xi_{\jj},\quad
\Omega^{(2)}=\Xi_{\jj+1}\Xi_{j+2}\ldots\Xi_{2J+1}\\
\Omega^{(3)}=\Omega_2(M)=\Xi_{\h+1}\Xi_{\h+2}\ldots\Xi_{2J+1},\quad
\Omega^{(4)}=\Xi_{1}\Xi_{2}\ldots\Xi_{\h}
\end{gather*}
and for any $\gamma_1\in\Omega^{(1)},\,\gamma_2\in\Omega^{(2)},\,\gamma_3\in\Omega^{(3)}$ the following inequalities hold
\begin{gather}\label{13-47}
\frac{M}{70A^2}\le\|\gamma_1\|\le1,03M^{1+2\epsilon_0},
\end{gather}
\begin{gather}\label{13-47-1}
\frac{N}{140A^2H_1(M)}\le\|\gamma_2\|\le\frac{73A^2N}{M},\quad\mbox{where}\quad H_1(M)=1,03M^{1+2\epsilon_0},
\end{gather}
\begin{gather}\label{13-48}
\frac{M}{150A^2}\le\|\gamma_3\|\le H_3(M),\quad\mbox{where}\quad H_3(M)=80A^{2,1}M^{1+2\epsilon_0}.
\end{gather}
\begin{proof}
In view of the Corollary \ref{conseq13-3} the inequalities \eqref{13-40} and \eqref{13-41} hold for the matrices $\gamma_1=\xi_1\xi_2\ldots\xi_{\jj}$ and  $\gamma_3=\xi_{\h+1}\xi_2\ldots\xi_{2J+1}.$ Using these inequalities we obtain \eqref{13-47} and \eqref{13-48}. Moreover, it follows from \eqref{continuant inequality} and \eqref{13-16} that
\begin{gather}\label{13-48-1}
\frac{N}{140A^2}\le\frac{\|\xi_1\xi_2\ldots\xi_{2J+1}\|}{2}\le\|\gamma_1\|\|\gamma_2\|\le\|\xi_1\xi_2\ldots\xi_{2J+1}\|\le1,01N.
\end{gather}
Substituting the bound \eqref{13-47} into \eqref{13-48-1} we obtain the inequality \eqref{13-47-1}. This completes the proof of the theorem.
\end{proof}
\end{Th}
\begin{Th}\label{theorem13-3-1}
Let $M^{(2)}\ge(M^{(1)})^{2\epsilon_0},$ and let the inequality \eqref{13-42} holds for $M=M^{(1)}$ and $M=M^{(1)}M^{(2)}.$ Then the ensemble $\Omega_N$ can be represented in the following form
\begin{gather*}
\Omega_N=\Omega^{(1)}\Omega^{(2)}\Omega^{(3)},
\end{gather*}
where
\begin{gather*}
\Omega^{(1)}=\Omega_1(M^{(1)})=\Xi_{1}\Xi_{2}\ldots\Xi_{\jj_1},\quad
\Omega^{(2)}=\Xi_{\jj_1+1}\Xi_{\jj_1+2}\ldots\Xi_{\jj_2}\\
\Omega^{(1)}\Omega^{(2)}=\Omega_1(M^{(1)}M^{(2)}),\\
\Omega^{(3)}=\Xi_{\jj_2+1}\Xi_{\jj_2+2}\ldots\Xi_{2J+1},\quad
\end{gather*}
where $\jj_1$ corresponds to $M^{(1)},$ and $\jj_2$  corresponds to $M^{(1)}M^{(2)}.$ And for any $\gamma_1\in\Omega^{(1)},\,\gamma_2\in\Omega^{(2)},\,\gamma_3\in\Omega^{(3)}$ the following inequalities hold
\begin{gather}\label{13-48-3}
\frac{M^{(1)}}{70A^2}\le\|\gamma_1\|\le1,03(M^{(1)})^{1+2\epsilon_0},
\end{gather}
\begin{gather}\label{13-48-4}
\frac{M^{(1)}M^{(2)}}{70A^2}\le\|\gamma_1\gamma_2\|\le1,03(M^{(1)}M^{(2)})^{1+2\epsilon_0},
\end{gather}
\begin{gather}\label{13-48-5}
\frac{M^{(2)}}{150A^2(M^{(1)})^{2\epsilon_0}}\le\|\gamma_2\|\le73A^2\frac{(M^{(2)})^{1+2\epsilon_0}}{(M^{(1)})^{2\epsilon_0}},
\end{gather}
\begin{gather}\label{13-48-6}
\frac{N}{150A^2H_1(M^{(1)}M^{(2)})}\le\|\gamma_3\|\le\frac{73A^2N}{M^{(1)}M^{(2)}}.
\end{gather}
\begin{proof}
To ensure that the partition of the ensemble is well-defined, it is enough to verify that $\jj_2>\jj_1.$  To do this it is sufficient to prove that for any $M^{(1)}$ in the interval $N_{\jj_1-1-J}\le M^{(1)}\le N_{\jj_1-J}$ the inequality $N_{\jj_1-J}\le M^{(1)}M^{(2)}$ holds. It follows from the conditions of the theorem that
\begin{gather}\label{13-48-7}
M^{(2)}\ge(M^{(1)})^{2\epsilon_0}\ge(M^{(1)})^{\frac{1}{1-\epsilon_0}-1}\ge\frac{N_{\jj_1-J}}{M^{(1)}},
\end{gather}
where the last inequality holds because of $M^{(1)}\ge N_{\jj-J}^{1-\epsilon_0}$ see \eqref{10-21}. Thus we have proved that the partition of the ensemble is well-defined. The bounds \eqref{13-48-3} and \eqref{13-48-4} follow from \eqref{13-47} and \eqref{13-48-6} follows from \eqref{13-47-1}. Finally, \eqref{13-48-5} follows from \eqref{13-48-3},\, \eqref{13-48-4} and \eqref{continuant inequality}. This completes the proof of the theorem.
\end{proof}
\end{Th}
\begin{Th}\label{theorem13-4}
For any $M^{(1)}$ and $M^{(3)}$ satisfying the hypotheses of Lemma \ref{lemma13-6} the ensemble $\Omega_N$ can be represented in the following form
\begin{gather*}
\Omega_N=\Omega^{(1)}\Omega^{(2)}\Omega^{(3)},
\end{gather*}
where
\begin{gather*}
\Omega^{(1)}=\Omega_1(M^{(1)})=\Xi_{1}\Xi_{2}\ldots\Xi_{\jj_1},\quad
\Omega^{(2)}=\Xi_{\jj_1+1}\Xi_{\jj_1+2}\ldots\Xi_{\h_2},\\
\Omega^{(3)}=\Omega_2(M^{(2)})=\Xi_{\h_2+1}\Xi_{\h_2+2}\ldots\Xi_{2J+1},\quad
\end{gather*}
and for any $\gamma_1\in\Omega^{(1)},\,\gamma_2\in\Omega^{(2)},\,\gamma_3\in\Omega^{(3)}$ the following inequalities hold
\begin{gather}\label{13-49}
\frac{M^{(1)}}{70A^2}\le\|\gamma_1\|\le1,03(M^{(1)})^{1+2\epsilon_0},
\end{gather}
\begin{gather}\label{13-49-1}
\frac{N}{3\cdot10^4A^5(M^{(1)}M^{(3)})^{1+2\epsilon_0}}\le\|\gamma_2\|\le\frac{11000A^4N}{M^{(1)}M^{(3)}},
\end{gather}
\begin{gather}\label{13-50}
\frac{M^{(3)}}{150A^2}\le\|\gamma_3\|\le80A^{2,1}(M^{(3)})^{1+2\epsilon_0}.
\end{gather}
\begin{proof}
The existence of the partition is ensured by the inequality $\jj_1<\h_3$ in Lemma \ref{lemma13-6} ($\jj_1=j^{(1)},\,\h_3=h^{(3)}$). It also follows from Lemma \ref{lemma13-6} that the inequalities \eqref{13-40} and \eqref{13-41} hold for $\gamma_1\in\Omega^{(1)}$ and $\gamma_3\in\Omega^{(3)}.$ Using these inequalities we obtain the estimates \eqref{13-49} and \eqref{13-50}. Next, in the same way as \eqref{13-48-1} we obtain
\begin{gather}\label{13-50-1}
\frac{N}{280A^2}\le\|\gamma_1\|\|\gamma_2\|\|\gamma_3\|\le1,01N.
\end{gather}
Substituting \eqref{13-49} and \eqref{13-50} into \eqref{13-50-1}, we obtain \eqref{13-49-1}. This completes the proof of the theorem.
\end{proof}
\end{Th}
\begin{Le}\label{lemma13-7}
For $M$ in the interval \eqref{13-42} the following inequality holds
\begin{gather}\label{13-51}
\frac{\log\log M^2}{\log(1-\epsilon_0)}-1\le j_0(M)+\frac{\log\log N}{\log(1-\epsilon_0)}.
\end{gather}
\begin{proof}
We consider two cases.
\begin{enumerate}
  \item Let
$\exp\left(\frac{10^5A^4}{\epsilon_0^2}\right)\le M\le N_1=N^{\frac{1}{2-\epsilon_0}}.$\par
Then it follows from Lemma \ref{lemma10-4} and \eqref{10-23} that firstly  $2\le\jj\le J+1$ and secondly
\begin{gather*}
N_{\jj-J}^{1-\epsilon_0}\le M\le N_{\jj-J}.
\end{gather*}
Hence $j_0(M)=J+1-\jj$ and applying \eqref{10-8} we obtain
\begin{gather}\label{13-52}
M\ge N_{\jj-J}^{1-\epsilon_0}=N^{\frac{1}{2-\epsilon_0}(1-\epsilon_0)^{2-\jj+J}}=
N^{\frac{1}{2-\epsilon_0}(1-\epsilon_0)^{1+j_0(M)}}\ge N^{\frac{1}{2}(1-\epsilon_0)^{1+j_0(M)}}.
\end{gather}
Taking a logarithm twice we obtain
\begin{gather*}
\log\log M^2\ge(1+j_0(M))\log(1-\epsilon_0)+\log\log N.
\end{gather*}
From this the inequality \eqref{13-51} follows.
  \item Let
$N^{\frac{1}{2-\epsilon_0}}=N_1\le M\le N\exp\left(-\frac{10^5A^4}{\epsilon_0^2}\right).$\par
Then $\jj>J+1$ and, hence, $j_0(M)=0.$ In view of Lemma \ref{lemma10-3} and \eqref{10-8} we obtain
\begin{gather}\notag
M\ge  N_{\jj-J}^{1-\epsilon_0}=N^{\left(1-\frac{1}{2-\epsilon_0}(1-\epsilon_0)^{\jj-J}\right)(1-\epsilon_0)}\ge
N^{\left(1-\frac{1}{2-\epsilon_0}(1-\epsilon_0)^{2}\right)(1-\epsilon_0)}\\\ge
N^{\frac{1}{2}(1-\epsilon_0)^{1+j_0(M)}},\label{13-53}
\end{gather}
since $1-\frac{1}{2-\epsilon_0}(1-\epsilon_0)^{2}\ge\frac{1}{2}.$ The inequality  \eqref{13-53} coincides with the bound \eqref{13-52}. Thus we obtain \eqref{13-51} in the same way.
\end{enumerate}
The lemma is proved.
\end{proof}
\end{Le}
For $M$ in the interval \eqref{13-42} we define the following function
\begin{gather}\label{13-54}
T(M)=\exp\left(-\left(\frac{\log\log M^2}{\log(1-\epsilon_0)}-1\right)^2\right).
\end{gather}
We observe that for $\epsilon_0\in(0,\frac{1}{2500})$ one has
\begin{gather}\label{13-54-1}
M^{-\epsilon_0}\le T(M), \quad\mbox{if}\quad M\ge \exp(\epsilon_0^{-5}).
\end{gather}
\begin{Th}\label{theorem13-5}
Let $M\ge \exp(\epsilon_0^{-5})$ and belongs to the interval \eqref{13-42}. Then the following bounds on the cardinality of the sets $|\Omega_1(M)|$ and $|\Omega_2(M)|$ hold
\begin{gather}\label{13-55}
M^{2\delta-\epsilon_0}\le M^{2\delta}T(M)\le\left|\Omega_1(M)\right|\le 9M^{2\delta+4\epsilon_0},
\end{gather}
\begin{gather}\label{13-56}
\left|\Omega_2(M)\right|\ge M^{2\delta}T(M)\ge M^{2\delta-\epsilon_0}.
\end{gather}
\begin{proof}
Using the definition of the set $|\Omega_1(M)|$ and the inequality
\begin{gather*}
N_{\jj-J}\le M^{\frac{1}{1-\epsilon_0}}\le M^{1+2\epsilon_0}
\end{gather*}
we obtain that the upper bound in \eqref{13-55} follows immediately  from Lemma \ref{lemma13-3}.
\begin{enumerate}
  \item We estimate the cardinality of the set $|\Omega_1(M)|$ from below.
Taking into account that ${M\le N_{\jj-J}}$ and $j_0(M)=j_0(0,\jj),$ it follows from the Theorem \ref{theorem13-1} that
\begin{gather}\label{13-57}
\left|\Omega_1(M)\right|\ge M^{2\delta}\exp\left(-\left(\frac{\log\log N}{\log(1-\epsilon_0)}+j_0(M)\right)^2\right).
\end{gather}
Using \eqref{13-35} and \eqref{13-51}, we obtain, in view of $j_0(M)\le J,$ that
\begin{gather*}
\frac{\log\log M^2}{\log(1-\epsilon_0)}-1\le j_0(M)+\frac{\log\log N}{\log(1-\epsilon_0)}\le0.
\end{gather*}
Hence,
\begin{gather*}
T(M)\le \exp\left(-\left(\frac{\log\log N}{\log(1-\epsilon_0)}+j_0(M)\right)^2\right).
\end{gather*}
Substituting the estimate into \eqref{13-57}, we obtain the lower bound in \eqref{13-55}.
  \item We estimate the cardinality of the set $|\Omega_2(M)|$ from below.
Taking into account that ${M\le \frac{N}{N_{\h-J}}}$ and $j_0(M)=j_0(\h,2J+1),$ it follows from the Theorem \ref{theorem13-1} that
\begin{gather}\label{13-58}
\left|\Omega_2(M)\right|\ge M^{2\delta}\exp\left(-\left(\frac{\log\log N}{\log(1-\epsilon_0)}+j_0(M)\right)^2\right).
\end{gather}
From this the estimate \eqref{13-56} follows in the same way.
\end{enumerate}
The theorem is proved.
\end{proof}
\end{Th}

\part{\large{Estimates of exponential sums and integrals. A generalization of the Bourgain-Kontorovich's method.}}
\section{General estimates of exponential sums over ensemble.}\label{Trigsums1}
Recall that to prove our main theorem \ref{uslov} we should obtain the maximum accurate bound on the integral
\begin{gather}\label{14-1}
\int_0^1\left|S_N(\theta)\right|^2d\theta=\int_0^1\left|\sum_{\gamma\in\Omega_{N} }e(\theta\|\gamma\|)\right|^2d\theta,
\end{gather}
where $N$ is a sufficiently large integer and for $\|\gamma\|$
\begin{equation}\label{14-2}
\gamma=
\begin{pmatrix}
a & b \\
c & d
\end{pmatrix}=
\begin{pmatrix}
0 & 1 \\
1 & d_1
\end{pmatrix}
\begin{pmatrix}
0 & 1 \\
1 & d_2
\end{pmatrix}\ldots
\begin{pmatrix}
0 & 1 \\
1 & d_k
\end{pmatrix}
\end{equation}
the norm $\|\gamma\|$  is defined by
\begin{gather}\label{14-3}
\|\gamma\|=d=<d_1,d_2,\ldots,d_k>=(0,1)
\begin{pmatrix}
a & b \\
c & d
\end{pmatrix}
\begin{pmatrix}
0   \\
1
\end{pmatrix}.
\end{gather}
Then
\begin{gather}\label{14-4}
S_N(\theta)=\sum_{\gamma\in\Omega_{N} }e(\theta\|\gamma\|)=\sum_{\gamma\in\Omega_{N} }e((0,1)\gamma
\begin{pmatrix}
0   \\
1
\end{pmatrix}
\theta).
\end{gather}
To estimate the sum \eqref{14-4} different methods were used in ~\cite{BK} depending on the value of $\theta.$  The basis of all these (and even new one) methods can be presented in an unified manner if special notations are used. Suppose that a partition of the ensemble
\begin{gather}\label{14-5}
\Omega_N=\Omega^{(1)}\Omega^{(2)}\Omega^{(3)}
\end{gather}
is given and that this partition has one of the following forms:
\begin{itemize}
  \item or as in Theorem \ref{theorem13-4} and then
\begin{gather}\label{14-6}
\Omega^{(1)}=\Omega_1(M^{(1)}),\quad\Omega^{(3)}=\Omega_2(M^{(3)}),
\end{gather}
  \item or as in Theorem \ref{theorem13-3-1} and then
\begin{gather}\label{14-7}
\Omega^{(1)}\Omega^{(2)}=\Omega_1(M^{(1)}M^{(2)}),
\end{gather}
  \item or as in Theorem \ref{theorem13-3} $\Omega_N=\Omega^{(1)}\Omega^{(3)},$ and then
\begin{gather}\label{14-8}
\Omega^{(1)}=\Omega_1(M^{(1)}),\quad\Omega^{(2)}=\{E\},
\end{gather}
where $E$~--is a unit matrix $2\times2,$ and $M^{(1)}=M.$
\end{itemize}
We recall that
\begin{gather}\label{14-10}
H_1(M)=1,03M^{1+2\epsilon_0},\quad H_3(M)=80A^{2,1}M^{1+2\epsilon_0}
\end{gather}
and for $\gamma_1\in\Omega^{(1)}$ one has
\begin{gather}\label{14-10-1}
\|\gamma_1\|\le H_1(M^{(1)})=H_1.
\end{gather}
For $n\in\{1,2,3\}$ we write
\begin{gather}\label{14-11}
\widetilde{\Omega}^{(n)}=
\left\{
              \begin{array}{ll}
                (0,1)\Omega^{(1)}=\left\{(0,1)g_1\Bigl| g_1\in\Omega^{(1)}\right\}, & \hbox{если $n=1$,} \\
                \Omega^{(2)}, & \hbox{если $n=2,$}\\
              \Omega^{(3)}\begin{pmatrix}0  \\1 \end{pmatrix}=\left\{g_3\begin{pmatrix}0  \\1 \end{pmatrix}\Bigl| g_3\in\Omega^{(3)}\right\}, & \hbox{если $n=3$.}
                 \end{array}
\right.,
\end{gather}
Let define the following function
\begin{gather}\label{14-12}
T(x)=\max\left\{0,1-|x|\right\},\quad S(x)=\left(\frac{\sin \pi x}{\pi x}\right)^2.
\end{gather}
It is  common knowledge ~\cite[(4.83)]{IK} that $\hat{S}(x)=T(x),$ where $\hat{f}(x)=\int_{-\infty}^{\infty}f(t)e(-xt)dt$ is the Fourier transform of the function $f(x).$ We consider $S_2(x)=3S(\frac{x}{2}).$ It is obvious that $S_2(x)$ is a nonnegative function and
\begin{gather}\label{14-13}
S_2(x)>1\, \mbox{for}\, x\in[-1,1].
\end{gather}
Since $\hat{S}_2(x)=6T(2x),$ we have $\hat{S}_2(x)\neq 0$ only for $|x|<\frac{1}{2}.$  We next undertake to estimate the sum of the form
\begin{gather}\label{14-14}
\sigma_{N,Z}=\sum_{\theta\in Z}\left|S_N(\theta)\right|,
\end{gather}
where $Z$ is a finite subset of the interval $[0,1].$
\begin{Le}\label{lemma-14-1}
For either $(\mu,\lambda)=(2,3)$ or $(\mu,\lambda)=(3,2)$ the following estimate holds
\begin{gather}\label{14-15}
\sigma_{N,Z}\le\left|\Omega^{(1)}\right|^{1/2}\sum_{g_{\lambda}\in\Om^{(\lambda)}}\Biggl(
\sum_{g_1\in\Z^2}\s\left(\frac{g_1}{H_1}\right)
\left|\sum_{\theta\in Z}\xi(\theta)\sum_{g_{\mu}\in\Om^{(\mu)}}e(g_1g_2g_3\theta)\right|^2\Biggr)^{1/2},
\end{gather}
where $\s(x,y)=S_2(x)S_2(y),$ and $\xi(\theta)$ is a complex number with $|\xi(\theta)|=1.$
\begin{proof}
The numbers $\xi(\theta)$ are defined by the relation $\left|S_N(\theta)\right|=\xi(\theta)S_N(\theta).$ Then we obtain from the formulae \eqref{14-4}, \eqref{14-5} and the definition \eqref{14-11} that
\begin{gather}\label{14-16}
\sigma_{N,Z}=\sum_{\theta\in Z}\xi(\theta)S_N(\theta)=\sum_{\theta\in Z}\xi(\theta)\sum_{g_{1}\in\Om^{(1)}}
\sum_{g_{2}\in\Om^{(2)}}\sum_{g_{3}\in\Om^{(3)}}e(g_1g_2g_3\theta),
\end{gather}
where $g_1,g_3$ are already vectors in $\Z^2.$ It follows from \eqref{14-16} that
\begin{gather}\label{14-17}
\sigma_{N,Z}\le\sum_{g_{\lambda}\in\Om^{(\lambda)}}\sum\limits_{g_1\in\Z^2\atop |g_1|\le H_1}\1_{g_{1}\in\Om^{(1)}}\left|
\sum_{\theta\in Z}\xi(\theta)\sum_{g_{\mu}\in\Om^{(\mu)}}e(g_1g_2g_3\theta)\right|.
\end{gather}
We note that in view of \eqref{14-10-1} the condition $|g_1|\le H_1$ does not impose any extra restrictions. Applying the Cauchy~-Schwarz inequality in \eqref{14-17} we obtain
\begin{gather}\label{14-18}
\sigma_{N,Z}\le\sum_{g_{\lambda}\in\Om^{(\lambda)}}
\left(\sum\limits_{g_1\in\Z^2\atop |g_1|\le H_1}\1_{g_{1}\in\Om^{(1)}}\right)^{1/2}
\left(\sum\limits_{g_1\in\Z^2\atop |g_1|\le H_1}
\left|\sum_{\theta\in Z}\xi(\theta)\sum_{g_{\mu}\in\Om^{(\mu)}}e(g_1g_2g_3\theta)\right|^2\right)^{1/2}.
\end{gather}
Considering that
$$\left(\sum\limits_{g_1\in\Z^2\atop |g_1|\le H_1}\1_{g_{1}\in\Om^{(1)}}\right)^{1/2}=\left|\Omega^{(1)}\right|^{1/2}$$
and taking into account that the function $\s(x,y)>1$ for {$(x,y)\in[-1,1]^2$} and is nonnegative, we obtain \eqref{14-15}. This completes the proof of the lemma.
\end{proof}
\end{Le}
To reduce our notations we in what follows will identify any variable $x$ with $x^{(1)}.$
\begin{Le}\label{lemma-14-2}
Under the hypotheses of Lemma \ref{lemma-14-1} the following bound holds
\begin{gather}\label{14-19}
\sigma_{N,Z}\le
10H_1\left|\Omega^{(1)}\right|^{1/2}\sum_{g_{\lambda}\in\Om^{(\lambda)}}
\Biggl(
\sum\limits_{g^{(1)}_{\mu},g^{(2)}_{\mu}\in\Om^{(\mu)}\atop\theta^{(1)},\theta^{(2)}\in Z}
\1_{\{\|z\|_{1,2}\le \frac{1}{2H_1}\}}
\Biggr)^{1/2},
\end{gather}
where $\|x\|_{1,2}=\max\{\|x_1\|,\|x_2\|\}$ for $x=(x_1,x_2)\in\rr^2,$  and
\begin{gather}\label{14-20}
z=
\left\{
              \begin{array}{ll}
              g_2^{(1)}g_3\theta^{(1)}-g_2^{(2)}g_3\theta^{(2)} & \hbox{if $\mu=2$,} \\
              g_2g_3^{(1)}\theta^{(1)}-g_2g_3^{(2)}\theta^{(2)}, & \hbox{if $\mu=3$.}
                 \end{array}
\right..
\end{gather}
\begin{proof}
We note that the quantity $z$ can be represented in a shorter form
\begin{gather}\label{14-21}
z=g_2^{(1)}g_3^{(1)}\theta^{(1)}-g_2^{(4-\mu)}g_3^{(\mu-1)}\theta^{(2)}
\end{gather}
Applying the relation $|x|^2=x\overline{x}$ and reversing orders we easily obtain that
\begin{gather}\label{14-22}
\sum_{g_1\in\Z^2}\s\left(\frac{g_1}{H_1}\right)
\left|\sum_{\theta\in Z}\xi(\theta)\sum_{g_{\mu}\in\Om^{(\mu)}}e(g_1g_2g_3\theta)\right|^2\le
\sum_{g^{(1)}_{\mu},g^{(2)}_{\mu}\in\Om^{(\mu)}}\sum_{\theta^{(1)},\theta^{(2)}\in Z}
\left|\sum_{g_1\in\Z^2}\s\left(\frac{g_1}{H_1}\right)e(g_1z)\right|.
\end{gather}
By application of the Poisson summation formula ~\cite[§4.3.]{IK}:
\begin{gather*}
\sum_{n\in\Z^2}f(n)=\sum_{k\in\Z^2}\hat{f}(k),
\end{gather*}
and writing $f(n)=\s\left(\frac{n}{H_1}\right)e(nz),$ we transform the inner sum in the right side of \eqref{14-22}:
\begin{gather}\label{14-23}
\sum_{g_1\in\Z^2}\s\left(\frac{g_1}{H_1}\right)e(g_1z)=
\sum_{k\in\Z^2}\int_{x\in\rr^2}\s\left(\frac{x}{H_1}\right)e(x(z-k))dx=H_1^2\sum_{k\in\Z^2}\hat{\s}((k-z)H_1).
\end{gather}
We note that the relation \eqref{14-23} can be obtained directly from ~\cite[(4.25)]{IK}. As $\hat{\s}(x,y)\neq0$ only if  $|x|\le\frac{1}{2}$ and $|y|\le\frac{1}{2},$ so the sum in the right side of \eqref{14-23} consists of at most one summand, hence,
\begin{gather}\label{14-24}
\left|\sum_{g_1\in\Z^2}\s\left(\frac{g_1}{H_1}\right)e(g_1z)\right|\le
36H_1^2\1_{\{\|z\|_{1,2}\le \frac{1}{2H_1}\}}.
\end{gather}
Substituting \eqref{14-24} into \eqref{14-22}, we obtain
\begin{gather}\label{14-25}
\sum_{g_1\in\Z^2}\s\left(\frac{g_1}{H_1}\right)
\left|\sum_{\theta\in Z}\xi(\theta)\sum_{g_{\mu}\in\Om^{(\mu)}}e(g_1g_2g_3\theta)\right|^2\le
36H_1^2
\sum\limits_{g^{(1)}_{\mu},g^{(2)}_{\mu}\in\Om^{(\mu)}\atop\theta^{(1)},\theta^{(2)}\in Z}
\1_{\{\|z\|_{1,2}\le \frac{1}{2H_1}\}}
\end{gather}
Substituting \eqref{14-25} into \eqref{14-15}, we obtain \eqref{14-19}. This completes the proof of the lemma.
\end{proof}
\end{Le}
To transform the right side of \eqref{14-19}, we are to specify the set $Z.$ It follows from the Dirichlet theorem that for any
$\theta\in[0,1]$ there exist $a,q\in\N\cup\{0\}$ and $\beta\in\rr,$ such that
\begin{gather}\label{14-26}
\theta=\frac{a}{q}+\beta,\;(a,q)=1,\; 0\le a\le q\le N^{1/2},\;\beta=\frac{K}{N},\; |K|\le\frac{N^{1/2}}{q},
\end{gather}
and $a=0$ or $a=q$ only if $q=1.$
We denote
\begin{gather}\label{14-27}
P_{Q_1,Q}^{(\beta)}=\left\{
\theta=\frac{a}{q}+\beta\;\Bigl|\;(a,q)=1,\;  0\le a\le q,\;Q_1\le q\le Q
\right\}.
\end{gather}
In what follows we always have $Z\subseteq P_{Q_1,Q}^{(\beta)}$ for some $Q_1,Q,\beta.$ We will write numbers $\theta^{(1)},\theta^{(2)}\in P_{Q_1,Q}^{(\beta)}$ in the following way
\begin{gather}\label{14-28}
\theta^{(1)}=\frac{a^{(1)}}{q^{(1)}}+\beta,\quad \theta^{(2)}=\frac{a^{(2)}}{q^{(2)}}+\beta.
\end{gather}
Let
\begin{gather}\label{14-29}
\NN_0(g_{\lambda})=\left\{ (g^{(1)}_{\mu},g^{(2)}_{\mu},\theta^{(1)},\theta^{(2)})\in
\Om^{(\mu)}\times\Om^{(\mu)}\times Z^2\Bigl|\,
\left\|g_2g_3\theta^{(1)}-g_2^{(4-\mu)}g_3^{(\mu-1)}\theta^{(2)}\right\|_{1,2}\le \frac{1}{2H_1},
\right\}
\end{gather}
\begin{gather}\label{14-30}
\NN(g_{\lambda})=\left\{ (g^{(1)}_{\mu},g^{(2)}_{\mu},\theta^{(1)},\theta^{(2)})\in
\Om^{(\mu)}\times\Om^{(\mu)}\times Z^2\Bigl|\,
\eqref{14-31}\, \mbox{и}\, \eqref{14-32}\, \mbox{hold}
\right\},
\end{gather}
where
\begin{gather}\label{14-31}
\|g_2g_3\frac{a^{(1)}}{q^{(1)}}-g_2^{(4-\mu)}g_3^{(\mu-1)}\frac{a^{(2)}}{q^{(2)}}\|_{1,2}\le \frac{74A^2\KK}{M^{(1)}},
\end{gather}
\begin{gather}\label{14-32}
|g_2g_3-g_2^{(4-\mu)}g_3^{(\mu-1)}|_{1,2}\le\min\left\{
\frac{73A^2N}{M^{(1)}},\; \frac{73A^2N}{M^{(1)}\overline{K}}+\frac{N}{\overline{K}}
\left\|g_2g_3\frac{a^{(1)}}{q^{(1)}}-g_2^{(4-\mu)}g_3^{(\mu-1)}\frac{a^{(2)}}{q^{(2)}}\right\|_{1,2}
\right\},
\end{gather}
and $\KK=\max\{1,|K|\}.$ We note that it follows from Lemma \ref{lemma-14-2} and the definition \eqref{14-29} that
\begin{gather}\label{14-33}
\sigma_{N,Z}\le
10H_1\left|\Omega^{(1)}\right|^{1/2}\sum_{g_{\lambda}\in\Om^{(\lambda)}}
\left|\NN_0(g_{\lambda})\right|^{1/2}.
\end{gather}
\begin{Le}\label{lemma-14-3}
For $Z\subseteq P_{Q_1,Q}^{(\beta)}$ and for $M^{(1)}$ in \eqref{13-42} such that $M^{(1)}>146A^2\KK,$ the following inequality holds
\begin{gather}\label{14-34}
\sigma_{N,Z}\le
10H_1\left|\Omega^{(1)}\right|^{1/2}\sum_{g_{\lambda}\in\Om^{(\lambda)}}
\left|\NN(g_{\lambda})\right|^{1/2}.
\end{gather}
\begin{proof}
In view of \eqref{14-33}, it is sufficient to prove that $\NN_0(g_{\lambda})\subseteq \NN(g_{\lambda}).$ Let $\Omega=\Omega^{(2)}\Omega^{(3)},$ then $\Omega_N=\Omega^{(1)}\Omega.$  For any partition \eqref{14-6},\,\eqref{14-7},\,\eqref{14-8} of the ensemble $\Omega_N$ using \eqref{continuant inequality} and the lower bound on $\|\gamma_1\|$ from the theorem corresponding to the partition, we obtain
\begin{gather}\label{14-35}
|g_2g_3-g_2^{(4-\mu)}g_3^{(\mu-1)}|_{1,2}\le\max\limits_{\gamma\in\Omega}\|\gamma\|\le\frac{73A^2N}{M^{(1)}}.
\end{gather}
Hence,
\begin{gather}\label{14-36}
|(g_2g_3-g_2^{(4-\mu)}g_3^{(\mu-1)})\beta|_{1,2}\le\frac{73A^2N}{M^{(1)}}\frac{\KK}{N}\le\frac{1}{2},
\end{gather}
so
\begin{gather}\label{14-37}
|(g_2g_3-g_2^{(4-\mu)}g_3^{(\mu-1)})\beta|_{1,2}=\|(g_2g_3-g_2^{(4-\mu)}g_3^{(\mu-1)})\beta\|_{1,2},
\end{gather}
and
\begin{gather}\label{14-38}
\|(g_2g_3-g_2^{(4-\mu)}g_3^{(\mu-1)})\beta\|_{1,2}\le\frac{73A^2\KK}{M^{(1)}}.
\end{gather}
It follows from the definition \eqref{14-29} and the bound \eqref{14-38} that
\begin{gather}\label{14-39}
\|g_2g_3\frac{a^{(1)}}{q^{(1)}}-g_2^{(4-\mu)}g_3^{(\mu-1)}\frac{a^{(2)}}{q^{(2)}}\|_{1,2}\le \frac{73A^2\KK}{M^{(1)}}+\frac{1}{M^{(1)}}<\frac{74A^2\KK}{M^{(1)}},
\end{gather}
that is, for $(g^{(1)}_{\mu},g^{(2)}_{\mu},\theta^{(1)},\theta^{(2)})\in\NN_0(g_{\lambda})$ the inequality \eqref{14-31} holds. Using \eqref{14-37} and \eqref{14-29}, we obtain
\begin{gather}\label{14-40}
|(g_2g_3-g_2^{(4-\mu)}g_3^{(\mu-1)})\frac{K}{N}|_{1,2}\le\frac{1}{H_1}+
\|g_2g_3\frac{a^{(1)}}{q^{(1)}}-g_2^{(4-\mu)}g_3^{(\mu-1)}\frac{a^{(2)}}{q^{(2)}}\|_{1,2},
\end{gather}
whence
\begin{gather}\label{14-41}
|g_2g_3-g_2^{(4-\mu)}g_3^{(\mu-1)}|_{1,2}\le
\frac{N}{M^{(1)}\K}+\frac{N}{\K}
\left\|g_2g_3\frac{a^{(1)}}{q^{(1)}}-g_2^{(4-\mu)}g_3^{(\mu-1)}\frac{a^{(2)}}{q^{(2)}}\right\|_{1,2}.
\end{gather}
The inequality  \eqref{14-32} follows from the estimates \eqref{14-41} and \eqref{14-35}. This completes the proof of the lemma.
\end{proof}
\end{Le}
We denote
\begin{gather}\label{14-42}
\M(g_{\lambda})=\left\{ (g^{(1)}_{\mu},g^{(2)}_{\mu},\theta^{(1)},\theta^{(2)})\in
\Om^{(\mu)}\times\Om^{(\mu)}\times Z^2\Bigl|\,
\eqref{14-43}\, \mbox{и}\, \eqref{14-44}\, \mbox{hold}
\right\},
\end{gather}
where
\begin{gather}\label{14-43}
|g_2g_3-g_2^{(4-\mu)}g_3^{(\mu-1)}|_{1,2}\le \frac{73A^2N}{M^{(1)}\KK},
\end{gather}
\begin{gather}\label{14-44}
\|g_2g_3\frac{a^{(1)}}{q^{(1)}}-g_2^{(4-\mu)}g_3^{(\mu-1)}\frac{a^{(2)}}{q^{(2)}}\|_{1,2}=0.
\end{gather}
\begin{Le}\label{lemma-14-4}
Let the hypotheses of one of the Theorem \ref{theorem13-4}, \ref{theorem13-3-1} or \ref{theorem13-3}, on which the partition of $\Omega_N$ in the form \eqref{14-5} is based, hold. Let $M^{(1)}$ be such that for any $\theta^{(1)},\theta^{(2)}\in Z$ the following inequality holds
\begin{gather}\label{14-45}
[q^{(1)},q^{(2)}]<\frac{M^{(1)}}{74A^2\KK}.
\end{gather}
Then the following bound holds
\begin{gather}\label{14-46}
\sigma_{N,Z}\le
10H_1\left|\Omega^{(1)}\right|^{1/2}\sum_{g_{\lambda}\in\Om^{(\lambda)}}
\left|\M(g_{\lambda})\right|^{1/2}.
\end{gather}
\begin{proof}
In view of \eqref{14-34} it is sufficient to prove that $\NN(g_{\lambda})\subseteq \M(g_{\lambda}).$ We note that to prove this it is sufficient to obtain that under the hypotheses of Lemma \ref{lemma-14-4} and $(g^{(1)}_{\mu},g^{(2)}_{\mu},\theta^{(1)},\theta^{(2)})\in\NN_0(g_{\lambda})$ the relation \eqref{14-44} holds. It follows from \eqref{14-31} and \eqref{14-45} that
\begin{gather}\label{14-47}
\|g_2g_3\frac{a^{(1)}}{q^{(1)}}-g_2^{(4-\mu)}g_3^{(\mu-1)}\frac{a^{(2)}}{q^{(2)}}\|_{1,2}\le \frac{74A^2\KK}{M^{(1)}}<
\frac{1}{[q^{(1)},q^{(2)}]},
\end{gather}
this implies that \eqref{14-44}. This completes the proof of the lemma.
\end{proof}
\end{Le}
Thus we reduced the problem of estimating $\sigma_{N,Z}$ to the evaluation the cardinality of one of the sets $\NN_0(g_{\lambda}),\,\NN(g_{\lambda}),\,\M(g_{\lambda}).$ The choice of the set will depend on $\mu.$ Let state one more lemma of a general nature. A similar statement was used by S.V.\,Konyagin in ~\cite[ 17]{Konyagin}.
\begin{Le}\label{lemma-14-5}
Let $W$ be a finite subset of the interval $[0,1]$ and let $|W|>10.$ Let $f:W\rightarrow \rr_{+}$ be a function such that, for any subset $Z\subseteq W$ the following bound holds
\begin{gather*}
\sum_{\theta\in Z}f(\theta)\le C_1|Z|^{1/2}+C_2,
\end{gather*}
where $C_1,C_2$ are non-negative constants not depending on the set $Z.$ Then the following estimate holds
\begin{gather}\label{14-49}
\sum_{\theta\in W}f^2(\theta)\ll C_1^2\log|W|+C_2\max_{\theta\in W}f(\theta)
\end{gather}
with the absolute constant in Vinogradov symbol.
\begin{proof}
Let number the value of $f(\theta)$ in the decreasing order
\begin{gather*}
f_1\ge\ldots\ge f_{|W|}>0.
\end{gather*}
If $C_1\ge C_2,$ in particular, if $C_2=0,$ we have $\sum_{\theta\in Z}f(n)\le 2C_1|Z|^{1/2}.$ Then for any $k$ such that $1\le k\le|W|,$ one has
\begin{gather*}
kf_k\le\sum_{n=1}^{k}f_n\le 2C_1k^{1/2}
\end{gather*}
and, hence, $f_k\le 2C_1k^{-1/2}.$ Thus
\begin{gather*}
\sum_{\theta\in W}f^2(\theta)=\sum_{n=1}^{|W|}f_n^2\le 8C_1^2\log|W|.
\end{gather*}
Therefore, it remains to consider the case $C_1<C_2.$ In a similar manner we obtain that
\begin{gather*}
kf_k\le\sum_{n=1}^{k}f_n\le C_1k^{1/2}+C_2\Rightarrow
f_k\le\frac{C_1}{k^{1/2}}+\frac{C_2}{k}.
\end{gather*}
Let
\begin{gather*}
M=\max_{\theta\in W}f(n),\, L=\frac{C_2^{4/3}}{(C_1M)^{2/3}}.
\end{gather*}
We consider three cases.
\begin{enumerate}
  \item If $L<10,$ then
\begin{gather}\label{14-50}
\sum_{\theta\in W}f^2(\theta)=\sum_{k\le|W|}f_k^2\ll\sum_{k\le|W|}\left(\frac{C^2_1}{k}+\frac{C^2_2}{k^2}\right)\ll
C_1^2\log|W|+C_2^2.
\end{gather}
It follows from the condition $L<10$ that $C_2^2\ll C_1M.$ Substituting this bound into \eqref{14-50} and taking into account $C_1<C_2,$  we obtain
\begin{gather}\label{14-51}
\sum_{\theta\in W}f^2(\theta)\ll
C_1^2\log|W|+C_2M.
\end{gather}
  \item If $10\le L<|W|,$ then
\begin{gather}\notag
\sum_{\theta\in W}f^2(\theta)=\sum_{k\le L}f_k^2+\sum_{L<k\le|W|}f_k^2\ll\max_{\theta\in W}f(n)\sum_{k\le L}f_k+
\sum_{L<k\le|W|}\left(\frac{C^2_1}{k}+\frac{C^2_2}{k^2}\right)\ll\\\ll
(C_1L^{1/2}+C_2)\max_{\theta\in W}f(\theta)+
C_1^2\log\frac{|W|}{L}+\frac{C_2^2}{L}.\label{14-52}
\end{gather}
It follows from the definition of $L$ that $C_1L^{1/2}M=\frac{C_2^2}{L}.$ Then by \eqref{14-52} we obtain
\begin{gather}\label{14-53}
\sum_{\theta\in W}f^2(\theta)\ll
C_1^2\log|W|+C_2M+(C_1C_2M)^{2/3}<C_1^2\log|W|+C_2M.
\end{gather}
  \item If $L\ge|W|,$ then
\begin{gather}\label{14-54}
\sum_{\theta\in W}f^2(\theta)\ll\max_{\theta\in W}f(n)\sum_{k\le |W|}f_k\le
(C_1|W|^{1/2}+C_2)M.
\end{gather}
It follows from $L\ge|W|$ that $|W|^{1/2}\le\frac{C_2^{2/3}}{(C_1M)^{1/3}}.$ Substituting this bound into \eqref{14-54}, we obtain
\begin{gather}\label{14-54}
\sum_{\theta\in W}f^2(\theta)\ll (C_1C_2M)^{2/3}+C_2M<C_1^2\log|W|+C_2M.
\end{gather}
\end{enumerate}
Thus, we have proved the desired formula \eqref{14-49} for all values of $L.$ This completes the proof of the lemma.
\end{proof}
\end{Le}

\section{Dedekind sums.}
The main result of this section will be used in §\ref{mu=3}. Let define the function $\varrho(x)$ in a following way $\varrho(x)=\frac{1}{2}-\{x\}$ for $x\in\rr/\Z$ and $\varrho(x)=0$ for $x\in\Z.$ This section is devoted to the estimates of generalized Dedekind sums of the following form
\begin{gather}\label{ded-0}
\sum_{0<n\le P}\varrho(y_{1}\frac{n}{P}+\frac{1}{R})\varrho(y_{2}\frac{n}{P}+\frac{1}{R}).
\end{gather}
The proof of the following statement is based on the Knuth's article ~\cite{DE-Knut}.
\begin{Le}\label{lemma-ded}
let $(y_1,y_2)=1,$ $P$ be a natural number, $R$ be a real number and $y_1,y_2<\frac{R}{10}.$  Then
\begin{gather}\label{ded-1}
\sum_{0<n\le P}
\1_{\{\|y_{1}\frac{n}{P}\|<\frac{1}{R}\}}\1_{\{\|y_{2}\frac{n}{P}\|<\frac{1}{R}\}}\le
4\frac{(y_1,P)+(y_2,P)}{R}+\frac{2P}{R}\min\{\frac{1}{y_1},\frac{1}{y_2}\}+\frac{4P}{R^2}+O(1).
\end{gather}
\begin{proof}
For any fixed real number $a$ in the interval $0<a<\frac{1}{2}$ one has
\begin{gather*}
\1_{\{\|x\|<a\}}=2a+\rho(x+a)-\varrho(x-a),\, \mbox{если}\, x\neq\pm a\\
\1_{\{\|x\|<a\}}<2a+\rho(x+a)-\varrho(x-a),\, \mbox{если}\, x=\pm a,
\end{gather*}
and, therefore,
\begin{gather}\label{15-40-1}
\1_{\{\|y_{1}\frac{n}{P}\|<\frac{1}{R}\}}\1_{\{\|y_{2}\frac{n}{P}\|<\frac{1}{R}\}}\le
\prod_{i=1,2}\left(\frac{2}{R}+\varrho(y_{i}\frac{n}{P}+\frac{1}{R})-\varrho(y_{i}\frac{n}{P}-\frac{1}{R})\right).
\end{gather}
Let $\varrho(y_{i}\frac{n}{P}\pm\frac{1}{R})=\varrho_{i}^{\pm},$ then we obtain
\begin{gather}\notag
\sum_{0<n\le P}\1_{\{\|y_{1}\frac{n}{P}\|<\frac{1}{R}\}}\1_{\{\|y_{2}\frac{n}{P}\|<\frac{1}{R}\}}\le
\frac{4P}{R^2}+\frac{2}{R}\sum_{0<n\le P}\left(\varrho_{1}^{+}-\varrho_{1}^{-}+\varrho_{2}^{+}-\varrho_{2}^{-}\right)+\\+
\sum_{0<n\le P}\left(\varrho_{1}^{+}\varrho_{2}^{+}+\varrho_{1}^{-}\varrho_{2}^{-}
-\varrho_{1}^{-}\varrho_{2}^{+}-\varrho_{1}^{+}\varrho_{2}^{-}\right)=
\frac{4P}{R^2}+\Sigma_1+\Sigma_2.
\label{15-40-2}
\end{gather}
To evaluate $\Sigma_1$ we need the following well known result.  For {$(p,q)=1$} and for any real $x$ one has
\begin{gather}\label{15-40-4}
\sum_{n=1}^{q}\varrho(\frac{p}{q}n+x)=\varrho(qx)
\end{gather}
see ~\cite[стр.170, лемма 483]{Landau}. Let estimate one of the four summands in the sum $\Sigma_1.$ Using \eqref{15-40-4}, we obtain
\begin{gather}\label{15-40-5}
\sum_{0<n\le P}\varrho(y_{1}\frac{n}{P}+\frac{1}{R})=(y_1,P)\varrho(\frac{P}{(y_1,P)}\frac{1}{R}).
\end{gather}
Substituting \eqref{15-40-5} into the definition of $\Sigma_1,$  we obtain
\begin{gather}\label{15-40-6}
|\Sigma_1|\le4\frac{(y_1,P)+(y_2,P)}{R}.
\end{gather}
To obtain an estimate on the sum $\Sigma_2$ we transform in the same manner each of the four summands of it. Consider the first summand
\begin{gather*}
\sum_{0<n\le P}\varrho(y_{1}\frac{n}{P}+\frac{1}{R})\varrho(y_{2}\frac{n}{P}+\frac{1}{R}).
\end{gather*}
\begin{description}
  \item[Step 1] We denote
\begin{gather*}
\delta_1=(y_1,P),\,\delta_2=(y_2,P),y_3=\frac{y_1}{\delta_1},\,y_4=\frac{y_2}{\delta_2},\,
P_1=\frac{P}{\delta_1},\,P_2=\frac{P}{\delta_1\delta_2}
\end{gather*}
and prove that
\begin{gather}\label{15-40-9}
\sum_{0<n\le P}\varrho(y_{1}\frac{n}{P}+\frac{1}{R})\varrho(y_{2}\frac{n}{P}+\frac{1}{R})=
\sum_{0<n\le P_2}\varrho(y_{3}\frac{n}{P_2}+\frac{\delta_2}{R})\varrho(y_{4}\frac{n}{P_2}+\frac{\delta_1}{R}).
\end{gather}
Actually, the change of variables $n=(k-1)P_1+m$ leads to
\begin{gather*}
\sum_{0<n\le P}\varrho(y_{1}\frac{n}{P}+\frac{1}{R})\varrho(y_{2}\frac{n}{P}+\frac{1}{R})=
\sum_{1\le k\le\delta_1}\sum_{1\le m\le P_1}\varrho(y_{1}\frac{k-1}{\delta_1}+\frac{y_1}{\delta_1}\frac{m}{P_1}+\frac{1}{R})
\varrho(y_{2}\frac{k-1}{\delta_1}+\frac{my_2}{P}+\frac{1}{R})=\\=
\sum_{1\le m\le P_1}\varrho(y_3\frac{m}{P_1}+\frac{1}{R})\sum_{1\le k\le\delta_1}
\varrho(y_{2}\frac{k}{\delta_1}+\frac{y_2(m-P_1)}{P}+\frac{1}{R}).
\end{gather*}
Applying formula \eqref{15-40-4} to the sum over $k$ we obtain
\begin{gather*}
\sum_{1\le k\le\delta_1}\varrho(y_{2}\frac{k}{\delta_1}+\frac{y_2(m-P_1)}{P}+\frac{1}{R})=
\varrho(\frac{y_2(m-P_1)}{P_1}+\frac{\delta_1}{R})=\varrho(\frac{y_2m}{P_1}+\frac{\delta_1}{R}),
\end{gather*}
and, hence,
\begin{gather}\label{15-40-10}
\sum_{0<n\le P}\varrho(y_{1}\frac{n}{P}+\frac{1}{R})\varrho(y_{2}\frac{n}{P}+\frac{1}{R})=
\sum_{1\le m\le P_1}\varrho(y_3\frac{m}{P_1}+\frac{1}{R})\varrho(\frac{y_2m}{P_1}+\frac{\delta_1}{R}).
\end{gather}
We note that $(y_3,P_1)=1,$ and as $(y_1,y_2)=1$ so $(y_2,P_1)=(y_2,P)=\delta_2.$ Changing the variables $m=(k-1)P_2+n$ and repeating the proof of \eqref{15-40-10}, we obtain
\begin{gather}\label{15-40-11}
\sum_{1\le m\le P_1}\varrho(y_3\frac{m}{P_1}+\frac{1}{R})\varrho(\frac{y_2m}{P_1}+\frac{\delta_1}{R})=
\sum_{0<n\le P_2}\varrho(y_{3}\frac{n}{P_2}+\frac{\delta_2}{R})\varrho(y_{4}\frac{n}{P_2}+\frac{\delta_1}{R}).
\end{gather}
The relation \eqref{15-40-9} follows from \eqref{15-40-10} and \eqref{15-40-11}.
  \item[Step 2] We find $x,y$ such that
\begin{gather}\label{15-40-12}
\frac{\delta_2}{R}=\frac{x}{P_2}+\theta_2,\,
\frac{\delta_1}{R}=\frac{y}{P_2}+\theta_1,\,
0\le\theta_1,\theta_2<\frac{1}{P_2},
\end{gather}
and prove that
\begin{gather}\label{15-40-13}
\sum_{0<n\le P_2}\varrho(\frac{ny_{3}}{P_2}+\frac{\delta_2}{R})\varrho(\frac{ny_{4}}{P_2}+\frac{\delta_1}{R})=
\sum_{0<n\le P_2}\varrho(\frac{ny_{3}}{P_2}+\frac{x}{P_2})\varrho(\frac{ny_{4}}{P_2}+\frac{y}{P_2})+O(1).
\end{gather}
Actually, in view of
\begin{gather*}
\varrho(\frac{a}{P_2})-\varrho(\frac{a}{P_2}+\theta)=
\left\{
              \begin{array}{ll}
                -\varrho(\theta), & \hbox{if $a\equiv 0 \pmod{P_2}$;} \\
                \theta, & \hbox{else,}
              \end{array}
\right.
\end{gather*}
we have
\begin{gather}\notag
\sum_{0<n\le P_2}\varrho(\frac{ny_{3}}{P_2}+\frac{x}{P_2})\varrho(\frac{ny_{4}}{P_2}+\frac{\delta_1}{R})-
\sum_{0<n\le P_2}\varrho(\frac{ny_{3}}{P_2}+\frac{\delta_2}{R})\varrho(\frac{ny_{4}}{P_2}+\frac{\delta_1}{R})=\\=
\theta_2\sum_{0<n\le P_2}\varrho(\frac{ny_{4}}{P_2}+\frac{\delta_1}{R})+O(1).\label{15-40-14}
\end{gather}
Evaluating by the formula \eqref{15-40-4} the sum in the right side of \eqref{15-40-14} we obtain
\begin{gather}\label{15-40-15}
\sum_{0<n\le P_2}\varrho(\frac{ny_{3}}{P_2}+\frac{\delta_2}{R})\varrho(\frac{ny_{4}}{P_2}+\frac{\delta_1}{R})=
\sum_{0<n\le P_2}\varrho(\frac{ny_{3}}{P_2}+\frac{x}{P_2})\varrho(\frac{ny_{4}}{P_2}+\frac{\delta_1}{R})+O(1).
\end{gather}
The similar transformations of the right side of \eqref{15-40-15} lead to the right side of \eqref{15-40-13}.
  \item[Step 3]
We make the change of variables $m\equiv ny_3+x\pmod{P_2}$ in the right side of \eqref{15-40-13}. Since $(y_3,P_2)=1,$ then $y_3^{-1}$ is defined modulo $P_2.$ Hence, $n\equiv my_3^{-1}-xy_3^{-1}\pmod{P_2}$ and
\begin{gather}\label{15-40-16}
\sum_{0<n\le P_2}\varrho(\frac{ny_{3}}{P_2}+\frac{x}{P_2})\varrho(\frac{ny_{4}}{P_2}+\frac{y}{P_2})=
\sum_{0<m\le P_2}\varrho(\frac{m}{P_2})\varrho(\frac{cm+z}{P_2}),
\end{gather}
where $c\equiv y_4y_3^{-1}\pmod{P_2},\,z=y-cx.$
\end{description}
Let
\begin{gather*}
V(z)=\sum_{0<m\le P_2}\varrho(\frac{m}{P_2})\varrho(\frac{cm+z}{P_2}),
\end{gather*}
then we have proved that
\begin{gather*}
\sum_{0<n\le P}\varrho(y_{1}\frac{n}{P}+\frac{1}{R})\varrho(y_{2}\frac{n}{P}+\frac{1}{R})=V(y-cx)+O(1),
\end{gather*}
where $c\equiv y_4y_3^{-1}\pmod{P_2}.$ Therefore,
\begin{gather*}
\sum_{0<n\le P}\left(\varrho_{1}^{+}\varrho_{2}^{+}+\varrho_{1}^{-}\varrho_{2}^{-}
-\varrho_{1}^{-}\varrho_{2}^{+}-\varrho_{1}^{+}\varrho_{2}^{-}\right)=
V(y-cx)+V(-y+cx)-V(y+cx)-V(-y-cx)+O(1).
\end{gather*}
It is proved in ~\cite[лемма 2]{DE-Knut} that $V(z)=V(-z),$ thus
\begin{gather}\label{15-40-17}
\sum_{0<n\le P}\left(\varrho_{1}^{+}\varrho_{2}^{+}+\varrho_{1}^{-}\varrho_{2}^{-}
-\varrho_{1}^{-}\varrho_{2}^{+}-\varrho_{1}^{+}\varrho_{2}^{-}\right)=
2V(cx-y)-2V(cx+y)+O(1).
\end{gather}
We note that by the symmetry of the left side of \eqref{15-40-16} with respect to $x$ and $y,$ one can assume, without loss of generality, that $x\ge y, 1\le c<P_2.$ If one of the numbers $x,y$ is equal to zero then we obtain by \eqref{15-40-17} that
\begin{gather*}
\sum_{0<n\le P}\left(\varrho_{1}^{+}\varrho_{2}^{+}+\varrho_{1}^{-}\varrho_{2}^{-}
-\varrho_{1}^{-}\varrho_{2}^{+}-\varrho_{1}^{+}\varrho_{2}^{-}\right)=O(1).
\end{gather*}
So, further we assume $x\ge y\ge1.$
\begin{description}
  \item[Step 4]
For $z>0$ we prove that
\begin{gather}\label{15-40-18}
V(z)=V(0)-\sum_{j=1}^{z}\varrho(\frac{c^{-1}j}{P_2})+O(1),
\end{gather}
where $cc^{-1}\equiv 1\pmod{P_2}.$
Since for an integer $a$ one has
\begin{gather*}
\varrho(\frac{a}{P_2})-\varrho(\frac{a-1}{P_2})=
\left\{
              \begin{array}{ll}
              \frac{1}{2}-\frac{1}{P_2}, & \hbox{if $a\equiv 0 \pmod{P_2}$;} \\
              \frac{1}{2}-\frac{1}{P_2}, & \hbox{if $a\equiv 1 \pmod{P_2}$;} \\
              -\frac{1}{P_2} , & \hbox{else,}
              \end{array}
\right.
\end{gather*}
we have
\begin{gather}\notag
V(z)-V(z-1)=\sum_{0<m\le P_2}\varrho(\frac{m}{P_2})\left(\varrho(\frac{cm+z}{P_2})-\varrho(\frac{cm+z-1}{P_2})\right)=\\=
-\frac{1}{P_2}\sum_{0<m\le P_2}\varrho(\frac{m}{P_2})+
\frac{1}{2}\left(\varrho(\frac{m_1}{P_2})+\varrho(\frac{m_2}{P_2})\right),\label{15-40-19}
\end{gather}
where $cm_1+z\equiv 0 \pmod{P_2}$ и $cm_2+z-1\equiv 0 \pmod{P_2}.$ Using the formula \eqref{15-40-4}, we obtain
\begin{gather}\label{15-40-20}
V(z)-V(z-1)=
-\frac{1}{2}\left(\varrho(\frac{c^{-1}z}{P_2})+\varrho(\frac{c^{-1}(z-1)}{P_2})\right).
\end{gather}
In a similar manner we obtain
\begin{gather}\notag
V(z-1)-V(z-2)=-\frac{1}{2}\left(\varrho(\frac{c^{-1}(z-1)}{P_2})+\varrho(\frac{c^{-1}(z-2)}{P_2})\right),\\
\vdots\notag\\
V(1)-V(0)=-\frac{1}{2}\left(\varrho(\frac{c^{-1}}{P_2})+0\right).\label{15-40-21}
\end{gather}
Adding up \eqref{15-40-20} and \eqref{15-40-21}, we obtain
\begin{gather*}
V(z)-V(0)=
-\sum_{j=1}^{z}\varrho(\frac{c^{-1}j}{P_2})+O(1).
\end{gather*}
Thus the relation \eqref{15-40-18} is proved.
\end{description}
Let $h=c^{-1}\equiv y_4^{-1}y_3\pmod{P_2},\,1\le h<P_2.$ Substituting \eqref{15-40-18} into \eqref{15-40-17}, we obtain
\begin{gather}\label{15-40-22}
\sum_{0<n\le P}\left(\varrho_{1}^{+}\varrho_{2}^{+}+\varrho_{1}^{-}\varrho_{2}^{-}
-\varrho_{1}^{-}\varrho_{2}^{+}-\varrho_{1}^{+}\varrho_{2}^{-}\right)=
2\sum_{cx-y\le j\le cx+y}\varrho(\frac{hj}{P_2})+O(1).
\end{gather}
We transform the sum in the right side of \eqref{15-40-22}. The change of variables $j=cx+n$ leads to
\begin{gather*}
\sum_{cx-y\le j\le cx+y}\varrho(\frac{hj}{P_2})=\sum_{-y\le n\le y}\varrho(\frac{hn+x}{P_2})=\\=
\sum_{0<n\le y}\left(\varrho(\frac{hn+x}{P_2})-\varrho(\frac{hn-x}{P_2})\right)+O(1)\le\\\le
\sum_{0<n\le y}\left(\frac{2x}{P_2}+\varrho(\frac{hn+x}{P_2})-\varrho(\frac{hn-x}{P_2})\right)+O(1)=
\sum_{0<n\le y}\1_{\{\|\frac{hn}{P_2}\|\le\frac{x}{P_2}\}}+O(1).
\end{gather*}
Hence,
\begin{gather}\label{15-40-23}
\sum_{0<n\le P}\left(\varrho_{1}^{+}\varrho_{2}^{+}+\varrho_{1}^{-}\varrho_{2}^{-}
-\varrho_{1}^{-}\varrho_{2}^{+}-\varrho_{1}^{+}\varrho_{2}^{-}\right)\le
2\sum_{0<n\le y}\1_{\{\|\frac{hn}{P_2}\|\le\frac{x}{P_2}\}}+O(1).
\end{gather}
By \eqref{15-40-12}, taking into account the inequality $x\ge y\ge1,$ we obtain
\begin{gather}\label{15-40-24}
\frac{x}{P_2}\le\frac{\delta_2}{R}<\frac{x+1}{P_2}\le\frac{2x}{P_2},\quad
\frac{y}{P_2}\le\frac{\delta_1}{R}<\frac{y+1}{P_2}\le\frac{2y}{P_2},
\end{gather}
and, therefore,
\begin{gather}\label{15-40-25}
\frac{\delta_2}{2R}<\frac{x}{P_2}\le\frac{\delta_2}{R},\quad
\frac{\delta_1}{2R}<\frac{y}{P_2}\le\frac{\delta_1}{R}.
\end{gather}
In particular, $\frac{\delta_2}{R}\ge\frac{1}{P_2}.$ Then one has $P_2\ge\frac{R}{\delta_2}>\frac{y_2}{\delta_2}=y_4$ and, in the same manner, $P_2>y_3.$ Further, by the definition one has $hy_4\equiv y_3 \pmod{P_2}.$ Then in view of $1\le h<P_2,\, 0<y_3,y_4<P_2$ there exists $k$ such that $$hy_4=y_3+kP_2$$ and $0\le k<y_4.$ We note that $(y_3,y_4)=1$ implies that $(k,y_4)=1.$ Hence,
\begin{gather}\label{15-40-26}
\frac{hn}{P_2}=\frac{n}{P_2}\frac{y_3+kP_2}{y_4}=\frac{nk+\frac{ny_3}{P_2}}{y_4}.
\end{gather}
We consider two cases depending on the value of $k.$
\begin{enumerate}
  \item
If $k=0,$ then $hy_4=y_3.$ So, because of $(y_4,y_3)=1,$ we obtain $y_4=1,$ that is, $y_2=\delta_2.$ Then it follows from \eqref{15-40-26} that
\begin{gather}\label{15-40-27}
\sum_{0<n\le y}\1_{\{\|\frac{hn}{P_2}\|\le\frac{x}{P_2}\}}=
\sum_{0<n\le y}\1_{\{\|\frac{ny_3}{P_2}\|\le\frac{x}{P_2}\}}.
\end{gather}
Using the trivial estimate of the right side of \eqref{15-40-27} and the bound \eqref{15-40-25} we obtain
\begin{gather}\label{15-40-28}
\sum_{0<n\le y}\1_{\{\|\frac{hn}{P_2}\|\le\frac{x}{P_2}\}}\le
y\le\frac{\delta_1P_2}{R}=\frac{P}{R\delta_2}=\frac{P}{Ry_2}.
\end{gather}
On the other hand, since $yy_3=y\frac{y_1}{\delta_1}\le\frac{\delta_1P_2}{R}\frac{y_1}{\delta_1}<\frac{P_2}{2},$ one has
\begin{gather}\label{15-40-29}
\sum_{0<n\le y}\1_{\{\|\frac{ny_3}{P_2}\|\le\frac{x}{P_2}\}}\le\frac{x}{y_3}\le\frac{\delta_2P_2}{R}\frac{\delta_1}{y_1}=
\frac{P}{Ry_1}.
\end{gather}
Hence, by \eqref{15-40-29} and \eqref{15-40-28} we obtain
\begin{gather}\label{15-40-30}
\sum_{0<n\le y}\1_{\{\|\frac{hn}{P_2}\|\le\frac{x}{P_2}\}}\le
\frac{P}{R}\min\{\frac{1}{y_1},\frac{1}{y_2}\}.
\end{gather}
  \item
If $k\neq0,$ then it is necessary to study the inequality
\begin{gather}\label{15-40-31}
\|\frac{nk+\frac{ny_3}{P_2}}{y_4}\|\le\frac{x}{P_2}.
\end{gather}
Because
\begin{gather*}
\frac{x}{P_2}=\frac{1}{y_4}\frac{xy_4}{P_2}\le\frac{1}{y_4}\frac{\delta_2y_4}{R}=
\frac{1}{y_4}\frac{y_2}{R}\le\frac{1}{10y_4},
\end{gather*}
\begin{gather*}
\frac{ny_3}{P_2}\le\frac{yy_3}{P_2}\le\frac{\delta_1}{R}\frac{y_1}{\delta_1}\le\frac{y_1}{R}\le\frac{1}{10}
\end{gather*}
the inequality \eqref{15-40-31} holds only if $y_4|nk.$ Since $(k,y_4)=1,$ one has $y_4|n.$ Therefore,
\begin{gather}\label{15-40-32}
\sum_{0<n\le y}\1_{\{\|\frac{hn}{P_2}\|\le\frac{x}{P_2}\}}\le
\sum_{0<n\le \frac{y}{y_4}}\1_{\{\|\frac{ny_3}{P_2}\|\le\frac{x}{P_2}\}}\le
\min\{\frac{y}{y_4},\frac{x}{y_3}\}\le
\frac{P}{R}\min\{\frac{1}{y_1},\frac{1}{y_2}\}.
\end{gather}
\end{enumerate}
Substituting \eqref{15-40-32} or \eqref{15-40-30} into \eqref{15-40-23}, we obtain
\begin{gather}\label{15-40-33}
\sum_{0<n\le P}\left(\varrho_{1}^{+}\varrho_{2}^{+}+\varrho_{1}^{-}\varrho_{2}^{-}
-\varrho_{1}^{-}\varrho_{2}^{+}-\varrho_{1}^{+}\varrho_{2}^{-}\right)\le
\frac{2P}{R}\min\{\frac{1}{y_1},\frac{1}{y_2}\}+O(1).
\end{gather}
Substituting \eqref{15-40-6} and \eqref{15-40-33} into \eqref{15-40-2}, we obtain \eqref{ded-1}. This completes the proof of the lemma.
\end{proof}
\end{Le}
\section{The case $\mu=3$.}\label{mu=3}
Let $Q_1,Q,\beta$ be given. For any $q$ in $Q_1\le q\le Q$ we define by any means the number $a_q,$ such that $(a_q,q)=1,\,  0\le a_q\le q.$ Let denote
\begin{gather}\label{15-0}
Z^{*}=\left\{
\theta=\frac{a_q}{q}+\beta\;\Bigl|\;Q_1\le q\le Q
\right\}.
\end{gather}
\begin{Le}\label{lemma-15-1}
Let $\mu=3,\,\lambda=2,\,$ then for any $Z\subseteq Z^{*}$ the following bound on the cardinality of the set $\M(g_{2})$ holds
\begin{gather}\label{15-1}
|\M(g_2)|\le|Z|\left|\Omega^{(3)}\right|\left(1+\frac{2H_3}{Q_1}\right)
\left(1+\frac{146A^2N}{\KK Q_1M^{(1)}\|g_2\|}\right).
\end{gather}
\begin{proof}
We use the partition of the ensemble $\Omega_N$ given by the formula \eqref{14-6} (that is, by Theorem \ref{theorem13-4}). For $\mu=3$ the equation \eqref{14-44} becomes the congruence
\begin{gather}\label{15-2}
\left(g_2(g_3^{(1)}\frac{a^{(1)}}{q^{(1)}}-g_3^{(2)}\frac{a^{(2)}}{q^{(2)}})\right)_{1,2}\equiv 0\pmod{1},
\end{gather}
where the subscripts "1,2" mean that the congruence \eqref{15-2} holds for both coordinates of the vector. Let $\q=[q^{(1)},q^{(2)}],$ then \eqref{15-2} can be written as
\begin{gather}\label{15-3}
\left(g_2(g_3^{(1)}\frac{a^{(1)}q^{(2)}}{(q^{(1)},q^{(2)})}-g_3^{(2)}\frac{a^{(2)}q^{(1)}}{(q^{(1)},q^{(2)})})\right)_{1,2}\equiv 0\pmod{\q}.
\end{gather}
Since $det\,g_2=1,$ the congruence \eqref{15-3} can be simplified to
\begin{gather}\label{15-4}
\left(g_3^{(1)}\frac{a^{(1)}q^{(2)}}{(q^{(1)},q^{(2)})}-g_3^{(2)}\frac{a^{(2)}q^{(1)}}{(q^{(1)},q^{(2)})}\right)_{1,2}\equiv 0\pmod{\q}.
\end{gather}
By setting $g_3^{(1)}=(x_1,x_2)^{t},\,g_3^{(2)}=(y_1,y_2)^{t}$ in \eqref{15-4} we obtain the congruence
\begin{gather}\label{15-5}
x_1\frac{a^{(1)}q^{(2)}}{(q^{(1)},q^{(2)})}\equiv y_1\frac{a^{(2)}q^{(1)}}{(q^{(1)},q^{(2)})}\pmod{\q}
\end{gather}
and, the same one for $x_2,y_2.$ But $(a^{(1)},\frac{q^{(1)}}{(q^{(1)},q^{(2)})})\le(a^{(1)},q^{(1)})=1$ and, therefore,
\begin{gather*}
x_1\equiv 0 \pmod{\frac{q^{(1)}}{(q^{(1)},q^{(2)})}},\quad
x_2\equiv 0 \pmod{\frac{q^{(1)}}{(q^{(1)},q^{(2)})}}
\end{gather*}
and the same for $y_1,y_2.$ At the same time $(x_1,x_2)=(y_1,y_2)=1$ as the component of the vectors $g_3^{(1)},\, g_3^{(2)},$ thus
\begin{gather}\label{15-6}
q^{(1)}=(q^{(1)},q^{(2)})=q^{(2)}=\q.
\end{gather}
Let fix $\q$, for which there are $|Z|$ choices, this gives the first factor in \eqref{15-1}. Then it follows from the conditions on the set $Z$ that $a^{(1)}=a^{(2)}$ and the congruence \eqref{15-4} can be simplified to
\begin{gather}\label{15-7}
(g_3^{(1)}-g_3^{(2)})_{1,2}\equiv 0\pmod{\q}.
\end{gather}
We choose and fix the vector $g_3^{(2)},$ for which there are $\left|\Omega^{(3)}\right|$ choices. This gives the second factor in \eqref{15-1}. In view of Theorem \ref{theorem13-4} we have $|g_3^{(1)}-g_3^{(2)}|_{1,2}\le H_3.$ Putting $z=\{z_1,z_2\}=g_3^{(1)}-g_3^{(2)}$ and using \eqref{15-7} we obtain that there are at most $\left(1+\frac{2H_3}{\q}\right)$ choices for $z_1.$ This gives the third factor in \eqref{15-1}. Finally, putting $\eta=\frac{73A^2N}{\KK M^{(1)}}$ and
$g_2=\begin{pmatrix}
a & b \\
c & d
\end{pmatrix},$
where $d=\|g_2\|,$ we write the inequality \eqref{14-43} in the form $|g_2(g_3^{(1)}-g_3^{(2)})|_{1,2}\le \eta,$ whence
\begin{gather*}
\frac{-\eta-cz_1}{d}\le z_2\le\frac{\eta-cz_1}{d}.
\end{gather*}
Therefore, taking into account the congruence \eqref{15-7}, we obtain that $z_2$ takes at most
\begin{gather*}
\frac{2\eta}{\q\|g_2\|}+1=\left(1+\frac{146A^2N}{\q\KK M^{(1)}\|g_2\|}\right)
\end{gather*}
values. This gives the third factor in \eqref{15-1}. This completes the proof of the lemma.
\end{proof}
\end{Le}
Let
\begin{gather}\label{15-8-0}
Q_0=\max\left\{\exp\left(\frac{10^5A^4}{\epsilon_0^2}\right),\exp(\epsilon_0^{-5})\right\}.
\end{gather}
\begin{Le}\label{lemma-15-2}
If
$
\KK^2Q^{3}\le\frac{N^{1-\epsilon_0}}{12000A^4},\, \KK Q\ge Q_0,
$
then the following bound holds
\begin{gather}\label{15-9}
\sum_{\theta\in P_{Q_1,Q}^{(\beta)}}\left|S_N(\theta)\right|^2\ll
|\Omega_N|^2\KK^{12\epsilon_0}Q^{20\epsilon_0}\frac{\KK^{4(1-\delta)}Q^{6(1-\delta)+1}}{\KK Q^2_1}
\end{gather}
\begin{proof}
We use the partition of the ensemble $\Omega_N$ given by the formula \eqref{14-6} (that is, by Theorem \ref{theorem13-4}). Let
$Z\subseteq Z^{*}$ be any subset. We put $\mu=3,\lambda=2$ and
\begin{gather}\label{15-10}
M^{(1)}=76A^2\KK Q^{2},\quad M^{(3)}=76A^2\KK Q.
\end{gather}
Then the condition  \eqref{14-45} holds. It follows from the statement of the lemma that inequalities \eqref{13-42} and \eqref{13-43} hold. So all conditions of Lemma \ref{lemma-14-4} hold and therefore the estimate \eqref{14-46} is valid. Since
$H_3\ge M^{(3)}\ge Q_1$ the third factor in \eqref{15-1} can be replaced by $\frac{3H_3}{Q_1}.$ Let consider the fourth factor in \eqref{15-1}.
Using the lower bound of \eqref{13-49-1}, we obtain:
\begin{gather}\label{15-11}
1+\frac{146A^2N}{\KK Q_1M^{(1)}\|g_2\|}\le1+
\frac{146A^2N}{\KK Q_1M^{(1)}}\frac{3\cdot10^4A^5(M^{(1)}M^{(3)})^{1+2\epsilon_0}}{N}\ll
\left(\KK^2Q^3\right)^{2\epsilon_0}\frac{Q}{Q_1},
\end{gather}
where we have used \eqref{15-10}. Substituting the result of the simplification into \eqref{15-1}, we obtain
\begin{gather}\label{15-12}
|\M(g_2)|\ll
|Z|\left|\Omega^{(3)}\right|\left(\KK^2Q^3\right)^{2\epsilon_0}\frac{H_3Q}{Q_1^2}.
\end{gather}
Substituting $H_3(M^{(3)})$ from \eqref{14-10}, we obtain
\begin{gather}\label{15-13}
|\M(g_2)|\ll
|Z|\left|\Omega^{(3)}\right|\left(\KK^3Q^4\right)^{2\epsilon_0}\frac{\KK Q^2}{Q_1^2}.
\end{gather}
Substituting the obtained bound on $|\M(g_2)|$ into \eqref{14-46}, we have
\begin{gather}\label{15-14}
\sigma_{N,Z}\ll
H_1\left|\Omega^{(1)}\right|^{1/2}\left|\Omega^{(2)}\right|\left|\Omega^{(3)}\right|^{1/2}|Z|^{1/2}
\left(\KK^3Q^4\right)^{\epsilon_0}\frac{\KK^{1/2}Q}{Q_1}.
\end{gather}
Using the estimate \eqref{13-55}, we obtain
\begin{gather*}
\sigma_{N,Z}\ll
|\Omega_N|\frac{(M^{(1)})^{1+2\epsilon_0}}{(M^{(1)}M^{(3)})^{\delta-\epsilon_0/2}}|Z|^{1/2}
\left(\KK^3Q^4\right)^{\epsilon_0}\frac{\KK^{1/2}Q}{Q_1}.
\end{gather*}
Substituting $M^{(1)},M^{(3)}$ into \eqref{15-10}, we have
\begin{gather}\label{15-15}
\sigma_{N,Z}\ll
|Z|^{1/2}|\Omega_N|
\KK^{6\epsilon_0}Q^{9.5\epsilon_0}\frac{\KK^{3/2-2\delta}Q^{3(1-\delta)}}{Q_1}.
\end{gather}
Applying Lemma \ref{lemma-14-5} with $W=Z^{*},C_2=0,$ we deduce from \eqref{15-15} that
\begin{gather}\label{15-16}
\sum_{\theta\in Z^*}\left|S_N(\theta)\right|^2\ll
|\Omega_N|^2\KK^{12\epsilon_0}Q^{20\epsilon_0}\frac{\KK^{4(1-\delta)}Q^{6(1-\delta)}}{\KK Q^2_1}.
\end{gather}
Using the trivial bound
\begin{gather}\label{15-16-1}
\sum_{\theta\in P_{Q_1,Q}^{(\beta)}}\left|S_N(\theta)\right|^2\le
Q\sum_{Q_1\le q\le Q}\max\limits_{1\le a\le q, (a,q)=1}\left|S_N(\frac{a}{q}+\frac{K}{N})\right|^2=Q\sum_{\theta\in Z^*}\left|S_N(\theta)\right|^2,
\end{gather}
where as $a_q$ we have chosen numerators for which the maximum is achieved, and substituting the obtained bound \eqref{15-16} into \eqref{15-16-1}, we obtained the desired estimate. This completes the proof of the lemma.
\end{proof}
\end{Le}
\begin{Le}\label{lemma-15-2-1}
If
$
\KK^2q^{2}\le\frac{N^{1-\epsilon_0}}{12000A^4},\, \KK q\ge Q_0,
$
then for $\theta=\frac{a}{q}+\frac{K}{N}$ the following bound holds
\begin{gather}\label{15-18}
\left|S_N(\theta)\right|\ll
|\Omega_N|(\KK q)^{6\epsilon_0}\frac{(\KK q)^{2(1-\delta)}}{\KK^{1/2}q}.
\end{gather}
\begin{proof}
Let $Z=\{\theta\}.$ We use the partition of the ensemble $\Omega_N$ given by the formula \eqref{14-6} (that is, by Theorem \ref{theorem13-4}). We put $\mu=3,\lambda=2$ and
\begin{gather}\label{15-19}
M^{(1)}=76A^2\KK q,\quad M^{(3)}=76A^2\KK q.
\end{gather}
Then the condition  \eqref{14-45} holds. It follows from the statement of the lemma that inequalities \eqref{13-42} and \eqref{13-43} hold. So all conditions of Lemma \ref{lemma-14-4} hold. Using Lemma \ref{lemma-15-1} with $|Z|=1$ and making all transformations in the same way as in  Lemma \ref{lemma-15-2}, we obtain
\begin{gather}\label{15-19}
\left|S_N(\theta)\right|=\sigma_{N,Z}\ll
|\Omega_N|(\KK q)^{6\epsilon_0}\frac{(\KK q)^{2(1-\delta)}}{\KK^{1/2}q}.
\end{gather}
This completes the proof of the lemma.
\end{proof}
\end{Le}


\begin{Le}\label{lemma-15-3}
If $\KK q\ge Q_0,$ then the following bound holds
\begin{gather}\label{15-22}
|S_N(\theta)|\ll|\Omega_N|\frac{(\KK q)^{2\epsilon_0}N^{1-\delta+\epsilon_0}}{\KK q}.
\end{gather}
\begin{proof}
We use the partition of the ensemble $\Omega_N$ given by the formula \eqref{14-8} (that is, by Theorem \ref{theorem13-3}).
We put $\mu=3,\lambda=2,$ $Z=\{\theta\}$ and
\begin{gather}\label{15-23}
M^{(1)}=76A^2\KK q.
\end{gather}
Then the condition  \eqref{14-45} holds and therefore the estimate \eqref{14-46} is valid. Since $\Omega^{(2)}=\{E\},$  the conditions \eqref{14-43} and \eqref{14-44} can be written as
\begin{gather}\label{15-24}
g_3^{(1)}\equiv g_3^{(2)} \pmod{q},\quad
|g_3^{(1)}-g_3^{(2)}|_{1,2}\le \frac{73A^2N}{M^{(1)}\KK}.
\end{gather}
Thus we obtain
\begin{gather}\notag
\left|\M(g_{\lambda})\right|\le
\sum_{g_3^{(2)}\in\Omega^{(3)}}\sum_{g_3^{(1)}\in\Z^2}\1_{\{
g_3^{(1)}\equiv g_3^{(2)} \pmod{q},\,|g_3^{(1)}-g_3^{(2)}|_{1,2}\le \frac{N}{\KK^2q}
\}}\le\\\le|\Omega^{(3)}|\left(1+\frac{N}{\KK^2q^2}\right)^2\ll\frac{|\Omega^{(3)}|N^2}{\KK^4q^4},\label{15-25}
\end{gather}
as, in view of \eqref{14-26} we have $\KK q\le N^{1/2}.$ Substituting \eqref{15-25} into \eqref{14-46}, we obtain
\begin{gather}\label{15-26}
|S_N(\theta)|\ll H_1\left|\Omega^{(1)}\right|^{1/2}|\Omega^{(3)}|^{1/2}\frac{N}{\KK^2q^2}\ll
|\Omega_N|\frac{(M^{(1)})^{1+2\epsilon_0}}{|\Omega_N|^{1/2}}\frac{N}{\KK^2q^2}.
\end{gather}
Using \eqref{15-23} and the lower bound of \eqref{13-39}, we have
\begin{gather*}
|S_N(\theta)|\ll|\Omega_N|\frac{(\KK q)^{1+2\epsilon_0}N^{1-\delta+\epsilon_0}}{\KK^2q^2}.
\end{gather*}
This completes the proof of the lemma.
\end{proof}
\end{Le}
\begin{Le}\label{lemma-15-4}
Let the following inequalities hold
$N^{\epsilon_0/2}\le Q^{1/2}\le Q_1\le Q,\,\KK Q\le N^{\alpha},$
and $\frac{1}{4}<\alpha\le\frac{1}{2}+\epsilon_0.$ Then the bound
\begin{gather}\label{15-27}
\sum_{\theta\in P_{Q_1,Q}^{(\beta)}}|S_N(\theta)|\ll
|\Omega_N|\left(
N^{1/2+\alpha/2-\delta+3\epsilon_0}Q+
N^{1-\delta+3\epsilon_0}\frac{Q^{1/2}}{\KK}\right)
\end{gather}
holds. Moreover, for any $Z\subseteq P_{Q_1,Q}^{(\beta)}$ the following bound holds
\begin{gather}\notag
\sum_{\theta\in Z}|S_N(\theta)|\ll
|\Omega_N||Z|^{1/2}\left(
\frac{N^{1-\delta+2\epsilon_0}}{\KK Q^{1/2}_1}+
N^{1-\frac{3-\alpha}{2}\delta+4,5\epsilon_0}+
N^{\frac{3+\alpha}{4}-\frac{3-\alpha}{2}\delta+4,5\epsilon_0}Q^{1/2}\right)+\\+
|\Omega_N|N^{\frac{1+\alpha}{2}-\delta+2.5\epsilon_0}Q.\label{15-27-1}
\end{gather}
\begin{proof}
The beginning of both proofs is the same. We use the partition of the ensemble $\Omega_N$ given by the formula \eqref{14-8} (that is, by Theorem \ref{theorem13-3}). We put $\mu=3,\lambda=2,Z\subseteq P_{Q_1,Q}^{(\beta)}$ and
\begin{gather}\label{15-28}
M^{(1)}=150A^2N^{1/2+\alpha/2}.
\end{gather}
Then $M^{(1)}\ge150A^2N^{\alpha}\ge150A^2\KK$ and so the conditions of Lemma \ref{lemma-14-3} hold. Therefore, the inequality \eqref{14-34} holds. It follows from \eqref{14-31} and \eqref{14-32} that
\begin{gather}\label{15-29}
\|g_3\frac{a^{(1)}}{q^{(1)}}-g_3^{(2)}\frac{a^{(2)}}{q^{(2)}}\|_{1,2}\le \frac{74A^2\KK}{M^{(1)}},
\end{gather}
\begin{gather}\label{15-30}
|g_3-g_3^{(2)}|_{1,2}\le\min\left\{
\frac{73A^2N}{M^{(1)}},\; \frac{73A^2N}{M^{(1)}\overline{K}}+\frac{N}{\overline{K}}
\left\|g_3\frac{a^{(1)}}{q^{(1)}}-g_3^{(2)}\frac{a^{(2)}}{q^{(2)}}\right\|_{1,2}.
\right\},
\end{gather}
Let $g_3^{(1)}=(x_1,x_2)^{t},\,g_3^{(2)}=(y_1,y_2)^{t},$ then it follows from Theorem \ref{theorem13-3} that
\begin{gather}\label{15-31}
\max\{y_1,y_2\}=y_2\le\frac{73A^2N}{M^{(1)}}.
\end{gather}
We denote
\begin{gather}\label{15-32}
Y=
\begin{pmatrix}
x_1 & y_1 \\
x_2 & y_2
\end{pmatrix},\quad
\mathcal{Y}=det(Y)=x_1y_2-y_1x_2.
\end{gather}
Then it follows from the triangle inequality that
\begin{gather}\label{15-33}
\|\mathcal{Y}\frac{a^{(1)}}{q^{(1)}}\|\le
\|y_2(x_1\frac{a^{(1)}}{q^{(1)}}-y_1\frac{a^{(2)}}{q^{(2)}})\|+\|y_1(y_2\frac{a^{(2)}}{q^{(2)}}-x_2\frac{a^{(1)}}{q^{(1)}})\|.
\end{gather}
Applying \eqref{15-29},\,\eqref{15-31} and  \eqref{15-28} we similarly estimate both summands in the right side of \eqref{15-33}
\begin{gather}\label{15-34}
\|\mathcal{Y}\frac{a^{(1)}}{q^{(1)}}\|\le
2\frac{74A^2\KK}{M^{(1)}}\frac{73A^2N}{M^{(1)}}<\frac{1}{2Q}<\frac{1}{q^{(1)}}.
\end{gather}
Hence, $\mathcal{Y}\equiv 0 \pmod{q^{(1)}}$ and similarly $\mathcal{Y}\equiv 0 \pmod{q^{(2)}}.$ So one has
$\mathcal{Y}\equiv 0 \pmod{\q},$ where $\q=[q^{(1)},q^{(2)}]$ and so $Q_1\le\q\le Q^2.$\par
To prove the inequality \eqref{15-27} we put $Z=P_{Q_1,Q}^{(\beta)}$ and use the bound \eqref{14-34}. The set $\NN(g_{2})$ can be represented as a union of the sets $\M_1,\M_2.$ For the first set $\mathcal{Y}=0,$ for the second one $\mathcal{Y}\neq0.$
\begin{enumerate}
  \item
  To prove the estimate \eqref{15-27} we use the following bounds
\begin{gather}\label{15-35}
|\M_1|\ll_{\epsilon}Q^{2+\epsilon}|\Omega^{(3)}|,
\end{gather}
\begin{gather}\label{15-36}
|\M_2|\ll_{\epsilon}
|\Omega^{(3)}|\left(\left(\frac{73A^2N}{M^{(1)}\KK}\right)^2+1\right)N^{\epsilon}Q,
\end{gather}
which will be proved below in Lemma \ref{lemma-15-5} and \ref{lemma-15-6} respectively. Hence,
\begin{gather}\label{15-36-1}
|\NN(g_{2})|^{1/2}\ll_{\epsilon}|\Omega^{(3)}|^{1/2}\left(Q^{1+\epsilon}+\frac{73A^2N^{1+\epsilon}Q^{1/2}}{M^{(1)}\KK}+N^{\epsilon}Q^{1/2}\right).
\end{gather}
Substituting \eqref{15-36-1} into \eqref{14-34}, we obtain for $\epsilon=\epsilon_0:$
\begin{gather}\label{15-36-2}
\sigma_{N,Z}\ll
10H_1\left|\Omega^{(1)}\right|^{1/2}
|\Omega^{(3)}|^{1/2}\left(Q^{1+\epsilon_0}+\frac{73A^2N^{1+\epsilon}Q^{1/2}}{M^{(1)}\KK}+N^{\epsilon_0}Q^{1/2}\right).
\end{gather}
We note that it follows from the condition of the lemma that $Q^{1/2}N^{\epsilon_0/2}\le Q,$ so the last inequality in \eqref{15-36-2} can be omitted. Using the bounds \eqref{13-55} and \eqref{13-39}, we obtain
\begin{gather}\label{15-36-4}
|\Omega^{(1)}|^{1/2}|\Omega^{(3)}|^{1/2}=|\Omega_N|\frac{1}{|\Omega_N|^{1/2}}\le
|\Omega_N|\frac{1}{N^{\delta-\epsilon_0/2}}.
\end{gather}
Substituting \eqref{15-36-4} into \eqref{15-36-2}, we obtain
\begin{gather}\label{15-36-5}
\sigma_{N,Z}\ll
|\Omega_N|\frac{1}{N^{\delta-\epsilon_0/2}}
\left((M^{(1)})^{1+2\epsilon_0}Q^{1+\epsilon_0}+\frac{N^{1+\epsilon_0}Q^{1/2}(M^{(1)})^{2\epsilon_0}}{\KK}\right)
\end{gather}
Substituting $M^{(1)}$ from \eqref{15-28}, we obtain
\begin{gather}\label{15-36-6}
\sigma_{N,Z}\ll|\Omega_N|\left(
N^{1/2+\alpha/2-\delta+2.5\epsilon_0}Q+
N^{1-\delta+3\epsilon_0}\frac{Q^{1/2}}{\KK}\right).
\end{gather}
Thus the bound \eqref{15-27} is proved.
\item
To prove the estimate \eqref{15-27-1} we use the bound \eqref{15-35} and
\begin{gather}\label{15-36-8}
|\M_2|\ll_{\epsilon}|Z|\left(
\left(\frac{73A^2N}{M^{(1)}}\right)^2\frac{|\Omega^{(3)}|}{\KK^2 Q_1}+\left(\frac{73A^2N}{M^{(1)}}\right)^2Q^{\epsilon}+
\left(\frac{73A^2N}{M^{(1)}}\right)Q^{1+\epsilon}\right),
\end{gather}
the last will be proved in Lemma \ref{lemma-15-7}. Hence,
\begin{gather}\label{15-36-9}
|\NN(g_{2})|^{1/2}\ll_{\epsilon}|Z|^{1/2}\left(
\frac{N}{M^{(1)}}\frac{|\Omega^{(3)}|^{1/2}}{\KK Q^{1/2}_1}+
\frac{NQ^{\epsilon}}{M^{(1)}}+
\left(\frac{N}{M^{(1)}}\right)^{1/2}Q^{1/2+\epsilon}\right)+Q^{1+\epsilon}|\Omega^{(3)}|^{1/2}.
\end{gather}
Substituting \eqref{15-36-9} into \eqref{14-34}, we obtain for $\epsilon=\epsilon_0:$
\begin{gather}\notag
\sigma_{N,Z}\ll
(M^{(1)})^{1+2\epsilon_0}\left|\Omega^{(1)}\right|^{1/2}|Z|^{1/2}\left(
\frac{N}{M^{(1)}}\frac{|\Omega^{(3)}|^{1/2}}{\KK Q^{1/2}_1}+
\frac{NQ^{\epsilon_0}}{M^{(1)}}+
\left(\frac{N}{M^{(1)}}\right)^{1/2}Q^{1/2+\epsilon_0}\right)+\\+
(M^{(1)})^{1+2\epsilon_0}\left|\Omega^{(1)}\right|^{1/2}Q^{1+\epsilon_0}|\Omega^{(3)}|^{1/2}.\label{15-36-10}
\end{gather}
Using the bounds \eqref{13-55} and \eqref{13-39}, we have
\begin{gather}\label{15-36-11}
\left|\Omega^{(1)}\right|^{1/2}=|\Omega_N|\frac{|\Omega^{(1)}|^{1/2}}{|\Omega_N|}\le
|\Omega_N|\frac{(M^{(1)})^{\delta+2\epsilon_0}}{N^{2\delta-\epsilon_0}}.
\end{gather}
Substituting \eqref{15-36-11} and \eqref{15-36-4} into \eqref{15-36-10}, we obtain
\begin{gather*}\label{15-36-12}
\sigma_{N,Z}\ll
|\Omega_N||Z|^{1/2}\Bigl(
\frac{N^{1-\delta+\epsilon_0/2}(M^{(1)})^{2\epsilon_0}}{\KK Q^{1/2}_1}+
N^{1-2\delta+\epsilon_0}(M^{(1)})^{\delta+4\epsilon_0}Q^{\epsilon_0}+\\+
N^{1/2-2\delta+\epsilon_0}(M^{(1)})^{1/2+\delta+4\epsilon_0}Q^{1/2+\epsilon_0}\Bigr)+
|\Omega_N|\frac{(M^{(1)})^{1+2\epsilon_0}}{N^{\delta-\epsilon_0/2}}
Q^{1+\epsilon}.
\end{gather*}
Substituting $M^{(1)}$ from \eqref{15-28}, we obtain
\begin{gather*}\label{15-36-12}
\sigma_{N,Z}\ll
|\Omega_N||Z|^{1/2}\left(
\frac{N^{1-\delta+2\epsilon_0}}{\KK Q^{1/2}_1}+
N^{1-\frac{3-\alpha}{2}\delta+4,5\epsilon_0}+
N^{\frac{3+\alpha}{4}-\frac{3-\alpha}{2}\delta+4,5\epsilon_0}Q^{1/2}\right)+\\+
|\Omega_N|N^{\frac{1+\alpha}{2}-\delta+2,5\epsilon_0}Q.
\end{gather*}
Thus the bound \eqref{15-27-1} is proved.
\end{enumerate}
This completes the proof of the lemma.
\end{proof}
\end{Le}
\begin{Le}\label{lemma-15-5}
Under the hypotheses of Lemma \ref{lemma-15-4} one has
\begin{gather}\label{15-37}
|\M_1|\ll_{\epsilon}Q^{2+\epsilon}|\Omega^{(3)}|.
\end{gather}
\begin{proof}
To simplify we denote $$R=(M^{(1)})^{1+2\epsilon_0},T=\frac{73A^2N}{M^{(1)}}.$$
We recall that $M^{(1)}$ is defined in \eqref{15-28}. It follows from the condition $\mathcal{Y}=0$ that $\frac{x_1}{x_2}=\frac{y_1}{y_2}.$ Since
\begin{gather*}
g_3^{(1)}=(x_1,x_2)^{t},\,g_3^{(2)}=(y_1,y_2)^{t},
\end{gather*}
one has $(x_1,x_2)=1,(y_1,y_2)=1$ and, hence, $x_1=y_1,x_2=y_2.$ In particular, ${|\beta(g_3^{(1)}-g_3^{(2)})|_{1,2}=0.}$ Then it follows from the inequality  \eqref{14-29} that
\begin{gather}\label{15-38}
\|y_1\left(\frac{a^{(1)}}{q^{(1)}}-\frac{a^{(2)}}{q^{(2)}}\right)\|\le
\|y_1\left(\theta^{(1)}-\theta^{(2)}\right)\|+0\le\frac{1}{H_1}<\frac{1}{R}
\end{gather}
and, in the same way, $\|y_2\left(\frac{a^{(1)}}{q^{(1)}}-\frac{a^{(2)}}{q^{(2)}}\right)\|\le\frac{1}{H_1}<\frac{1}{R}.$ Thus we obtain
\begin{gather}\label{15-39}
|\M_1|\le \sum_{g_3^{(2)}\in\Om^{(3)}}\sum_{Q_1\le q^{(1)},q^{(2)}\le Q}\mathop{{\sum}^*}_{0\le a^{(1)}\le q^{(1)}}
\mathop{{\sum}^*}_{0\le a^{(2)}\le q^{(2)}}
\1_{\{\|y_{1,2}\left(\frac{a^{(1)}}{q^{(1)}}-\frac{a^{(2)}}{q^{(2)}}\right)\|<\frac{1}{R}\}}.
\end{gather}
We denote $P=[q^{(1)},q^{(2)}]$ and write $\frac{a^{(1)}}{q^{(1)}}-\frac{a^{(2)}}{q^{(2)}}=\frac{n}{P}.$ WE note that
\begin{gather*}
\#\{a^{(1)},a^{(2)}\,|\,
a^{(1)}\frac{P}{q^{(1)}}-a^{(2)}\frac{P}{q^{(2)}}=n,\,1\le a^{(1)}\le q^{(1)},\,1\le a^{(2)}\le q^{(2)}
\}\le(q^{(1)},q^{(2)}),
\end{gather*}
then \eqref{15-39} can be written as
\begin{gather}\notag
|\M_1|\le \sum_{g_3^{(2)}\in\Om^{(3)}}\sum_{Q_1\le q^{(1)},q^{(2)}\le Q}(q^{(1)},q^{(2)})\sum_{|n|<P}
\1_{\{\|y_{1}\frac{n}{P}\|<\frac{1}{R}\}}\1_{\{\|y_{2}\frac{n}{P}\|<\frac{1}{R}\}}\le\\\le
2\sum_{g_3^{(2)}\in\Om^{(3)}}\sum_{Q_1\le q^{(1)},q^{(2)}\le Q}(q^{(1)},q^{(2)})\sum_{0<n\le P}
\1_{\{\|y_{1}\frac{n}{P}\|<\frac{1}{R}\}}\1_{\{\|y_{2}\frac{n}{P}\|<\frac{1}{R}\}}.\label{15-40}
\end{gather}
Applying Lemma \ref{lemma-ded}, we obtain
\begin{gather}\label{15-40}
|\M_1|\le
8\sum_{g_3^{(2)}\in\Om^{(3)}}\sum_{Q_1\le q^{(1)},q^{(2)}\le Q}(q^{(1)},q^{(2)})
\left(\frac{(y_1,P)+(y_2,P)}{R}+\frac{P}{R}\min\{\frac{1}{y_1},\frac{1}{y_2}\}+\frac{P}{R^2}+O(1)\right).
\end{gather}
Using \eqref{13-47-1}, we have
\begin{gather*}
\frac{P}{R}\min\{\frac{1}{y_1},\frac{1}{y_2}\}=\frac{P}{Ry_2}\le\frac{P}{R}\frac{140A^2H_1}{N}\le
\frac{150A^2P}{N}=O(1).
\end{gather*}
Since $Q\le N^{1/2}\le R,$ one has $\frac{P}{R^2}\le\frac{Q^2}{R^2}=O(1).$ Hence,
\begin{gather}\label{15-41}
|\M_1|\le
8\sum_{g_3^{(2)}\in\Om^{(3)}}\sum_{Q_1\le q^{(1)},q^{(2)}\le Q}(q^{(1)},q^{(2)})
\left(\frac{(y_1,P)+(y_2,P)}{R}+O(1)\right).
\end{gather}
Let estimate the first summand in \eqref{15-41}:
\begin{gather}
\frac{1}{R}\sum_{Q_1\le q^{(1)},q^{(2)}\le Q}(q^{(1)},q^{(2)})
\sum_{y_1,y_2<T}(y_1,P)\le\frac{T}{R}\sum_{q^{(1)},q^{(2)}\le Q}(q^{(1)},q^{(2)})\sum_{y_1<T}(y_1,P).
\end{gather}
Because of
\begin{gather*}
\sum_{y_1<T}(y_1,P)\le\sum_{d|P\atop d<T}d\sum_{y_1<T\atop y_1\equiv0\pmod{d}}1\ll_{\epsilon}TP^{\epsilon}\ll_{\epsilon} TQ^{2\epsilon},
\end{gather*}
we obtain
\begin{gather}\label{15-41-1}
\frac{1}{R}\sum_{Q_1\le q^{(1)},q^{(2)}\le Q}(q^{(1)},q^{(2)})
\sum_{y_1,y_2<T}(y_1,P)\ll_{\epsilon}\frac{T^2Q^{2\epsilon}}{R}\sum_{q^{(1)},q^{(2)}\le Q}(q^{(1)},q^{(2)}).
\end{gather}
Since
\begin{gather}\label{15-41-2}
\sum_{q^{(1)},q^{(2)}\le Q}(q^{(1)},q^{(2)})\le\sum_{q^{(1)}\le Q}\sum_{d|q^{(1)}}d\sum_{q^{(2)}\le Q\atop q^{(2)}\equiv0\pmod{d}}1\le
2\sum_{q^{(1)}\le Q}\sum_{d|q^{(1)}}Q\ll_{\epsilon}Q^{2+\epsilon},
\end{gather}
one has
\begin{gather}\label{15-41-3}
\frac{1}{R}\sum_{Q_1\le q^{(1)},q^{(2)}\le Q}(q^{(1)},q^{(2)})
\sum_{y_1,y_2<T}(y_1,P)\ll_{\epsilon}\frac{T^2Q^{2+3\epsilon}}{R}.
\end{gather}
Substituting \eqref{15-41-3} into \eqref{15-41} and using \eqref{15-41-2}, we obtain
\begin{gather}\label{15-41-4}
|\M_1|
\ll_{\epsilon}\frac{T^2Q^{2+\epsilon}}{R}+Q^{2+\epsilon}|\Omega^{(3)}|\ll_{\epsilon}Q^{2+\epsilon}|\Omega^{(3)}|.
\end{gather}
This completes the proof of the lemma.
\end{proof}
\end{Le}
\begin{Le}\label{lemma-15-6}
Under the hypotheses of Lemma \ref{lemma-15-4} one has
\begin{gather}\label{15-46}
|\M_2|\ll_{\epsilon}
|\Omega^{(3)}|\left(\left(\frac{73A^2N}{M^{(1)}\KK}\right)^2+1\right)N^{\epsilon}Q.
\end{gather}
\begin{proof}
Since $|\mathcal{Y}|\le\max\{x_1y_2,x_2y_1\}\le x_2y_2,$ then applying the bound \eqref{15-31} to each factor we obtain
\begin{gather}\label{15-48}
|\mathcal{Y}|\le\left(\frac{73A^2N}{M^{(1)}}\right)^2.
\end{gather}
As $\q|\mathcal{Y}$ and $\mathcal{Y}\neq0,$ one has
\begin{gather}\label{15-49}
\q\le\min\{Q^2,|\mathcal{Y}|\}\le\min\{Q^2,\left(\frac{73A^2N}{M^{(1)}}\right)^2\}\le
\sqrt{Q^2\left(\frac{73A^2N}{M^{(1)}}\right)^2}=Q\frac{73A^2N}{M^{(1)}}.
\end{gather}
It follows from the conditions of the lemma and from \eqref{15-28} that
\begin{gather}\label{15-49-1}
\q\frac{74A^2\KK}{M^{(1)}}\le Q\frac{73A^2N}{M^{(1)}}\frac{74A^2\KK}{M^{(1)}}<1.
\end{gather}
By \eqref{15-49-1} we obtain that the right side of \eqref{15-29} is less than $\frac{1}{\q},$ whence
\begin{gather}\label{15-50}
\|g_3^{(1)}\frac{a^{(1)}}{q^{(1)}}-g_3^{(2)}\frac{a^{(2)}}{q^{(2)}}\|_{1,2}=0.
\end{gather}
It follows from the equation \eqref{15-50} that the congruence \eqref{15-4} holds, from which we deduced in Lemma \ref{lemma-15-1} that
\begin{gather}\label{15-51}
q^{(1)}=q^{(2)}=\q,\quad (g_3^{(1)}a^{(1)}-g_3^{(2)}a^{(2)})_{1,2}\equiv 0\pmod{\q}.
\end{gather}
Then by \eqref{15-30} we obtain
\begin{gather}\label{15-51-1}
|g_3^{(1)}-g_3^{(2)}|_{1,2}\le\frac{73A^2N}{M^{(1)}\KK}.
\end{gather}
Hence, if we fix $g_3^{(1)}\in\Om^{(3)},$ then for $g_3^{(2)}$ there will be at most
$\left(\frac{73A^2N}{M^{(1)}\KK}\right)^2+1$ choices. So we have defined  $\mathcal{Y}.$ But then there are at most $\tau(\mathcal{Y})$ choices for the number $q^{(1)}=q^{(2)}=\q$ since $\q|\mathcal{Y}.$ It is known (see \cite[p.91]{Korobov}) that
$$\tau(\mathcal{Y})\le(\epsilon\log2)^{-\exp(1/\epsilon)}\mathcal{Y}^{\epsilon}$$
for any $\epsilon>0.$ With $q^{(1)}$ being fixed there are at most $q^{(1)}\le Q$ choices for $a^{(1)}.$ We note that if $g_3^{(1)},\,g_3^{(2)},\,a^{(1)}$ are fixed then $a^{(2)}$ is uniquely determined by \eqref{15-51}. Actually, we write
\begin{gather}\label{15-56}
g_3^{(1)}=(x_1,x_2)^{t},\,g_3^{(2)}=(y_1,y_2)^{t},
\end{gather}
then it follows from \eqref{15-51} that
\begin{gather}\label{15-57}
x_1a^{(1)}\equiv y_1a^{(2)}\pmod{\q},\,x_2a^{(1)}\equiv y_2a^{(2)}\pmod{\q}.
\end{gather}
Since $(a^{(1)},\q)=(a^{(2)},\q)=1,$ we have
\begin{gather}\label{15-56-1}
\delta_1=(x_1,\q)=(y_1,\q),\,\delta_2=(x_2,\q)=(y_2,\q),
\end{gather}
and $(\delta_1,\delta_2)=1,$ as $(x_1,x_2)=1.$ Then one can obtain from \eqref{15-57} that
\begin{gather}\label{15-58}
a^{(2)}\equiv A_1a^{(1)}\pmod{\frac{\q}{\delta_1}},\,a^{(2)}\equiv A_2a^{(1)}\pmod{\frac{\q}{\delta_2}}.
\end{gather}
These two congruences are equivalent to the congruence modulo $[\frac{\q}{\delta_1},\frac{\q}{\delta_2}]=\q,$ and hence $a^{(2)}$ is uniquely determined. Finally we obtain
\begin{gather}\notag
|\M_2|\le \sum_{g_3^{(1)}\in\Om^{(3)}}\sum_{g_3^{(2)}\in\Om^{(3)}}\sum_{q^{(1)}|\mathcal{Y}}\sum_{0<a^{(1)}\le q^{(1)}}
\1_{\{|g_3^{(1)}-g_3^{(2)}|_{1,2}\le\frac{73A^2N}{M^{(1)}\KK}\}}\le\\\ll_{\epsilon}
|\Omega^{(3)}|\left(\left(\frac{73A^2N}{M^{(1)}\KK}\right)^2+1\right)N^{\epsilon}Q.\label{15-52}
\end{gather}
This completes the proof of the lemma.
\end{proof}
\end{Le}
\begin{Le}\label{lemma-15-7}
Under the hypotheses of Lemma \ref{lemma-15-4} one has
\begin{gather}\label{15-55}
|\M_2|\ll_{\epsilon}|Z|\left(
\left(\frac{73A^2N}{M^{(1)}}\right)^2\frac{|\Omega^{(3)}|}{\KK^2 Q_1}+\left(\frac{73A^2N}{M^{(1)}}\right)^2Q^{\epsilon}+
\left(\frac{73A^2N}{M^{(1)}}\right)Q^{1+\epsilon}\right).
\end{gather}
\begin{proof}
Repeating the arguments in the proof of Lemma \ref{lemma-15-6}, we obtain \eqref{15-51} \eqref{15-51-1}.
To simplify we write $T=\frac{73A^2N}{M^{(1)}}$ and use the notation \eqref{15-56} and \eqref{15-56-1}. If $g_3^{(1)},\,g_3^{(2)},\,a^{(1)}$ are fixed, then it has been proved in Lemma \ref{lemma-15-6} that $a^{(2)}$ is uniquely determined by \eqref{15-51}. It follows from the congruence \eqref{15-51} that ${x_1y_2\equiv x_2y_1\pmod{\q}.}$ Thus there are $|Z|$ choices for $\frac{a^{(1)}}{q^{(1)}}$ and we obtain the inequality
\begin{gather}\label{15-59}
|\M_2|\le|Z|\sum_{g_3^{(1)}\in\Omega^{(3)}}\sum_{g_3^{(2)}\in\Omega^{(3)}\atop |g_3^{(1)}-g_3^{(2)}|_{1,2}\le\frac{T}{\KK}}\1_{\{x_1y_2\equiv x_2y_1\pmod{\q}\}}
\end{gather}
We denote
$$x_3=\frac{x_1}{\delta_1},\,x_4=\frac{x_2}{\delta_2},\,y_3=\frac{y_1}{\delta_1},\,y_4=\frac{y_2}{\delta_2},\,
p=\frac{\q}{\delta_1\delta_2},$$
then
\begin{gather}\label{15-69}
|\M_2|\le|Z|\sum_{g_3^{(1)}\in\Omega^{(3)}}\sum_{y_3\le\frac{T}{\delta_1}\atop |y_3-x_3|\le\frac{T}{\KK\delta_1}}
\sum_{y_4\le\frac{T}{\delta_2}\atop |y_4-x_4|\le\frac{T}{\KK\delta_2}}
\1_{\{x_3y_4\equiv x_4y_3\pmod{p}\}}.
\end{gather}
Because of $(x_3,p)=(x_4,p)=1$ the congruence can be written as $y_4\equiv cy_3\pmod{p},$ where $c\equiv x_3^{-1}x_4\pmod{p}.$ Then,using the result for the number of solutions of such congruences in ~\cite[p.18]{Korobov}, we obtain
\begin{gather}\label{15-70}
\sum_{y_3}\sum_{y_4}\1_{\{x_3y_4\equiv x_4y_3\pmod{p}\}}=\sum_{y_3}\sum_{y_4}\1_{\{y_4\equiv cy_3\pmod{p}\}}=
\frac{1}{p}\frac{T^2}{\KK^2\delta_1\delta_2}+O(s(\frac{c}{p})\log^2p),
\end{gather}
where $s(\alpha)=\sum\limits_{1\le i \le s}a_i$ is the sum of partial quotients of the number $\alpha=[0;a_1,\ldots,a_s].$ Substituting \eqref{15-70} into \eqref{15-69}, we have
\begin{gather}\label{15-71}
|\M_2|\le|Z|\left(\frac{|\Omega^{(3)}|T^2}{\KK^2Q_1}+\log^2Q\sum_{x_1\le T}\sum_{x_2\le T}s(\frac{x_3^{-1}x_2}{q/\delta_1})\right).
\end{gather}
Using the following result of Knuth and Yao ~\cite{KnuthYao},
$$\sum_{a\le b}s(a/b)\ll b\log^2 b,$$
we obtain
\begin{gather}\label{15-72}
\sum_{x_1\le T}\sum_{x_2\le T}s(\frac{x_3^{-1}x_2}{q/\delta_1})\le\sum_{x_1\le T}
\left(\frac{T}{q/\delta_1}+1\right)\frac{q}{\delta_1}\log^2 q\le(T^2+Tq)\log^2 q.
\end{gather}
Substituting \eqref{15-72} into \eqref{15-71}, we obtain
\begin{gather}\label{15-73}
|\M_2|\ll_{\epsilon}|Z|\left(\frac{|\Omega^{(3)}|T^2}{\KK^2Q_1}+T^2Q^{\epsilon}+TQ^{1+\epsilon}\right).
\end{gather}
This completes the proof of the lemma.
\end{proof}
\end{Le}
\begin{Sl}\label{sled15-1}
Under the hypotheses of Lemma \ref{lemma-15-4}, one has
\begin{gather}\label{15-74-0}
\sum_{\theta\in P_{Q_1,Q}^{(\beta)}}|S_N(\theta)|^2\ll
|\Omega_N|^2\left(C_1^2Q^{\epsilon_0}+C_2\right),
\end{gather}
where
\begin{gather}\label{15-75}
C_1=\frac{N^{1-\delta+2\epsilon_0}}{\KK Q^{1/2}_1}+
N^{1-\frac{3-\alpha}{2}\delta+4,5\epsilon_0}+
N^{\frac{3+\alpha}{4}-\frac{3-\alpha}{2}\delta+4,5\epsilon_0}Q^{1/2},\quad
C_2=\frac{N^{\frac{3+\alpha}{2}-2\delta+3,5\epsilon_0}Q}{(\KK Q_1)^{1-2\epsilon_0}}.
\end{gather}
\begin{proof}
It was proved in Lemma \ref{lemma-15-4} that for any $Z\subseteq P_{Q_1,Q}^{(\beta)}$ one has
\begin{gather*}
\sum_{\theta\in Z}|S_N(\theta)|\ll
|\Omega_N||Z|^{1/2}C_1+|\Omega_N|N^{\frac{1+\alpha}{2}-\delta+2.5\epsilon_0}Q.
\end{gather*}
Applying Lemma \ref{lemma-14-5} with $W=P_{Q_1,Q}^{(\beta)},f(\theta)=\frac{|S_N(\theta)|}{|\Omega_N|},$ we obtain
\begin{gather}\label{15-76}
\sum_{\theta\in P_{Q_1,Q}^{(\beta)}}|S_N(\theta)|^2\ll
|\Omega_N|^2C_1^2Q^{\epsilon_0}+|\Omega_N|N^{\frac{1+\alpha}{2}-\delta+2.5\epsilon_0}Q\max_{\theta\in P_{Q_1,Q}^{(\beta)}}|S_N(\theta)|.
\end{gather}
To estimate the maximum we apply Lemma \ref{lemma-15-3} and obtain \eqref{15-74-0}. This completes the proof of the corollary.
\end{proof}
\end{Sl}
\begin{Sl}\label{sled15-2}
Under the hypotheses of Lemma \ref{lemma-15-4} with $\alpha=\frac{1+\epsilon_0}{2},$  then
\begin{gather}\label{15-74}
\mathop{{\sum}^*}_{Q_1\le q\le Q\atop 1\le a\le q }\left|S_N(\frac{a}{q}+\frac{K}{N})\right|^2\ll
|\Omega_N|^2\left(
\frac{N^{2(1-\delta)-1/4+5\epsilon_0}Q}{(\KK Q_1)^{1-2\epsilon_0}}+
\frac{N^{2(1-\delta)+4\epsilon_0}Q^{1/2}}{\KK^{2-2\epsilon_0}Q_1^{1-2\epsilon_0}}
\right),
\end{gather}
\begin{proof}
Estimating $|S_N(\frac{a}{q}+\frac{K}{N})|$ by the maximum over $a,q,$ we obtain
\begin{gather}\label{15-75}
\mathop{{\sum}^*}_{Q_1\le q\le Q\atop 1\le a\le q }\left|S_N(\frac{a}{q}+\frac{K}{N})\right|^2\le
\max\limits_{Q_1\le q\le Q\atop 1\le a\le q, (a,q)=1}\left|S_N(\frac{a}{q}+\frac{K}{N})\right|
\mathop{{\sum}^*}_{Q_1\le q\le Q\atop 1\le a\le q }\left|S_N(\frac{a}{q}+\frac{K}{N})\right|.
\end{gather}
To estimate the maximum we apply Lemma \ref{lemma-15-3}. Using the first item of Lemma \ref{lemma-15-4} with $\alpha=\frac{1+\epsilon_0}{2},$ we obtain
\begin{gather*}
\mathop{{\sum}^*}_{Q_1\le q\le Q\atop 1\le a\le q }\left|S_N(\frac{a}{q}+\frac{K}{N})\right|^2\ll
|\Omega_N|^2\left(
\frac{N^{2(1-\delta)-1/4+5\epsilon_0}Q}{(\KK Q_1)^{1-2\epsilon_0}}+
\frac{N^{2(1-\delta)+4\epsilon_0}Q^{1/2}}{\KK^{2-2\epsilon_0}Q_1^{1-2\epsilon_0}}\right).
\end{gather*}
This completes the proof of the corollary.
\end{proof}
\end{Sl}
\section{The case $\mu=2$.}
We recall that $Q_0$ is defined by the equation \eqref{15-8-0}.
\begin{Le}\label{lemma-16-1}
Under the hypotheses of Lemma \ref{BK lemma 7.1}, if the inequalities
$$
(\KK Q)^{18/5+21\epsilon_0}Q<N,\quad \KK Q\ge Q_0,$$
hold, then the following bound is valid
\begin{gather}\label{16-1}
\sum_{\theta\in P_{Q_1,Q}^{(\beta)}}\left|S_N(\theta)\right|^2\ll|\Omega_N|^2
\frac{\KK^{\frac{36}{5}(1-\delta)} Q^{\frac{46}{5}(1-\delta)+1}\KK^{34\epsilon_0}Q^{46\epsilon_0}}{(\KK Q_1)^2}.
\end{gather}
\begin{proof}
It follows from \eqref{15-16-1} that it is sufficient to prove the bound
\begin{gather}\label{16-1-2}
\sum_{\theta\in Z^{*}}\left|S_N(\theta)\right|^2\ll|\Omega_N|^2
\frac{\KK^{\frac{36}{5}(1-\delta)} Q^{\frac{46}{5}(1-\delta)}\KK^{34\epsilon_0}Q^{46\epsilon_0}}{(\KK Q_1)^2}
\end{gather}
for $Z^{*}$ in \eqref{15-0}. Let $Z\subseteq Z^{*}$ be any subset. We use the partition of the ensemble $\Omega_N$ given by the formula \eqref{14-7} (that is, by Theorem \ref{theorem13-3-1}). We put $\mu=2,\lambda=3$ and
\begin{gather}\label{16-2}
M^{(1)}=75A^2\KK Q^2,\quad M^{(2)}=(75A^2\KK Q)^{13/5+11\epsilon_0},
\end{gather}
then the condition \eqref{14-45} and all conditions of Theorem \ref{theorem13-3-1} hold. Hence, all conditions of Lemma \ref{lemma-14-4} hold and therefore the estimate \eqref{14-46} is valid. For the equation \eqref{14-44} can be written as the congruence
\begin{gather}\label{16-3}
\left((g_2^{(1)}\frac{a^{(1)}}{q^{(1)}}-g_2^{(2)}\frac{a^{(2)}}{q^{(2)}})g_3\right)_{1,2}\equiv 0\pmod{1}.
\end{gather}
In the same way, as it was done in Lemma \ref{lemma-15-1}, we obtain
\begin{gather*}
q^{(1)}=q^{(2)}=\q
\end{gather*}
and, hence, $a^{(1)}=a^{(2)},$ since the definition of the set $Z^{*}$ in \eqref{15-0}. Then the relations \eqref{14-43} and \eqref{14-44} imply that
\begin{gather}\label{16-4}
\left|(g_2^{(1)}-g_2^{(2)})g_3\right|_{1,2}\le \frac{73A^2N}{M^{(1)}\KK},\quad
\left((g_2^{(1)}-g_2^{(2)})g_3\right)_{1,2}\equiv 0\pmod{\q}.
\end{gather}
We write
\begin{gather}\label{16-5}
\eta'=g_2^{(2)}g_3,\,\eta=g_3,\,\gamma=g_2^{(1)},\,X=\|\eta'\|,\,Y=\|g_2^{(2)}\|,\,K_1=\KK\frac{XM^{(1)}}{73A^2N}.
\end{gather}
Without loss of generality we may suppose that $\|g_2^{(1)}\|\le\|g_2^{(2)}\|,$ then  $\|\gamma\|\le Y.$ It also follows from the properties of ensemble that $\|\gamma\|\asymp Y.$ Moreover, it follows from the bounds proved in Theorem \ref{theorem13-3-1}, that
\begin{gather}\label{16-5-1}
\frac{\KK}{(M^{(1)})^{4\epsilon_0}(M^{(2)})^{2\epsilon_0}}\ll K_1\ll\KK\frac{(M^{(2)})^{2\epsilon_0}}{(M^{(1)})^{2\epsilon_0}},\quad
\frac{(M^{(2)})}{(M^{(1)})^{2\epsilon_0}}\ll Y\ll\frac{(M^{(2)})^{1+2\epsilon_0}}{(M^{(1)})^{2\epsilon_0}}.
\end{gather}
The relations \eqref{16-4} can be written as
\begin{gather}\label{16-7}
|\gamma\eta-\eta'|_{1,2}<\frac{X}{K_1},\quad(\gamma\eta-\eta')_{1,2}\equiv0\pmod{q}.
\end{gather}
Let verify whether the conditions of Lemma \ref{BK lemma 7.1} hold. To do this it is sufficient to confirm that $Y<X,\,(QK_1)^{13/5}<Y,$ that is, it is enough to check that
\begin{gather}\label{16-8-1}
(\KK Q)^{13/5}<(M^{(2)})^{1-26\epsilon_0/5}(M^{(1)})^{16\epsilon_0/5}.
\end{gather}
The inequality \eqref{16-8-1} is satisfied by the choice of parameters $M^{(1)},M^{(2)}.$ Thus to estimate the cardinality of the set $\M(g_3)$ Lemma \ref{BK lemma 7.1} can be applied as follows. We fix the vector $g_2^{(2)},$ for which there are $\left|\Omega^{(2)}\right|$ choices. We also fix $\frac{a^{(1)}}{q^{(1)}},$ for which there are $|Z|$ choices. It has been proved that then $\frac{a^{(2)}}{q^{(2)}}$ is fixed. Hence, for a fixed $g_3$ we have
\begin{gather}\label{16-8}
|\M(g_3)|\ll\frac{Y^2}{(K_1Q_1)^2}\left|\Omega^{(2)}\right||Z|.
\end{gather}
Using the bounds \eqref{16-5-1}, we obtain
\begin{gather}\label{16-9}
|\M(g_3)|\ll\frac{(M^{(2)})^{2+8\epsilon_0}(M^{(1)})^{4\epsilon_0}}{(\KK Q_1)^2}\left|\Omega^{(2)}\right||Z|\ll
\frac{(\KK Q)^{26/5+43\epsilon_0}(\KK Q^2)^{4\epsilon_0}}{(\KK Q_1)^2}\left|\Omega^{(2)}\right||Z|.
\end{gather}
Substituting \eqref{16-9} into  \eqref{14-46}, we have
\begin{gather}\label{16-10}
\sum_{\theta\in Z}\left|S_N(\theta)\right|\ll|Z|^{1/2}
H_1\left|\Omega^{(1)}\right|^{1/2}\left|\Omega^{(3)}\right|\left|\Omega^{(2)}\right|^{1/2}
\frac{(\KK Q)^{13/5+22\epsilon_0}(\KK Q^2)^{2\epsilon_0}}{\KK Q_1}.
\end{gather}
Using the lower bound of \eqref{13-55} and $\left|\Omega^{(1)}\right|,\,\left|\Omega^{(2)}\right|,$ we obtain
\begin{gather}\label{16-11-0}
\sum_{\theta\in Z}\left|S_N(\theta)\right|\ll|Z|^{1/2}
|\Omega_N|\frac{(M^{(1)})^{1-\delta+2,5\epsilon_0}}{(M^{(2)})^{\delta-\epsilon_0/2}}
\frac{(\KK Q)^{13/5+22\epsilon_0}(\KK Q^2)^{2\epsilon_0}}{\KK Q_1}.
\end{gather}
Substituting \eqref{16-2} into \eqref{16-11-0}, we have
\begin{gather}\label{16-12}
\sum_{\theta\in Z}\left|S_N(\theta)\right|\ll|Z|^{1/2}|\Omega_N|
\frac{\KK^{\frac{18}{5}(1-\delta)} Q^{\frac{23}{5}(1-\delta)}\KK^{17\epsilon_0}Q^{23\epsilon_0}}{\KK Q_1}.
\end{gather}
Applying Lemma \ref{lemma-14-5} with $W=Z^*,c_2=0,f(\theta)=\frac{|S_N(\theta)|}{|\Omega_N|},$ we obtain \eqref{16-1-2}. This completes the proof of the lemma.
\end{proof}
\end{Le}
\begin{Le}\label{lemma-16-1-0}
Under the hypotheses of Lemma \ref{BK lemma 7.1}, if the inequality $\KK q>Q_0$ holds, then the following bound is valid
\begin{gather}\label{16-1-0}
|S_N(\theta)|\ll\frac{|\Omega_N|}{\KK q}(\KK q)^{\frac{18}{5}(1-\delta)+19\epsilon_0}.
\end{gather}
\begin{proof}
We use the partition of the ensemble $\Omega_N$ given by the formula \eqref{14-7} (that is, by Theorem \ref{theorem13-3-1}). We put $\mu=2,\lambda=3$ and
\begin{gather}\label{16-2-0}
M^{(1)}=75A^2\KK q,\quad M^{(2)}=(75A^2\KK q)^{13/5+11\epsilon_0},
\end{gather}
then the condition \eqref{14-45}of Lemma \ref{lemma-14-4} holds. Suppose that the inequality
\begin{gather}\label{16-3-0}
(\KK q)^{18/5+19\epsilon_0}<N
\end{gather}
is valid. Then all conditions of Theorem \ref{theorem13-3-1} hold and, hence, all conditions of Lemma \ref{lemma-14-4} hold.
So one has the bound \eqref{14-46}. Below we will use notations \eqref{16-5}. Verification of the feasibility of conditions of Lemma \ref{BK lemma 7.1}. can be done in the same manner. Thus to estimate the cardinality of the set $\M(g_3)$ Lemma \ref{BK lemma 7.1} can be applied
\begin{gather}\label{16-8-0}
|\M(g_3)|\ll\frac{Y^2}{(K_1q)^2}\left|\Omega^{(2)}\right|.
\end{gather}
Using the bounds \eqref{16-5-1}, we have
\begin{gather}\label{16-9-0}
|\M(g_3)|\ll\frac{(M^{(2)})^{2+8\epsilon_0}(M^{(1)})^{4\epsilon_0}}{(\KK q)^2}\left|\Omega^{(2)}\right|\ll
\frac{(\KK q)^{26/5+43\epsilon_0}(\KK q)^{4\epsilon_0}}{(\KK q)^2}\left|\Omega^{(2)}\right|.
\end{gather}
Substituting \eqref{16-9-0} into \eqref{14-46}, we obtain
\begin{gather}\label{16-10}
\left|S_N(\theta)\right|\ll
H_1\left|\Omega^{(1)}\right|^{1/2}\left|\Omega^{(3)}\right|\left|\Omega^{(2)}\right|^{1/2}
(\KK q)^{8/5+24\epsilon_0}.
\end{gather}
Using the lower bound of \eqref{13-55} and $\left|\Omega^{(1)}\right|,\,\left|\Omega^{(2)}\right|,$ we obtain
\begin{gather}\label{16-11}
\left|S_N(\theta)\right|\ll
|\Omega_N|(\KK q)^{\frac{18}{5}(1-\delta)-1+19\epsilon_0},
\end{gather}
and the inequality \eqref{16-1-0} is proved under the condition \eqref{16-3-0}.\par
Let next the inequality \eqref{16-3-0} be false, that is,
\begin{gather}\label{16-10-0}
N\le(\KK q)^{18/5+19\epsilon_0}.
\end{gather}
Then the conditions of Lemma \ref{lemma-15-3} hold. Applying \eqref{15-22} and taking into account \eqref{16-10-0}, we obtain
\begin{gather}\label{16-11}
|S_N(\theta)|\ll
|\Omega_N|(\KK q)^{\frac{18}{5}(1-\delta)-1+9\epsilon_0}.
\end{gather}
This completes the proof of the lemma.
\end{proof}
\end{Le}

\section{Estimates for integrals of $|S_N(\theta)|^2$.}
\begin{Le}\label{lemma-17-1}
The following inequality holds
\begin{gather}\label{17-1}
\int_0^1\left|S_N(\theta)\right|^2d\theta\le \frac{1}{N}
\mathop{{\sum}^*}_{0\le a\le q\le N^{1/2}}\int\limits_{|K|\le\frac{N^{1/2}}{q}}
\left|S_N(\frac{a}{q}+\frac{K}{N})\right|^2dK,
\end{gather}
where $\mathop{{\sum}^*}$ means that the sum is taken over $a$ and $q$ being coprime for $q\ge1,$ and $a=0,1$ for $q=1.$
\begin{proof}
It follows from the Dirichlet theorem that for any $\theta\in[0,1]$ there exist $a,q\in\N$ and $\beta\in\rr$ such that
\begin{gather*}
\theta=\frac{a}{q}+\beta,\;(a,q)=1,\; 0\le a\le q\le N^{1/2},\;|\beta|\le\frac{1}{qN^{1/2}},
\end{gather*}
so
\begin{gather}\label{17-2}
\int_0^1\left|S_N(\theta)\right|^2d\theta=\int_0^1
\left|S_N(\frac{a}{q}+\beta)\right|^2d(\frac{a}{q}+\beta)\le
\mathop{{\sum}^*}_{0\le a\le q\le N^{1/2}}\int\limits_{|\beta|\le\frac{1}{qN^{1/2}}}
\left|S_N(\frac{a}{q}+\beta)\right|^2d\beta.
\end{gather}
The change of variables $j=cx+n$ in \eqref{17-2} leads to the inequality \eqref{17-1}. This completes the proof of the lemma.
\end{proof}
\end{Le}
We recall that
\begin{gather*}
Q_0=\max\left\{\exp\left(\frac{10^5A^4}{\epsilon_0^2}\right),\exp(\epsilon_0^{-5})\right\}.
\end{gather*}
\begin{Le}\label{lemma-17-2}
The following inequality holds
\begin{gather}\notag
\int_0^1\left|S_N(\theta)\right|^2d\theta\le
2Q_0^2\frac{|\Omega_N|^2}{N}+
\frac{1}{N}\mathop{{\sum}^*}_{0\le a\le q\le N^{1/2}\atop q>Q_0}\int\limits_{\frac{Q_0}{q}\le|K|\le\frac{N^{1/2}}{q}}
\left|S_N(\frac{a}{q}+\frac{K}{N})\right|^2dK+\\
\frac{1}{N}\mathop{{\sum}^*}_{0\le a\le q\le Q_0}\int\limits_{\frac{Q_0}{q}\le|K|\le\frac{N^{1/2}}{q}}
\left|S_N(\frac{a}{q}+\frac{K}{N})\right|^2dK
+\frac{1}{N}\mathop{{\sum}^*}_{1\le a\le q\le N^{1/2}\atop q>Q_0}\int\limits_{|K|\le\frac{Q_0}{q}}
\left|S_N(\frac{a}{q}+\frac{K}{N})\right|^2dK.\label{17-4}
\end{gather}
\begin{proof}
To simplify we write $f(K)=\left|S_N(\frac{a}{q}+\frac{K}{N})\right|^2,$ then
\begin{gather}\notag
\mathop{{\sum}^*}_{0\le a\le q\le N^{1/2}}\int\limits_{|K|\le\frac{N^{1/2}}{q}}f(K)dK=
\mathop{{\sum}^*}_{0\le a\le q\le N^{1/2}\atop q>Q_0}\int\limits_{\frac{Q_0}{q}<|K|\le\frac{N^{1/2}}{q}}f(K)dK+
\mathop{{\sum}^*}_{0\le a\le q\le Q_0}\int\limits_{\frac{Q_0}{q}<|K|\le\frac{N^{1/2}}{q}}f(K)dK+\\
\mathop{{\sum}^*}_{1\le a\le q\le N^{1/2}\atop q>Q_0}\int\limits_{|K|\le\frac{Q_0}{q}}f(K)dK+
\mathop{{\sum}^*}_{1\le a\le q\le Q_0}\int\limits_{|K|\le\frac{Q_0}{q}}f(K)dK.\label{17-5}
\end{gather}
We estimate the fourth integral trivially
\begin{gather}\label{17-6}
\sum_{q\le Q_0}\mathop{{\sum}^*}_{0\le a\le q}\int\limits_{|K|\le\frac{Q_0}{q}}f(K)dK\le 2Q_0^2|\Omega_N|^2.
\end{gather}
Substituting \eqref{17-6} into \eqref{17-5} and using \eqref{17-1}, we obtain \eqref{17-4}. This completes the proof of the lemma.
\end{proof}
\end{Le}
First we estimate the third integral in the right side of \eqref{17-4}. It is convenient to use the following notation
\begin{gather}\label{17-7}
\gamma=1-\delta,\quad\xi_1=N^{2\gamma+6\epsilon_0}.
\end{gather}


\begin{Le}\label{lemma-17-3-0}
For $\gamma<\frac{1}{8}$ and $\epsilon_0\in(0,\frac{1}{2500})$  the following inequality holds
\begin{gather}\label{17-7-1}
\frac{1}{N}\mathop{{\sum}^*}_{1\le a\le q\le N^{1/2}\atop q>\xi_1}\int\limits_{|K|\le\frac{Q_0}{q}}
\left|S_N(\frac{a}{q}+\frac{K}{N})\right|^2dK\ll\frac{|\Omega_N|^2}{N}.
\end{gather}
\begin{proof}
We denote by $I$ the integral in the left side of \eqref{17-7-1}. Since ${\KK\ge1,\,q>Q_0,}$ the conditions of Lemma \ref{lemma-15-3} hold. Applying this lemma we obtain
\begin{gather}\label{17-7-2}
I\ll_{A,\epsilon_0}|\Omega_N|^2\frac{N^{2\gamma+2\epsilon_0}}{q^{2-4\epsilon_0}}
\int\limits_{|K|\le\frac{Q_0}{q}}\frac{dK}{\KK^{2-4\epsilon_0}}\ll
|\Omega_N|^2\frac{N^{2\gamma+2\epsilon_0}}{q^{2-4\epsilon_0}}\frac{Q_0}{q},
\end{gather}
since $|K|\le\frac{Q_0}{q}\le1$ and, hence, $\KK=1.$ Substituting \eqref{17-7-2} into \eqref{17-7-1}, we obtain
\begin{gather}\label{17-7-3}
\frac{1}{N}\mathop{{\sum}^*}_{1\le a\le q\le N^{1/2}\atop q>\xi_1}I\ll\frac{|\Omega_N|^2}{N}
\sum_{\xi_1<q\le N^{1/2}}\frac{N^{2\gamma+2\epsilon_0}}{q^{2-4\epsilon_0}}.
\end{gather}
By the choice of the parameter $\xi_1,$ we have
\begin{gather}\label{17-7-4}
\sum_{\xi_1<q\le N^{1/2}}\frac{N^{2\gamma+2\epsilon_0}}{q^{2-4\epsilon_0}}\le
\frac{N^{2\gamma+2\epsilon_0}}{\xi_1^{1-4\epsilon_0}}\le N^{-\epsilon_0}\ll1.
\end{gather}
Substituting \eqref{17-7-4} into \eqref{17-7-3}, we obtain \eqref{17-7-1}. This completes the proof of the lemma.
\end{proof}
\end{Le}

\begin{Le}\label{lemma-17-3-1}
For $\gamma<\frac{1}{8}$ and $\epsilon_0\in(0,\frac{1}{2500})$  the following inequality holds
\begin{gather}\label{17-7-5}
\frac{1}{N}\mathop{{\sum}^*}_{1\le a\le q\le \xi_1\atop q>Q_0}\int\limits_{|K|\le\frac{Q_0}{q}}
\left|S_N(\frac{a}{q}+\frac{K}{N})\right|^2dK\ll\frac{|\Omega_N|^2}{N}.
\end{gather}
\begin{proof}
We denote by $I$ the integral in the left side of \eqref{17-7-5}. Thus ${q>Q_0,}$ so one has $\KK=1$ and, hence, all conditions of Lemma \ref{lemma-15-2-1} hold. Applying this lemma we obtain
\begin{gather}\label{17-7-6}
I\ll|\Omega_N|^2q^{4\gamma-2+12\epsilon_0}
\int\limits_{|K|\le\frac{Q_0}{q}}\frac{dK}{(\KK)^{1-4\gamma-12\epsilon_0}}\ll
|\Omega_N|^2q^{4\gamma-3+12\epsilon_0},
\end{gather}
since $|K|\le\frac{Q_0}{q}\le1$ and, hence, $\KK=1.$ Substituting \eqref{17-7-6} into \eqref{17-7-5}, we have
\begin{gather}\label{17-7-7}
\frac{1}{N}\mathop{{\sum}^*}_{1\le a\le q\le\xi_1\atop q>Q_0}I\ll\frac{|\Omega_N|^2}{N}
\sum_{1\le q\le\xi_1}q^{4\gamma-2+12\epsilon_0}.
\end{gather}
By the choice of parameters $\delta,\epsilon_0$ we obtain \eqref{17-7-5}. This completes the proof of the lemma.
\end{proof}
\end{Le}
Thus we have estimated the third integral in the right side of \eqref{17-4}. Next we estimate the second one.
\begin{Le}\label{lemma-17-3-2}
Under the hypotheses of Lemma \ref{BK lemma 7.1} for $\gamma<\frac{1}{8}$ and $\epsilon_0\in(0,\frac{1}{2500})$  the following inequality holds
\begin{gather}\label{17-7-8}
\frac{1}{N}\mathop{{\sum}^*}_{0\le a\le q\le Q_0}\int\limits_{\frac{Q_0}{q}\le|K|\le\frac{N^{1/2}}{q}}
\left|S_N(\frac{a}{q}+\frac{K}{N})\right|^2dK\ll\frac{|\Omega_N|^2}{N}.
\end{gather}
\begin{proof}
We denote by $I$ the integral in the left side of  \eqref{17-7-8}.  Applying Lemma \ref{lemma-16-1-0}, we obtain
\begin{gather}\label{17-7-9}
I\ll|\Omega_N|^2q^{\frac{36}{5}\gamma-2+38\epsilon_0}
\int\limits_{\frac{Q_0}{q}\le|K|\le\frac{N^{1/2}}{q}}(\KK)^{\frac{36}{5}\gamma-2+38\epsilon_0}dK\ll
\frac{|\Omega_N|^2}{q},
\end{gather}
by the choice of the parameter $\gamma.$ Summing \eqref{17-7-9} over $0\le a\le q\le Q_0,$ we obtain \eqref{17-7-8}. This completes the proof of the lemma.
\end{proof}
\end{Le}

It remains to estimate the first integral in the right side of \eqref{17-4}, that is,
\begin{gather}\label{17-7-10}
\frac{1}{N}\mathop{{\sum}^*}_{0\le a\le q\le N^{1/2}\atop q>Q_0}\int\limits_{\frac{Q_0}{q}\le|K|\le\frac{N^{1/2}}{q}}
\left|S_N(\frac{a}{q}+\frac{K}{N})\right|^2dK.
\end{gather}
The following lemmas will be devoted to this. We partition the range of summation and integration over $q,\,K$ into nine subareas:


\begin{center}
  \includegraphics[width=450pt,height=350pt]{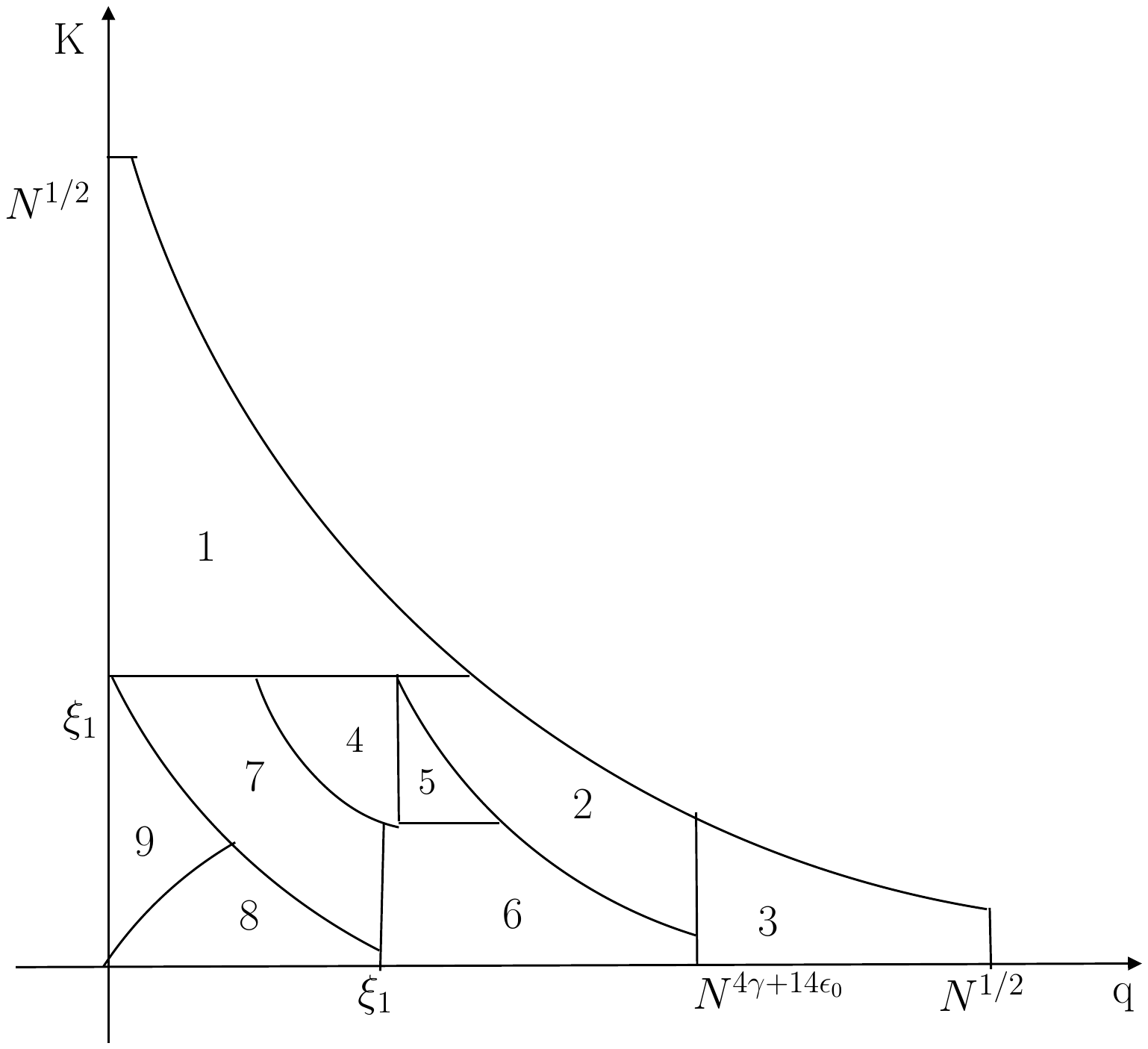}
\end{center}
Lemma \ref{lemma-17-3} corresponds to the domain 1, Lemma \ref{lemma-17-6} corresponds to the domain 2, Lemma \ref{lemma-17-7} corresponds to the domain 3, Lemma \ref{lemma-17-6-1} corresponds to the domain 4, Lemma \ref{lemma-17-6-2} corresponds to the domain 5, Lemma \ref{lemma-17-8} corresponds to the domain 6, Lemma \ref{lemma-17-9} corresponds to the domain 7, Lemma \ref{lemma-17-10} corresponds to the domain 8, Lemma \ref{lemma-17-11} corresponds to the domain 9.
\begin{Le}\label{lemma-17-3}
For $\gamma\le\frac{13-\sqrt{145}}{8}-5\epsilon_0,\,\epsilon_0\in(0,\frac{1}{2500})$ the following inequality holds
\begin{gather}\label{17-7-11}
\frac{1}{N}\mathop{{\sum}^*}_{1\le a\le q\le N^{1/2}\atop q>Q_0}\int\limits_{\xi_1\le|K|\le\frac{N^{1/2}}{q}}
\left|S_N(\frac{a}{q}+\frac{K}{N})\right|^2dK\ll\frac{|\Omega_N|^2}{N}.
\end{gather}
\begin{proof}
We denote by $I$ the integral in the left side of \eqref{17-7-11}.  Applying Lemma \ref{lemma-15-3}, we obtain
\begin{gather}\label{17-7-12}
I\ll|\Omega_N|^2\frac{N^{2\gamma+2\epsilon_0}}{q^{2-4\epsilon_0}}
\int\limits_{\xi_1\le|K|\le\frac{N^{1/2}}{q}}\frac{dK}{\KK^{2-4\epsilon_0}}\ll
|\Omega_N|^2\frac{N^{2\gamma+2\epsilon_0}}{q^{2-4\epsilon_0}\xi_1^{1-4\epsilon_0}}.
\end{gather}
Substituting \eqref{17-7-12} into the left side of \eqref{17-7-11}, we have
\begin{gather}\label{17-7-13}
\frac{1}{N}\mathop{{\sum}^*}_{1\le a\le q\le N^{1/2}\atop q>Q_0}\int\limits_{\xi_1\le|K|\le\frac{N^{1/2}}{q}}
\left|S_N(\frac{a}{q}+\frac{K}{N})\right|^2dK\ll
\frac{|\Omega_N|^2}{N}
\sum_{Q_0<q\le N^{1/2}}\frac{N^{2\gamma+2\epsilon_0}}{q^{1-4\epsilon_0}\xi_1^{1-4\epsilon_0}}.
\end{gather}
By the choice of the parameter $\xi_1,$ we obtain
\begin{gather}\label{17-7-14}
\sum_{Q_0<q\le N^{1/2}}\frac{N^{2\gamma+2\epsilon_0}}{q^{1-4\epsilon_0}\xi_1^{1-4\epsilon_0}}\ll
\frac{N^{2\gamma+4\epsilon_0}}{\xi_1^{1-4\epsilon_0}}\ll N^{-\epsilon_0/2}\ll1.
\end{gather}
Substituting \eqref{17-7-14} into \eqref{17-7-13}, we obtain \eqref{17-7-11}. This completes the proof of the lemma.
\end{proof}
\end{Le}


Let
\begin{gather*}
c_1=c_1(N),\,c_2=c_2(N),\,Q_0\le c_1<c_2\le  N^{1/2},
\end{gather*}
and let
\begin{gather*}
f_1=f_1(N,q),\,f_2=f_2(N,q),\,\frac{Q_0}{q}\le f_1<f_2\le  \frac{N^{1/2}}{q},\\
m_1=\min\{f_1(N,N_j),f_1(N,N_{j+1})\},\,m_2=\max\{f_2(N,N_j),f_2(N,N_{j+1})\}.
\end{gather*}
We recall that the sequence $\{N_j\}_{-J-1}^{J+1}$ was defined in \eqref{10-8}.
\begin{Le}\label{lemma-17-5}
If the functions $f_1(N,q),f_2(N,q)$ are monotonic for $q,$ then the following inequality holds
\begin{gather}\notag
\mathop{{\sum}^*}_{c_1\le q\le c_2\atop 1\le a\le q }\,\int\limits_{f_1\le|K|\le f_2}
\left|S_N(\frac{a}{q}+\frac{K}{N})\right|^2dK\le\\\le
\sum_{j:\,c_1^{1-\epsilon_0}\le N_j\le c_2\,}\int\limits_{m_1\le|K|\le m_2}
\mathop{{\sum}^*}_{N_j\le q\le N_{j+1}\atop 1\le a\le q }\left|S_N(\frac{a}{q}+\frac{K}{N})\right|^2dK.\label{17-7-15}
\end{gather}
\begin{proof}
The interval $[c_1,c_2]$ can be covered  by the intervals $[N_j,N_{j+1}],$ then
\begin{gather*}
\sum_{c_1\le q\le c_2}\le\sum_{j:\,c_1^{1-\epsilon_0}\le N_j\le c_2}\sum_{N_j\le q\le N_{j+1}}.
\end{gather*}
Interchanging the order of summation over $q$ and integration over $K,$ we obtain \eqref{17-7-15}. This completes the proof of the lemma.
\end{proof}
\end{Le}
To simplify we denote $Q_1=N_j,Q=N_{j+1}.$ Using the relation of Lemma \ref{lemma10-2}, we obtain
\begin{gather*}
\frac{Q}{Q_1}\le Q^{\epsilon_0}\le N^{\epsilon_0/2},\quad c_1^{1-\epsilon_0}\le Q_1\le c_2,\quad c_1\le Q\le c_2^{1+2\epsilon_0}.
\end{gather*}
Further, we will use this bounds without reference to them. We recall that ${\KK=\max\{1,|K|\}.}$ We note that we will always have $m_2\ge1,$ thus for $\eta<1$ we have
\begin{gather}\label{int-K-1}
\int\limits_{m_1\le|K|\le m_2}\frac{dK}{\KK^{\eta}}\ll m_2^{1-\eta}.
\end{gather}
For $\eta>1$ we always have
\begin{gather}\label{int-K-2}
\int\limits_{m_1\le|K|\le m_2}\frac{dK}{\KK^{\eta}}<\frac{1}{m_1^{\eta-1}}.
\end{gather}
However, if $m_1\le1,$ then for $\eta>1$ one has
\begin{gather}\label{int-K-3}
\int\limits_{m_1\le|K|\le m_2}\frac{dK}{\KK^{\eta}}\ll1.
\end{gather}
\begin{Le}\label{lemma-17-6}
For $\gamma\le\frac{13-\sqrt{145}}{8}-5\epsilon_0,\,\epsilon_0\in(0,\frac{1}{2500})$ the following inequality holds
\begin{gather}\label{17-8}
\frac{1}{N}\mathop{{\sum}^*}_{1\le a\le q\le N^{4\gamma+14\epsilon_0}\atop q>N^{2\gamma+9\epsilon_0}}
\int\limits_{\frac{N^{4\gamma+15\epsilon_0}}{q}\le|K|\le\frac{N^{1/2}}{q}}
\left|S_N(\frac{a}{q}+\frac{K}{N})\right|^2dK\ll\frac{|\Omega_N|^2}{N}.
\end{gather}
\begin{proof}
We use Lemma \ref{lemma-17-5} with
\begin{gather*}
c_1=N^{2\gamma+9\epsilon_0},\,c_2=N^{4\gamma+14\epsilon_0},
f_1=\frac{N^{4\gamma+15\epsilon_0}}{q},\,f_2=\frac{N^{1/2}}{q},\,
m_1=\frac{N^{4\gamma+15\epsilon_0}}{Q},\,m_2=\frac{N^{1/2}}{Q_1}.
\end{gather*}
We note that $KQ\le N^{1/2}\frac{Q}{Q_1}\le N^{1/2+\epsilon_0/2},$ thus, applying Corollary \ref{sled15-2}, we obtain
\begin{gather}\label{17-10}
\mathop{{\sum}^*}_{N_j\le q\le N_{j+1}\atop 1\le a\le q }\left|S_N(\frac{a}{q}+\frac{K}{N})\right|^2\ll
|\Omega_N|^2\left(
\frac{N^{2\gamma-1/4+5\epsilon_0}Q}{(\KK Q_1)^{1-2\epsilon_0}}+
\frac{N^{2\gamma+4\epsilon_0}Q^{1/2}}{\KK^{2-2\epsilon_0}Q_1^{1-2\epsilon_0}}\right).
\end{gather}
Using \eqref{int-K-1} while integrating the first summand over $K$ and using \eqref{int-K-2} while integrating the second one, we obtain
\begin{gather*}
\int\limits_{m_1\le|K|\le m_2}\mathop{{\sum}^*}_{N_j\le q\le N_{j+1}\atop 1\le a\le q }
\left|S_N(\frac{a}{q}+\frac{K}{N})\right|^2dK\ll
|\Omega_N|^2\Bigl(
N^{2\gamma-\frac{1}{4}+7\epsilon_0}\frac{Q}{Q_1}+
N^{2\gamma+4\epsilon_0+(4\gamma+15\epsilon_0)(-1+2\epsilon_0)}\frac{Q^{3/2}}{Q_1}\Bigr).
\end{gather*}
Thus $Q\le N^{4\gamma+15\epsilon_0},\,\frac{Q}{Q_1}\le N^{\epsilon_0/2}$ and $\gamma<\frac{1}{8}-4\epsilon_0,$ we obtain
\begin{gather}\label{17-11}
\int\limits_{m_1\le|K|\le m_2}\mathop{{\sum}^*}_{N_j\le q\le N_{j+1}\atop 1\le a\le q }
\left|S_N(\frac{a}{q}+\frac{K}{N})\right|^2dK\ll
|\Omega_N|^2 N^{-0,1\epsilon_0}.
\end{gather}
Since the number of summands in the sum over $j$ is less than $\log\log N,$ one has
\begin{gather*}
\sum_{j:\,c_1^{1-\epsilon_0}\le N_j\le c_2\,}\int\limits_{m_1\le|K|\le m_2}
\mathop{{\sum}^*}_{N_j\le q\le N_{j+1}\atop 1\le a\le q }\left|S_N(\frac{a}{q}+\frac{K}{N})\right|^2dK\ll
|\Omega_N|^2 .
\end{gather*}
This completes the proof of the lemma.
\end{proof}
\end{Le}
\begin{Le}\label{lemma-17-7}
For $\gamma\le\frac{13-\sqrt{145}}{8}-5\epsilon_0,\,\epsilon_0\in(0,\frac{1}{2500})$  the following inequality holds
\begin{gather}\label{17-12}
\frac{1}{N}\mathop{{\sum}^*}_{1\le a\le q\le N^{1/2}\atop q>N^{4\gamma+14\epsilon_0}}
\int\limits_{\frac{Q_0}{q}\le|K|\le\frac{N^{1/2}}{q}}
\left|S_N(\frac{a}{q}+\frac{K}{N})\right|^2dK\ll\frac{|\Omega_N|^2}{N}.
\end{gather}
\begin{proof}
We use Lemma \ref{lemma-17-5} with
\begin{gather*}
c_1=N^{4\gamma+14\epsilon_0},\,c_2=N^{1/2},\,
f_1=\frac{Q_0}{q},\,f_2=\frac{N^{1/2}}{q},\,
m_1=\frac{Q_0}{Q},\,m_2=\frac{N^{1/2}}{Q_1}.
\end{gather*}
Applying Corollary \ref{sled15-2} we obtain the bound \eqref{17-10}. Using \eqref{int-K-1} while integrating the first summand over $K$ and using \eqref{int-K-3} while integrating the second one, we obtain
\begin{gather*}
\int\limits_{m_1\le|K|\le m_2}\mathop{{\sum}^*}_{N_j\le q\le N_{j+1}\atop 1\le a\le q }
\left|S_N(\frac{a}{q}+\frac{K}{N})\right|^2dK\ll
|\Omega_N|^2\Bigl(
N^{2\gamma-1/4+7\epsilon_0}\frac{Q}{Q_1}+
N^{2\gamma+4\epsilon_0}\frac{Q}{Q_1^{1-2\epsilon_0}Q^{1/2}}\Bigr).
\end{gather*}
Thus $\frac{Q}{Q_1^{1-2\epsilon_0}}\le N^{1,5\epsilon_0},\,Q^{1/2}\ge N^{2\gamma+7\epsilon_0}$ and $\gamma<\frac{1}{8}-4\epsilon_0,$ we have
\begin{gather*}
\int\limits_{m_1\le|K|\le m_2}\mathop{{\sum}^*}_{N_j\le q\le N_{j+1}\atop 1\le a\le q }
\left|S_N(\frac{a}{q}+\frac{K}{N})\right|^2dK\ll
|\Omega_N|^2 N^{-\epsilon_0}.
\end{gather*}
Since the number of summands in the sum over $j$ is less than $\log\log N,$ the lemma is proved.
\end{proof}
\end{Le}
\begin{Le}\label{lemma-17-6-1}
For $\gamma\le\frac{13-\sqrt{145}}{8}-5\epsilon_0,\,\epsilon_0\in(0,\frac{1}{2500})$ the following inequality holds
\begin{gather}\label{17-8-1}
\frac{1}{N}\mathop{{\sum}^*}_{1\le a\le q\le N^{2\gamma+9\epsilon_0}\atop q>N^{\gamma+6\epsilon_0}}
\int\limits_{\frac{N^{3\gamma+12\epsilon_0}}{q}\le|K|\le \xi_1}
\left|S_N(\frac{a}{q}+\frac{K}{N})\right|^2dK\ll\frac{|\Omega_N|^2}{N}.
\end{gather}
\begin{proof}
We use Lemma \ref{lemma-17-5} with
\begin{gather*}
c_1=N^{\gamma+6\epsilon_0},\,c_2=N^{2\gamma+9\epsilon_0},
f_1=\frac{N^{3\gamma+12\epsilon_0}}{q},\,f_2=\xi_1,\,
m_1=\frac{N^{3\gamma+12\epsilon_0}}{Q},\,m_2=\xi_1.
\end{gather*}
Applying Corollary \ref{sled15-2} we obtain the bound \eqref{17-10}. Using \eqref{int-K-1} while integrating the first summand over $K$ and using \eqref{int-K-2} while integrating the second one, we obtain
\begin{gather*}
\int\limits_{m_1\le|K|\le m_2}\mathop{{\sum}^*}_{N_j\le q\le N_{j+1}\atop 1\le a\le q }
\left|S_N(\frac{a}{q}+\frac{K}{N})\right|^2dK\ll
|\Omega_N|^2\Bigl(
N^{2\gamma-1/4+7\epsilon_0}+
N^{2\gamma+5\epsilon_0+(3\gamma+12\epsilon_0)(-1+2\epsilon_0)}\frac{Q^{3/2}}{Q_1}\Bigr).
\end{gather*}
Thus $\frac{Q}{Q_1}\le Q^{\epsilon_0},\,Q\le N^{2\gamma+10\epsilon_0}$ and $\gamma<\frac{1}{8}-4\epsilon_0,$ we obtain
\begin{gather}\label{17-11-2}
\int\limits_{m_1\le|K|\le m_2}\mathop{{\sum}^*}_{N_j\le q\le N_{j+1}\atop 1\le a\le q }
\left|S_N(\frac{a}{q}+\frac{K}{N})\right|^2dK\ll
|\Omega_N|^2 N^{-0,1\epsilon_0}.
\end{gather}
Since the number of summands in the sum over $j$ is less than $\log\log N,$ the lemma is proved.
\end{proof}
\end{Le}
\begin{Le}\label{lemma-17-6-2}
For $\gamma\le\frac{13-\sqrt{145}}{8}-5\epsilon_0,\,\epsilon_0\in(0,\frac{1}{2500})$ the following inequality holds
\begin{gather}\label{17-8-2}
\frac{1}{N}\mathop{{\sum}^*}_{1\le a\le q\le N^{3\gamma+11\epsilon_0}\atop q>N^{2\gamma+9\epsilon_0}}
\int\limits_{N^{\gamma+4\epsilon_0}\le|K|\le\frac{N^{4\gamma+15\epsilon_0}}{q}}
\left|S_N(\frac{a}{q}+\frac{K}{N})\right|^2dK\ll\frac{|\Omega_N|^2}{N}.
\end{gather}
\begin{proof}
We use Lemma \ref{lemma-17-5} with
\begin{gather*}
c_1=N^{2\gamma+9\epsilon_0},\,c_2=N^{3\gamma+11\epsilon_0},
f_1=N^{\gamma+4\epsilon_0},\,f_2=\frac{N^{4\gamma+15\epsilon_0}}{q},\,
m_1=N^{\gamma+4\epsilon_0},\,m_2=\frac{N^{4\gamma+15\epsilon_0}}{Q_1}.
\end{gather*}
Applying Corollary \ref{sled15-2} we obtain the bound \eqref{17-10}. Using \eqref{int-K-1} while integrating the first summand over $K$ and using \eqref{int-K-2} while integrating the second one, we obtain
\begin{gather*}
\int\limits_{m_1\le|K|\le m_2}\mathop{{\sum}^*}_{N_j\le q\le N_{j+1}\atop 1\le a\le q }
\left|S_N(\frac{a}{q}+\frac{K}{N})\right|^2dK\ll
|\Omega_N|^2\Bigl(
N^{2\gamma-1/4+7\epsilon_0}+
N^{2\gamma+6,5\epsilon_0+(\gamma+4\epsilon_0)(-1+2\epsilon_0)}\frac{1}{Q^{1/2}}\Bigr).
\end{gather*}
Thus $Q\ge N^{\gamma+9\epsilon_0}$ and $\gamma<\frac{1}{8}-4\epsilon_0,$ we obtain
\begin{gather}\label{17-11-2}
\int\limits_{m_1\le|K|\le m_2}\mathop{{\sum}^*}_{N_j\le q\le N_{j+1}\atop 1\le a\le q }
\left|S_N(\frac{a}{q}+\frac{K}{N})\right|^2dK\ll
|\Omega_N|^2 N^{-0,1\epsilon_0}.
\end{gather}
Since the number of summands in the sum over $j$ is less than $\log\log N,$ the lemma is proved.
\end{proof}
\end{Le}
\begin{Le}\label{lemma-17-8}
For $\gamma\le\frac{13-\sqrt{145}}{8}-5\epsilon_0,\,\epsilon_0\in(0,\frac{1}{2500})$ the following inequality holds
\begin{gather}\label{17-13-1}
\frac{1}{N}\mathop{{\sum}^*}_{1\le a\le q\le N^{4\gamma+14\epsilon_0}\atop q>\xi_1}
\int\limits_{\frac{Q_0}{q}\le|K|\le\min\{N^{\gamma+6\epsilon_0},\frac{N^{4\gamma+15\epsilon_0}}{q}\}}
\left|S_N(\frac{a}{q}+\frac{K}{N})\right|^2dK\ll\frac{|\Omega_N|^2}{N}.
\end{gather}
\begin{proof}
We use Lemma \ref{lemma-17-5} with
\begin{gather*}
c_1=\xi_1,\,c_2=N^{4\gamma+14\epsilon_0},\,
f_1=\frac{Q_0}{q},\,f_2=\min\{N^{\gamma+6\epsilon_0},\frac{N^{4\gamma+15\epsilon_0}}{q}\},\,
m_1=\frac{Q_0}{Q},\,m_2=\min\{N^{\gamma+6\epsilon_0},\frac{N^{4\gamma+15\epsilon_0}}{Q_1}\}
\end{gather*}
and estimate the inner sum in the obtained relation. Applying Corollary  \ref{sled15-1} with $\alpha=4\gamma+15,5\epsilon_0,$ we obtain
\begin{gather}\label{17-14}
\mathop{{\sum}^*}_{N_j\le q\le N_{j+1}\atop 1\le a\le q }\left|S_N(\frac{a}{q}+\frac{K}{N})\right|^2\ll
|\Omega_N|^2\left(C_1^2Q^{\epsilon_0}+C_2\right),
\end{gather}
where
\begin{gather}\label{17-14-1}
C^2_1\le\frac{N^{2\gamma+4\epsilon_0}}{\KK^2 Q_1}+
N^{-4\gamma^2+7\gamma-1+25\epsilon_0}+
N^{-4\gamma^2+9\gamma-\frac{3}{2}+33\epsilon_0}Q,\quad
C_2=\frac{N^{4\gamma-\frac{1}{2}+12\epsilon_0}Q}{(\KK Q_1)^{1-2\epsilon_0}},
\end{gather}
Since the number of summands in the sum over $j$ is less than $\log\log N,$ to prove the lemma it is sufficient after taking an integral to obtain a quantity less than $N^{-0,1\epsilon_0}.$ We estimate separately the integral of both summands in \eqref{17-14}. We note that for $\eta\le1$ one has
\begin{gather}\label{17-18}
\int\limits_{m_1\le|K|\le m_2}Q_1^{\eta}dK<Q_1^{\eta}\min\{N^{\gamma+6\epsilon_0},\frac{N^{4\gamma+15\epsilon_0}}{Q_1}\}\le
Q_1^{\eta}N^{(\gamma+6\epsilon_0)(1-\eta)}\left(\frac{N^{4\gamma+15\epsilon_0}}{Q_1}\right)^{\eta}=
N^{(3\gamma+9\epsilon_0)\eta+\gamma+6\epsilon_0}.
\end{gather}
\begin{enumerate}
  \item The estimate of $C_1^2Q^{\epsilon_0}.$ Using \eqref{int-K-3} while integrating the first summand over $K$ and using  \eqref{17-18} while integrating the second and the third one, we obtain
\begin{gather}\notag
\int\limits_{m_1\le|K|\le m_2}C_1^2Q^{\epsilon_0}dK\ll
\frac{N^{2\gamma+4\epsilon_0}Q^{\epsilon_0}}{Q_1}+
N^{-4\gamma^2+8\gamma-1+31\epsilon_0}+\\+
N^{-4\gamma^2+9\gamma-\frac{3}{2}+33\epsilon_0}N^{4\gamma+16\epsilon_0}.\label{17-15}
\end{gather}
Since $Q>\xi_1,$ and substituting $\xi_1$ in the first summand we obtain
\begin{gather*}
\frac{N^{2\gamma+4\epsilon_0}Q^{\epsilon_0}}{Q_1}<\frac{N^{2\gamma+4\epsilon_0}Q^{2\epsilon_0}}{Q}<
\frac{N^{2\gamma+5\epsilon_0}}{N^{2\gamma+6\epsilon_0}}<N^{-0,1\epsilon_0}.
\end{gather*}
In view of the conditions on $\gamma,\epsilon_0$ the second and the third summands are less than $N^{-0,1\epsilon_0}.$
 \item The estimate of $C_2.$ Using \eqref{int-K-3} while integrating over $K$ we obtain
\begin{gather*}
\int\limits_{m_1\le|K|\le m_2}C_2dK\ll
N^{4\gamma-\frac{1}{2}+12\epsilon_0}\frac{Q}{Q_1}N^{\epsilon_0}\le N^{4\gamma-\frac{1}{2}+14\epsilon_0}<N^{-0,1\epsilon_0}.
\end{gather*}
The last inequality holds in view of the conditions on $\gamma,\epsilon_0.$
\end{enumerate}
This completes the proof of the lemma.
\end{proof}
\end{Le}


\begin{Le}\label{lemma-17-9}
For $\gamma\le\frac{13-\sqrt{145}}{8}-5\epsilon_0,\,\epsilon_0\in(0,\frac{1}{2500})$ the following inequality holds
\begin{gather}\label{17-17}
\frac{1}{N}\mathop{{\sum}^*}_{1\le a\le q\le \xi_1\atop q>Q_0}
\int\limits_{\frac{\xi_1}{q}\le|K|\le\min\{\xi_1,\frac{N^{3\gamma+12\epsilon_0}}{q}\}}
\left|S_N(\frac{a}{q}+\frac{K}{N})\right|^2dK\ll\frac{|\Omega_N|^2}{N}.
\end{gather}
\begin{proof}
We use Lemma \ref{lemma-17-5} with
\begin{gather*}
c_1=Q_0,\,c_2=\xi_1,\,
f_1=\frac{\xi_1}{q},\,f_2=\min\{\xi_1,\frac{N^{3\gamma+12\epsilon_0}}{q}\},\,
m_1=\frac{\xi_1}{Q},\,m_2=\min\{\xi_1,\frac{N^{3\gamma+12\epsilon_0}}{Q_1}\}.
\end{gather*}
and estimate the inner sum in the obtained relation. Applying Corollary \ref{sled15-1} with $\alpha=3\gamma+12,5\epsilon_0,$ we obtain
\begin{gather}\label{17-17-1}
\mathop{{\sum}^*}_{N_j\le q\le N_{j+1}\atop 1\le a\le q }\left|S_N(\frac{a}{q}+\frac{K}{N})\right|^2\ll
|\Omega_N|^2\left(C_1^2Q^{\epsilon_0}+C_2\right),
\end{gather}
where
\begin{gather}\label{17-17-2}
C^2_1\le\frac{N^{2\gamma+4\epsilon_0}}{\KK^2 Q_1}+
N^{-3\gamma^2+6\gamma-1+22\epsilon_0}+
N^{-3\gamma^2+7,5\gamma-1,5+28\epsilon_0}Q,\quad
C_2=\frac{N^{3,5\gamma-\frac{1}{2}+10\epsilon_0}Q}{(\KK Q_1)^{1-2\epsilon_0}}.
\end{gather}
As usual it is sufficient to to obtain a quantity less than $N^{-0,1\epsilon_0}$ after taking an integral over $K.$  The integral of $C_2$ can be estimated in the same way as in Lemma \ref{lemma-17-8}. Thus we have to estimate only the integral of $C_1^2Q^{\epsilon_0}.$ We note that for $\eta\le1$ one has
\begin{gather}\label{17-17-3}
\int\limits_{m_1\le|K|\le m_2}Q_1^{\eta}dK<Q_1^{\eta}\min\{\xi_1,\frac{N^{3\gamma+12\epsilon_0}}{Q_1}\}\le
Q_1^{\eta}\xi_1^{(1-\eta)}\left(\frac{N^{3\gamma+12\epsilon_0}}{Q_1}\right)^{\eta}=
N^{(\gamma+6\epsilon_0)\eta+2\gamma+6\epsilon_0}.
\end{gather}
We estimate separately the integral of each of the three summands in $C_1^2Q^{\epsilon_0}.$ Using \eqref{int-K-2} while integrating the first summand over $K$ and using \eqref{17-17-3} while integrating the second and the third one, we obtain
\begin{gather*}
\int\limits_{m_1\le|K|\le m_2}C_1^2Q^{\epsilon_0}dK\ll
\frac{N^{2\gamma+4\epsilon_0}Q^{\epsilon_0}}{Q_1}\frac{Q}{\xi_1}+
N^{-3\gamma^2+8\gamma-1+28\epsilon_0}+
N^{-3\gamma^2+10,5\gamma-1,5+41\epsilon_0}\xi_1Q.
\end{gather*}
The first summand is less than $N^{-0,1\epsilon_0}$ since the definition of $\xi_1.$ The second and the third summand are less than $N^{-0,1\epsilon_0}$ because of the conditions on $\gamma,\epsilon_0.$ This completes the proof of the lemma.
\end{proof}
\end{Le}
We denote
\begin{gather}\label{17-21}
\nu=\frac{\sqrt{369}-7}{20}=0,61\ldots,\quad Q_C=\xi_1^{\frac{1}{\nu+1}}.
\end{gather}

\begin{Le}\label{lemma-17-10}
For $\gamma\le\frac{5}{\sqrt{369}+23}-8\epsilon_0,\,\epsilon_0\in(0,\frac{1}{2500})$ the following inequality holds
\begin{gather}\label{17-22}
\frac{1}{N}\mathop{{\sum}^*}_{1\le a\le q\le \xi_1\atop q>Q_0}
\int\limits_{\frac{Q_0}{q}\le|K|\le\min\{q^{\nu},\frac{\xi_1}{q}\}}
\left|S_N(\frac{a}{q}+\frac{K}{N})\right|^2dK\ll\frac{|\Omega_N|^2}{N}.
\end{gather}
\begin{proof}
We partition the left side of \eqref{17-22} into two sums
\begin{gather}\label{17-23}
\mathop{{\sum}^*}_{1\le a\le q\le \xi_1\atop q>Q_0}\int\limits_{\frac{Q_0}{q}\le|K|\le\min\{q^{\nu},\frac{\xi_1}{q}\}}=
\mathop{{\sum}^*}_{1\le a\le q\le Q_C\atop q>Q_0}\int\limits_{\frac{Q_0}{q}\le|K|\le q^{\nu}}+
\mathop{{\sum}^*}_{1\le a\le q\le \xi_1\atop q>Q_C}\int\limits_{\frac{Q_0}{q}\le|K|\le\frac{\xi_1}{q}}
\end{gather}
For both of these sums we apply Lemma \ref{lemma-17-5} with
\begin{gather*}
c_1=Q_0,\,c_2=Q_C,\,
f_1=\frac{Q_0}{q},\,f_2=q^{\nu},\,
m_1=\frac{Q_0}{Q_1},\,m_2=Q^{\nu}
\end{gather*}
for the first sum and with
\begin{gather*}
c_1=Q_C,\,c_2=\xi_2,\,
f_1=\frac{Q_0}{q},\,f_2=\frac{\xi_1}{q},\,
m_1=\frac{Q_0}{Q_1},\,m_2=\frac{\xi_1}{Q_1}
\end{gather*}
for the second one. To estimate both sums we use Lemma \ref{lemma-15-2}:
\begin{gather}\label{17-22-1}
\mathop{{\sum}^*}_{N_j\le q\le N_{j+1}\atop 1\le a\le q }\left|S_N(\frac{a}{q}+\frac{K}{N})\right|^2\ll
|\Omega_N|^2\frac{\KK^{4\gamma+12\epsilon_0}Q^{6\gamma+1+20\epsilon_0}}{\KK Q^2_1}.
\end{gather}
The conditions of Lemma \ref{lemma-15-2} hold in both domains because of the choice of the parameter $\xi_1$ and the conditions on $\gamma,\epsilon_0:$
\begin{gather*}
\KK^2Q^3\le Q_C^{3+2\nu+5\epsilon_0}\le N^{1-1,5\epsilon_0},\quad
\KK^2Q^3\le \xi_1^2Q^{1+2\epsilon_0}\le \xi_1^{3+5\epsilon_0}\le N^{1-1,5\epsilon_0}.
\end{gather*}
Using \eqref{int-K-1} while integrating the first summand over $K$ we obtain
\begin{gather}\label{17-23}
\int\limits_{m_1\le|K|\le m_2}
\mathop{{\sum}^*}_{N_j\le q\le N_{j+1}\atop 1\le a\le q }\left|S_N(\frac{a}{q}+\frac{K}{N})\right|^2\ll
|\Omega_N|^2Q^{6\gamma-1+22\epsilon_0}Q^{\nu(4\gamma+12\epsilon_0)}.
\end{gather}
For the sum over $j$ to be bounded by a constant, it is sufficient to have
\begin{gather}\label{17-23-1}
\gamma\le\frac{1-(23+12\nu)\epsilon_0}{6+4\nu}\,\Rightarrow\,6\gamma-1+22\epsilon_0+\nu(4\gamma+12\epsilon_0)\le -0,1\epsilon_0.
\end{gather}
Substituting the value of $\nu,$ we obtain that it is sufficient to have $\gamma\le\frac{5}{\sqrt{369}+23}-5\epsilon_0.$
Using \eqref{int-K-1} while integrating over $K,$ we obtain for the second sum
\begin{gather*}
\int\limits_{m_1\le|K|\le m_2}
\mathop{{\sum}^*}_{N_j\le q\le N_{j+1}\atop 1\le a\le q }\left|S_N(\frac{a}{q}+\frac{K}{N})\right|^2\ll
|\Omega_N|^2Q^{6\gamma-1+22\epsilon_0}\frac{\xi_1^{4\gamma+12\epsilon_0}}{Q^{4\gamma+12\epsilon_0}}\le\\\le
Q_C^{2\gamma-1+10\epsilon_0}\xi_1^{4\gamma+12\epsilon_0}=\xi_1^{4\gamma+12\epsilon_0+(2\gamma-1+10\epsilon_0)/(\nu+1)}
\end{gather*}
For the sum over $j$ to be bounded by a constant, it is sufficient to have
\begin{gather*}
4\gamma+12\epsilon_0+(2\gamma-1+10\epsilon_0)/(\nu+1)<-0,1\epsilon_0,
\end{gather*}
which is equivalent to
\begin{gather*}
\gamma\le\frac{1-(23+12\nu)\epsilon_0}{6+4\nu}.
\end{gather*}
This inequality can be checked in the same manner as \eqref{17-23-1}. This completes the proof of the lemma.
\end{proof}
\end{Le}

\begin{Le}\label{lemma-17-11}
For $\gamma\le\frac{5}{\sqrt{369}+23}-8\epsilon_0,\,\epsilon_0\in(0,\frac{1}{2500})$  the following inequality holds
\begin{gather}\label{17-24}
\frac{1}{N}\mathop{{\sum}^*}_{1\le a\le q\le Q_C\atop q>Q_0}
\int\limits_{q^{\nu}\le|K|\le\frac{\xi_1}{q}}
\left|S_N(\frac{a}{q}+\frac{K}{N})\right|^2dK\ll\frac{|\Omega_N|^2}{N}.
\end{gather}
\begin{proof}
We use Lemma \ref{lemma-17-5} with
\begin{gather*}
c_1=Q_0,\,c_2=Q_C,\,
f_1=q^{\nu},\,f_2=\frac{\xi_1}{q},\,
m_1=Q_1^{\nu},\,m_2=\frac{\xi_1}{Q_1}
\end{gather*}
Let verify the conditions of Lemma \ref{lemma-16-1}:
\begin{gather*}
(\KK Q)^{18/5+21\epsilon_0}Q\le\left(\frac{\xi_1Q}{Q_1}\right)^{18/5+21\epsilon_0}Q\le
\xi_1^{18/5+21\epsilon_0}Q_C^{1+5\epsilon_0}\le
\xi_1^{\frac{18\nu+23}{5(\nu+1)}+24,5\epsilon_0}<N
\end{gather*}
Substituting $\xi_1,$ we obtain that the last inequality follows from
\begin{gather*}
\gamma\le\frac{5(1+\nu)}{46+36\nu}-5\epsilon_0,
\end{gather*}
which holds in view of the definition of $\nu$ and the conditions on $\gamma.$ Applying Lemma \ref{lemma-16-1}, we obtain
\begin{gather}\label{17-25}
\mathop{{\sum}^*}_{N_j\le q\le N_{j+1}\atop 1\le a\le q }\left|S_N(\frac{a}{q}+\frac{K}{N})\right|^2\ll
|\Omega_N|^2
\KK^{\frac{36}{5}\gamma+34\epsilon_0-2}Q^{\frac{46}{5}\gamma+48\epsilon_0-1}.
\end{gather}
Using \eqref{int-K-2} while integrating over $K,$ we obtain
\begin{gather*}
\int\limits_{m_1\le|K|\le m_2}
\mathop{{\sum}^*}_{N_j\le q\le N_{j+1}\atop 1\le a\le q }\left|S_N(\frac{a}{q}+\frac{K}{N})\right|^2dK\ll
|\Omega_N|^2 Q^{\frac{46}{5}\gamma+48\epsilon_0-1}
Q_1^{\nu(\frac{36}{5}\gamma+34\epsilon_0-1)}\le\\\le
|\Omega_N|^2 Q^{\frac{46}{5}\gamma-1+\nu(\frac{36}{5}\gamma-1)+72\epsilon_0}
\end{gather*}
For the sum over $j$ to be bounded by a constant, it is sufficient to have
\begin{gather*}
\gamma\le\frac{5(1+\nu)}{46+36\nu}-6\epsilon_0.
\end{gather*}
This inequality holds in view of the definition of $\nu$ and the conditions on $\gamma.$ This completes the proof of the lemma.
\end{proof}
\end{Le}

\section{The proof of Theorem \ref{uslov}.}
Let $\gamma<\frac{5}{\sqrt{369}+23}.$ We choose $\epsilon_0$ such that $\epsilon_0\in(0,\frac{1}{2500})$ and
$\gamma\le\frac{5}{\sqrt{369}+23}-8\epsilon_0.$ Then it follows from Lemma \ref{lemma-17-3}~--\ref{lemma-17-11} that the first integral in the right side of \eqref{17-4} is less than $\frac{|\Omega_N|^2}{N},$ that is,
\begin{gather}\label{17-30}
\frac{1}{N}\mathop{{\sum}^*}_{0\le a\le q\le N^{1/2}\atop q>Q_0}\int\limits_{\frac{Q_0}{q}\le|K|\le\frac{N^{1/2}}{q}}
\left|S_N(\frac{a}{q}+\frac{K}{N})\right|^2dK\ll\frac{|\Omega_N|^2}{N}.
\end{gather}.
Substituting \eqref{17-30} and the results of Lemma \ref{lemma-17-3-0}~--\ref{lemma-17-3-2} in Lemma \ref{lemma-17-2}, we obtain
\begin{gather}\label{17-31}
\int_0^1\left|S_N(\theta)\right|^2d\theta\ll\frac{|\Omega_N|^2}{N}\quad\mbox{for}\quad
\gamma<\frac{27-\sqrt{633}}{16}.
\end{gather}
Thus the inequality \eqref{8-7} is proved. This, as it was proved in the section §\ref{section idea BK}, is enough for proving Theorem \ref{uslov}. This completes the proof of the theorem.
\begin{Zam}
As mentioned in Remark \ref{zam1} the condition \eqref{BK lemma 7.1 condition} can be replaced by
\begin{equation}\label{BK lemma 7.1 condition-2}
(qK)^{\frac{64}{25}+4\epsilon_0}<Y<X.
\end{equation}
After the necessary changes of Lemma \ref{lemma-16-1} we obtain that the optimal value of the parameter $\nu$ is $\nu=\frac{\sqrt{2274}-18}{50}=0,59\ldots.$  Hence, the statements of Lemma \ref{lemma-17-10},\,\eqref{lemma-17-11} hold for $\gamma\le\frac{25}{114+2\sqrt{2274}}-10\epsilon_0,\,\epsilon_0\in(0,\frac{1}{2500}).$ Then it follows from Lemma \ref{lemma-17-3}~--\ref{lemma-17-11} that the statement of Theorem \ref{uslov} hold for
\begin{equation}\label{KFcondition1-2}
\delta_{\A}>1-\frac{25}{114+2\sqrt{2274}}=0,8805\ldots.
\end{equation}
\end{Zam}


\textit{D.A. Frolenkov}\\
\textit{Steklov Mathematical Institute}\\
\textit{Gubkina str., 8, }\\
\textit{119991, Moscow, Russia}\\
\textit{e-mail: frolenkov\underline{ }adv@mail.ru}\\
\\
\textit{I.D. Kan}\\
\textit{Department of Number theory}\\
\textit{Moscow State University}\\
\textit{e-mail: igor.kan@list.ru}

\end{document}